\documentclass[11pt,leqno]{article}
\usepackage{a4wide,amsmath}
\usepackage{amsfonts,amsbsy,amscd}
\usepackage{amssymb}
\usepackage{mathrsfs}
\usepackage{mathcomp}
\usepackage{color}
\usepackage{tikz}
\usepackage[colorlinks]{hyperref}
\newtheorem{theorem}{Theorem}[section]
\newtheorem{corollary}[theorem]{Corollary}
\newtheorem{lemma}[theorem]{Lemma}
\newtheorem{proposition}[theorem]{Proposition}

\newtheorem{definition}[theorem]{Definition}
\newtheorem{remark}[theorem]{Remark}

\numberwithin{equation}{section}

\newcommand{\RR}{\mathbb{R}}
\newcommand{\CC}{\mathbb{C}}

\newcommand{\HH}{\mathbb{H}}

\newcommand{\ee}{\mathrm{e}}
\newcommand{\ii}{\mathrm{i}}
\renewcommand{\colon}{:\,}

\newcommand{\eqdef}{\stackrel{{\rm {def}}}{=}}   % Composed symbols

\newcommand{\Sgn}{\mathop{\rm {sign}}}

\newcommand{\Square}{$\sqcap$\hskip -1.5ex $\sqcup$}
\newcommand{\Blacksquare}{\vrule height 1.7ex width 1.7ex depth 0.2ex }
\newcommand{\proof}{{\em Proof. }}
\newcommand{\qed}{$\;$\Blacksquare}

\definecolor{violet}{rgb}{0.5,0,0.5}
\definecolor{orange}{cmyk}{0,0.3,0.7,0}

\title{\null\vspace*{-2.0cm}
      The Heston stochastic volatility model\\
      has a boundary trace at zero volatility%
\footnote{A part of this research was performed while
          the second author (P.T.) was a visiting professor at
          Toulouse School of Economics, I.M.T.,
          Universit\'e de Toulouse -- Capitole,
          Toulouse, France.\vspace{2.0mm}}
\vspace*{0.5cm}
}
 
\author{%
        B\'en\'edicte Alziary%
\thanks{{\it e-mail:} {\tt benedicte.alziary@ut-capitole.fr},$\;$}
\vspace*{0.3cm}
\\
        Toulouse School of Economics, I.M.T.,
        Universit\'e de Toulouse -- Capitole \\
        1, Esplanade de l'Universit\'e,
        F--31000 Toulouse Cedex, France \\
\vspace*{0.3cm}
\\
        and \\
\vspace*{0.3cm}
\and
        Peter Tak\'a\v{c}%
\thanks{{\it e-mail:} {\tt peter.takac@uni-rostock.de}.$\;$}
\vspace*{0.3cm}
\\
        Universit\"at Rostock,
        Institut f\"ur Mathematik \\
        Ulmenstra{\ss}e~69, Haus~3,
        D--18057 Rostock, Germany
\vspace*{0.5cm}
}
\date\today

\begin{document}
\baselineskip=16pt plus 1pt minus 1pt

\maketitle
%\vspace{3cm}
\baselineskip=14pt plus 1pt minus 1pt
%\tableofcontents
\baselineskip=16pt plus 1pt minus 1pt
 
%\newpage

\begin{abstract}
We establish boundary regularity results in H\"older spaces
for the degenerate parabolic problem obtained from
the Heston stochastic volatility model in Mathematical Finance set up
in the spatial domain (upper half\--plane)
$\mathbb{H} = \RR\times (0,\infty)\subset \RR^2$.
Starting with nonsmooth initial data $u_0\in H$, we take advantage of
smoothing properties of the parabolic semigroup
$\ee^{-t\mathcal{A}}\colon H\to H$, $t\in \RR_+$,
generated by the Heston model, to derive the smoothness of the solution
$u(t) = \ee^{-t\mathcal{A}} u_0$ for all $t>0$.
The existence and uniqueness of a weak solution is obtained
in a Hilbert space $H = L^2(\mathbb{H};\mathfrak{w})$
with very weak growth restrictions at infinity and on the boundary
$\partial\HH = \RR\times \{ 0\}\subset \RR^2$
of the half\--plane $\mathbb{H}$.
We investigate the influence of the boundary behavior
of the initial data $u_0\in H$ on the boundary behavior of
$u(t)$ for $t>0$.
\end{abstract}

\par\vfill
\vspace*{0.5cm}
\noindent
\begin{tabular}{lll}
{\bf 2020 Mathematics Subject Classification:}
& Primary   & 35B65, 35K65;\\
& Secondary & 35K15, 91G80;\\
\end{tabular}
 
\par\vspace*{0.5cm}
\noindent
\begin{tabular}{ll}
{\bf Key words:}
& degenerate parabolic equation; weighted Sobolev space; \\
& holomorphic semigroup; parabolic smoothing effect; \\
& dynamic boundary conditions; Heston's stochastic volatility model \\
\end{tabular}
 
%%%%%%%%%%%%%%%%%%%%%%%%%%%%%%%%%%%%%%%%%%%%%%%%%%%%%%%%%%%%%%%%%%%%%%%%
%\baselineskip=14pt plus 1pt minus 1pt
\baselineskip=16pt plus 1pt minus 1pt
\parskip=2mm plus .5mm minus .5mm
%%%%%%%%%%%%%%%%%%%%%%%%%%%%%%%%%%%%%%%%%%%%%%%%%%%%%%%%%%%%%%%%%%%%%%%%
\newpage
 
%%%%%%%%%%%%%%%%%%%%%%%%%%%%%%%%%%%%%%%%%%%%%%%%%%%%%%%%%%%%%%%%%%%%%%%
%%%%%    INTRODUCTION    %%%%%%%%%%%%%%%%%%%%%%%%%%%%%%%%%%%%%%%%%%%%%%
%%%%%%%%%%%%%%%%%%%%%%%%%%%%%%%%%%%%%%%%%%%%%%%%%%%%%%%%%%%%%%%%%%%%%%%

\section{Introduction}
\label{s:Intro}

The {\it\bfseries Heston stochastic volatility model\/} for
pricing the {\it European call options on stocks\/}
({\sc S.~L.\ Heston} \cite{Heston})
has been of considerable interest to economists and mathematicians
for almost three decades.
Numerous articles have been written about mathematical treatment
and solvability of this model in a number of settings.
In our present work we focus on the degenerate parabolic problem
with prescribed initial and boundary conditions.
The question of existence, uniqueness, and regularity of
a weak solution to this problem is studied in
{\sc P.~M.~N.\ Feehan} and {\sc C.~A.\ Pop} \cite{Feehan-Pop-13},
{\sc C.\ Chiarella}, {\sc B.\ Kang}, and {\sc G.~H. Meyer}
\cite{Chiarella-K-M-2015},
{\sc G.~H. Meyer} \cite{Meyer-2015}, and
{\sc B.\ Alziary} and {\sc P.\ Tak\'a\v{c}}
\cite[Sect.~4, pp.\ 16--17]{AlziaryTak},
to mention only a few.
The analyticity of the solution in both, space and time variables,
has been established in
\cite[Sect.~4, Theorem~4.2, pp.\ 16--17]{AlziaryTak}.
As a consequence, the completeness of the market
(cf.\ {\sc T.\ Bj\"ork} \cite[Sect.~8, pp.\ 115--124]{Bjoerk-3rd}
 and
 {\sc M.~H.~A.\ Davis} and {\sc J.\ Ob{\l}{\'o}j}
 \cite{Davis-Obloj})
described in Heston's model is verified in
\cite[Sect.~5, Theorem~5.2, p.~19]{AlziaryTak}.
Thanks to the importance of Heston's model in Mathematical Finance,
there is a strong interest in efficient numerical methods
applicable to computing the solution of this degenerate parabolic problem
(\cite{Chiarella-K-M-2015, Meyer-2015}).
A major obstacle to an efficient numerical method is
the degeneracy of the diffusion coefficient at low volatility; see e.g.\
{\sc B.\ D{\"u}ring} and {\sc M.~Fourni{\'e}} \cite{Duering-Fournie}
and
{\sc S.\ Ikonen} and {\sc J.\ Toivanen} \cite{Ikonen-Toivanen}.
This degeneracy causes serious problems in formulating and justifying
the correct boundary conditions on the portion of the boundary
with vanishing volatility, denoted by $\partial\mathbb{H}$.
A numerical scheme using a finite difference method in the domain
$\HH = \RR\times (0,\infty)\subset \RR^2$ with the boundary
$\partial\HH = \RR\times \{ 0\}$
has to be designed with a mesh of points much too fine
near the boundary $\partial\HH$, so that it finally becomes
rather inefficient and unprecise there.
This is one of the reasons why in this article we investigate
the limiting boundary behavior of the solution of Heston's model
as the volatility approaches zero.
We obtain a limiting partial differential equation of first order
on the boundary $\partial\HH$, Eq.~\eqref{e:Cauchy:xi=0},
thus specifying also the boundary conditions on $\partial\HH$.
It is worth of noticing that this equation on the space\--time domain
$\partial\HH\times (0,\infty)$
is coupled with the degenerate parabolic equation \eqref{e:Cauchy}
inside the domain $\HH\times (0,\infty)$
solely through a linear term with the partial derivative
with respect to the volatility
(the volatility approaching zero)
that appears in Eq.~\eqref{sol:Cauchy:xi=0}.
This feature of Heston's model is used in the recent work by
{\sc F.\ Baustian}, {\sc K.\ Filipov{\'a}}, and {\sc J.\ Posp\'{\i}\v{s}il}
\cite{Baust-Pospis-19}
with an orthogonal polynomial expansion in the spatial domain~$\HH$.
Orthogonal polynomial expansions have been used recently also in
{\sc D.\ Ackerer} and {\sc D.\ Filipovi{\'c}} \cite{Ackerer-Filip}
for numerical approximations.
Earlier, the authors
\cite[Sect.~11, pp.\ 48--51]{AlziaryTak}
have used orthogonal polynomial expansions with
Hermite and Laguerre polynomials in Gal\"erkin's method
to approximate functions in $L^2(\HH)$ by analytic functions.

Our derivation of Eq.~\eqref{e:Cauchy:xi=0} on
$\partial\HH\times (0,\infty)$ is motivated by
the limiting behavior of the diffusion part
(second\--order partial derivatives)
in Eq.~\eqref{e:Cauchy}.
The limit, equal to zero on $\partial\HH\times (0,\infty)$,
has been obtained in
{\sc P.~M.~N.\ Feehan} and {\sc C.~A.\ Pop}
\cite{Feehan-Pop-13}, Lemma 3.1, Eq.\ (3.1), on p.~4409
(see also
 {\sc P.\ Daskalopoulos} and {\sc R.\ Hamilton}
 \cite{Daska-Hamil-98}, Prop.\ I.12.1 on p.~940)
for the corresponding (stationary) elliptic problem
with the Heston operator $\mathcal{A}$ given by
Eq.~\eqref{e:Heston-oper}.
However, in order to fulfill the regularity hypothesis required in
\cite[Lemma 3.1]{Feehan-Pop-13},
we need to establish a new regularity result for the weak solution
$u(\,\cdot\,, \,\cdot\,, t)\colon [0,\infty)]\to H$
of the Heston model
(see Proposition~\ref{prop-Lions} with $f\equiv 0$)
which is given by the $C^0$-semigroup of bounded linear operators
$\ee^{-t\mathcal{A}}\colon H\to H$, $t\in \RR_+$,
determined by the homogeneous initial value problem \eqref{e:Cauchy},
that is to say,
\begin{math}
  u(\,\cdot\,,\,\cdot\,,t)\equiv u(t) = \ee^{-t\mathcal{A}} u_0\in H ,\
  t\in \RR_+ ,
\end{math}
with an arbitrary initial value $u_0\in H$.
The underlying Hilbert space $H$
is a weighted $L^2$\--type Lebesgue space
$H = L^2(\mathbb{H};\mathfrak{w})$.
Our regularity result is based on the smoothing property of
the holomorphic semigroup $\ee^{-t\mathcal{A}}$, $t\in \RR_+$,
acting on $H$, see Theorem~\ref{thm-Regular}.
This result contains a number of local and global
partial regularity results which are new, as well.
We stress the main difference between the classical H{\"o}lder\--type
regularity treated in
{\sc P.~M.~N.\ Feehan} and {\sc C.~A.\ Pop}
\cite[Theorem~1.1 on p.~4409]{Feehan-Pop-13}
and the regularity obtained by parabolic smoothing:
The H{\"o}lder\--type regularity in \cite{Feehan-Pop-13}
assumes the same spatial regularity
already for the initial value $u(0) = u_0$
(in a suitable weighted H{\"o}lder space).
As a consequence, analogous regularity for the solution
$u(t)$ is proved (by Schauder estimates)
at all times $t\in (0,T)$ in a bounded time interval.
In contrast, we begin with {\em nonsmooth\/} initial data
$u_0\in H$ at $t=0$;
then we apply the {\em parabolic smoothing\/} of the $C^0$-semigroup
$\ee^{-t\mathcal{A}}$ for $t\in (0,\infty)$, thus arriving at
$u(t)\in \mathcal{D}(\mathcal{A}^k)\subset H$
for all $t\in (0,\infty)$ and every $k=1,2,3,\dots$.
Since the domain $\mathcal{D}(\mathcal{A}^k)$ of the $k$-th power of
the Heston operator $\mathcal{A}$ is the image (range) of
the $k$-th power of the bounded inverse
$(\lambda I + \mathcal{A})^{-1}\colon H\to H$
(the resolvent of $-\mathcal{A}$),
the solution $u(t) = \ee^{-t\mathcal{A}} u_0$
has higher smoothness for all $t > 0$.
This smoothing effect is essential for applications
in Mathematical Finance where the initial data $u_0\in H$
are typically not continuously differentiable
\begin{math}
(
  u_0\in W^{1,\infty}(\HH) \setminus C^1(\HH)
).
\end{math}
Indeed, for our derivation of the limiting equation
\eqref{e:Cauchy:xi=0} on
$\partial\HH\times (0,\infty)$ from equation \eqref{e:Cauchy}
we need H{\"o}lder regularity of type $C^{2+\alpha}$
over the closure of the open half\--plane $\HH$
(cf.\ \cite[Lemma 3.1]{Feehan-Pop-13}).

The proof of our main result, Theorem~\ref{thm-Regular},
makes use of the factorization\hfil\break
\begin{math}
  (\lambda I + \mathcal{A})^{-j}
  (\lambda I + \mathcal{A})^k \ee^{-t\mathcal{A}}
\end{math}
of the bounded linear operator
$\ee^{-t\mathcal{A}}\colon H\to H$ for $t > 0$;
with $k=1,2,3$ and $j=0,1,\dots,k$.
Thanks to the smoothing effect, the latter factor,
\begin{math}
  (\lambda I + \mathcal{A})^k \ee^{-t\mathcal{A}} ,
\end{math}
is a bounded linear operator on $H$ for each $t > 0$,
whereas the former factor,
$(\lambda I + \mathcal{A})^{-j}$,
is a bounded linear operator from $H$ to the domain
$\mathcal{D}(\mathcal{A}^j)$ of the $j$-th power of
the Heston operator $\mathcal{A}$.
We use the resolvent
$(\lambda I + \mathcal{A})^{-1}\colon H\to H$
of $-\mathcal{A}$ in order to describe the function space
$\mathcal{D}(\mathcal{A}^j)$
(endowed with the graph norm)
for $1\leq j\leq k\leq 3$.
This factorization
(in Section~\ref{s:smooth_Hoelder})
is split into three consecutive steps in Paragraphs
{\S}\ref{ss:smooth_Hoelder-1}, {\S}\ref{ss:smooth_Hoelder-2}, and
{\S}\ref{ss:smooth_Hoelder-3}, with the auxiliary functions
$f_{j,k}(\,\cdot\,,\,\cdot\,,t)\equiv f_{j,k}(t)\in H$
for $0\leq j\leq k\leq 3$ defined in eqs.\
\eqref{def:f_0,k(t)} and \eqref{def:f_j,k(t)}
for a given $u_0\in H$ and $t > 0$.
Clearly, for $j=k$; $k=1,2,3$, and $t > 0$ we obtain
$f_{k,k}(t) = u(t)$ which yields the desired regularity of
the solution $u(t)$ for $t > 0$ as stated in
Theorem~\ref{thm-Regular}.

Our second theorem
(Theorem~\ref{thm-MaxPrinciple})
is a {\it weak maximum principle\/}
for the initial value Cauchy problem \eqref{e:Cauchy}
in the unbounded space\--time domain $\HH\times (0,T)$.
As it is typical for parabolic problems posed in
an unbounded spatial domain
(the open half\--plane $\HH\subset \RR^2$ in our case),
the growth of the solution $u(t)\equiv u(x,\xi,t)$
has to be limited with respect to the space variable
$(x,\xi)\in \HH$ as $x\to \pm\infty$ or
$\xi\to 0+$ or $\xi\to +\infty$, uniformly for all $t\in (0,T)$.
We find a positive ``majorizing'' function
$\mathfrak{h}_0\colon \HH\to (0,\infty)$
in Eq.~\eqref{def:h_0} that provides the required limit on
the solution $u(x,\xi,t)$ in Theorem~\ref{thm-MaxPrinciple}.
This theorem has an important corollary applicable to
a typical initial value problem in Mathematical Finance
(see Corollary~\ref{cor-MaxPrinciple}).
The majorizing function
\begin{math}
%\label{ineq:Cauchy:u_0}
  U(x,\xi)\eqdef K_1\, \ee^{x + \varpi\xi} + K_0 ,
\end{math}
for all $(x,\xi)\in \HH$, provides an important upper bound
(independent from time)
as the {\it volatility\/} $\xi\in (0,\infty)$ approaches zero
($\xi\to 0+$).
Here, $K_0, K_1\in \RR_+$ are arbitrary constants, and\/
$\varpi\in \RR$ is another constant restricted by
inequalities in \eqref{ineq:varpi}.
This choice of the majorizing function and the initial data
covers the most typical alternatives for {\it derivative contracts\/}
(which are determined by the choice of the initial data $u_0$);
see e.g.\
{\sc J.-P.\ Fouque}, {\sc G.\ Papanicolaou}, and {\sc K.~R.\ Sircar}
\cite[{\S}1.2, pp.\ 8--12]{FouqPapaSir}.
The case of $u_0(x,\xi)\equiv u_0(x)$
being independent from the volatility $\xi\in (0,\infty)$
is of special interest (e.g., European call and put options);
we may set $\varpi = 0$.
Derivative contracts do not seem to include the volatility level
since volatility does not produce any direct returns such as
dividends or interest.
Volatility does not show long term upwards trends like equities,
but typically shows periods of high volatility occurring within
a short period of time (i.e. volatility ``jumps'') and then shows
a downward trend to return to the long run medium level.

This article is organized as follows.
We begin with basic notations and function spaces of
H{\"o}lder, Lebesgue, and Sobolev types, which involve weights.
Most of these spaces were originally introduced in
{\sc P.\ Daskalopoulos} and {\sc P.~M.~N.\ Feehan} \cite{Daska-Feehan-14}
and \cite[Sect.~2, p.~5048]{Daska-Feehan-16} and
{\sc P.~M.~N.\ Feehan} and {\sc C.~A.\ Pop} \cite{Feehan-Pop-15}.
The mathematical problem resulting from
{\sc S.~L.\ Heston}'s \cite{Heston} model in Mathematical Finance
(described in Appendix~\ref{s:econom-Heston} in ``economic'' terms)
is formulated in Section~\ref{s:MathProblem}.
The details of this formulation, especially a justification of
the boundary conditions and restrictions imposed on some important constants
(e.g., the volatility $\sigma > 0$ of the volatility,
 the rate of mean reversion $\kappa > 0$, and
 the long\--term variance $\theta > 0$),
such as the well\--known {\it Feller condition\/},
can be found in our previous work
\cite[Sect.~2, pp.\ 6--13]{AlziaryTak}.
Our main results are collected in Section~\ref{s:Main}, in
Theorems \ref{thm-Regular} and~\ref{thm-MaxPrinciple}.
In addition, also Proposition~\ref{prop-Lions}
(existence and uniqueness),
Corollary~\ref{cor-Regular} (boundary behavior), and
Corollary~\ref{cor-MaxPrinciple} (maximum principle)
are of importance.
Our strategy of the proof of Theorem~\ref{thm-Regular}
(laid out above) is described in all details in
Section~\ref{s:smooth_Heston}.
The first part of this strategy, obtaining H{\"o}lder regularity,
is implemented in Section~\ref{s:smooth_Hoelder}.
The proof of Theorem~\ref{thm-Regular}
(and that of Corollary~\ref{cor-Regular}, as well)
is completed in Section~\ref{s:proofMain}.
The main part of this article ends up with the proofs of
Theorem~\ref{thm-MaxPrinciple} and Corollary~\ref{cor-MaxPrinciple}
in Section~\ref{s:MaxPrinciple}.
We have postponed some rather technical results about
weighted Sobolev spaces and boundary traces until
Appendix~\ref{s:Sobolev_trace}.
Most of our regularity results gradually derived in
Section~\ref{s:smooth_Heston}
take advantage of difficult elliptic Schauder\--type estimates
for the degenerate Heston operator $\mathcal{A}$
in weighted H{\"o}lder spaces over the half\--plane $\HH$
obtained in a series of articles by
{\sc P.~M.~N.\ Feehan} and {\sc C.~A.\ Pop}
\cite{Feehan-Pop-14, Feehan-Pop-15, Feehan-Pop-17}.
For reader's convenience, we restate these results
in Appendix~\ref{s:ellipt_regul}.

%%%%%%%%%%%%%%%%%%%%%%%%%%%%%%%%%%%%%%%%%%%%%%%%%%%%%%%%%%%%%%%%%%%%%%%
%%%%%    Basic notation, function spaces    %%%%%%%%%%%%%%%%%%%%%%%%%%%
%%%%%%%%%%%%%%%%%%%%%%%%%%%%%%%%%%%%%%%%%%%%%%%%%%%%%%%%%%%%%%%%%%%%%%%

\section{Basic notations, function spaces}
\label{s:Notation}

We use the standard notation
$\mathbb{R} = (-\infty,+\infty)$, $\RR_+ = [0,\infty)\subset \RR$ and
$\HH = \RR\times (0,\infty)\subset \RR^2$ with the closure
$\overline{\HH} = \RR\times \RR_+$
for the open and closed upper half\--planes, respectively.
As usual, for $x\in \RR$ we abbreviate
$x^{+}\eqdef \max\{ x,\,0\}$ and $x^{-}\eqdef \max\{ -x,\,0\}$.
The complex plane is denoted by $\mathbb{C} = \RR + \ii\RR$.
The complex conjugate of a number $z\in \CC$ is denoted by $\bar{z}$,
so that the absolute value of $z$ is given by
$|z| = (z\bar{z})^{1/2}$.

The basic function space, $H$, in our treatment of the Heston model
is defined as follows:
We define the weight
$\mathfrak{w}\colon \mathbb{H}\to (0,\infty)$ by
\begin{equation}
\label{def:w}
  \mathfrak{w}(x,\xi)\eqdef
  \xi^{\beta - 1}\, \ee^{ - \gamma |x| - \mu\xi }
  \quad\mbox{ for $(x,\xi)\in \HH$, }
\end{equation}
where $\beta, \gamma, \mu\in (0,\infty)$
are suitable positive constants that will be specified later,
in Section~\ref{s:MathProblem}
(see also Appendix~\ref{s:Sobolev_trace}).
However, it is already clear that if we want that the weight
$\mathfrak{w}(x,\xi)$ tends to zero as $\xi\to 0+$,
we have to assume $\beta > 1$.
Similarly, if we want that the function
$u_0(x,\xi) = K (\mathrm{e}^x - 1)^{+}$ of $(x,\xi)\in \HH$
(an initial condition in Heston's model)
belongs to $H$, we must require $\gamma > 2$.
Then $H = L^2(\mathbb{H};\mathfrak{w})$
is the {\em complex\/} Hilbert space of
all complex\--valued Lebesgue\--measurable functions
$f\colon \HH\to \CC$ with the finite norm
\begin{equation*}
%\label{e:norm_H}
  \| f\|_H\eqdef
  \left( \int_{\HH} |f(x,\xi)|^2\, \mathfrak{w}(x,\xi)
         \,\mathrm{d}x \,\mathrm{d}\xi
  \right)^{1/2} < \infty \,.
\end{equation*}
This norm is induced by the inner product
\begin{equation*}
%\label{def:inn_prod}
  (f,g)_H\equiv
  (f,g)_{ L^2(\mathbb{H};\mathfrak{w}) } \eqdef
  \int_{\mathbb{H}} f\, \bar{g}\cdot \mathfrak{w}(x,\xi)
    \,\mathrm{d}x \,\mathrm{d}\xi
    \quad\mbox{ for }\, f,g\in H \,.
\end{equation*}

The domain, $V$, of the sesquilinear form that defines the Heston operator
is the weighted Sobolev space of all functions $f\in H$, such that
the first\--order partial derivatives
(in the sense of distributions),
\begin{math}
  f_x    \equiv \frac{\partial f}{\partial x} \,,\
  f_{\xi}\equiv \frac{\partial f}{\partial\xi} \,,
\end{math}
satisfy
\begin{equation*}
  [f]_V^2\eqdef
    \int_{\HH} \left( |f_x|^2 + |f_{\xi}|^2\right)
               \cdot \xi\cdot \mathfrak{w}(x,\xi)
               \,\mathrm{d}x \,\mathrm{d}\xi < \infty \,.
\end{equation*}
The Hilbert norm $\|\cdot\|_V$ on
$V = H^1(\mathbb{H};\mathfrak{w})$,
$\|f\|_V^2 = \|f\|_H^2 + [f]_V^2$ for $f\in V$,
is induced by the inner product
\begin{equation*}
%\label{def:inn_prod-H^1}
\begin{aligned}
  (f,g)_V\equiv
  (f,g)_{ H^1(\mathbb{H};\mathfrak{w}) }
& \eqdef \int_{\mathbb{H}}
  \left( f_x\, \bar{g}_x + f_{\xi}\, \bar{g}_{\xi} \right)
  \cdot \xi\cdot \mathfrak{w}(x,\xi) \,\mathrm{d}x \,\mathrm{d}\xi
  + \int_{\mathbb{H}} f\, \bar{g} \cdot \mathfrak{w}(x,\xi)
    \,\mathrm{d}x \,\mathrm{d}\xi
\end{aligned}
\end{equation*}
for $f,g\in H^1(\mathbb{H};\mathfrak{w})$.
In particular, the Sobolev imbedding $V\hookrightarrow H$
is bounded (i.e., continuous).

We will see later that the domain of the Heston operator is contained
in a local version of the following weighted Sobolev space,
$H^2(\mathbb{H};\mathfrak{w})$, of all functions $f\in V$,
such that also the second\--order partial derivatives
(in the sense of distributions),
\begin{math}
  f_{xx}    \equiv \frac{\partial^2 f}{\partial x^2} \,,
  f_{x\xi}  \equiv \frac{\partial^2 f}{\partial x\, \partial\xi} \,,
  f_{\xi\xi}\equiv \frac{\partial^2 f}{\partial\xi^2} \,,
\end{math}
satisfy
\begin{equation*}
  \int_{\HH} \left( |f_{xx}|^2 + |f_{x\xi}|^2 + |f_{\xi\xi}|^2\right)
             \cdot \xi^2\cdot \mathfrak{w}(x,\xi)
             \,\mathrm{d}x \,\mathrm{d}\xi < \infty \,.
\end{equation*}
In addition, we require that the Hilbert norm of $f$ on
$H^2(\mathbb{H};\mathfrak{w})$, as defined below, is finite,
\begin{equation*}
%\label{def:norm-H^2}
\begin{aligned}
  \| f\|_{ H^2(\mathbb{H};\mathfrak{w}) }^2
& \eqdef \int_{\mathbb{H}}
    \left( |f_{xx}|^2 + |f_{x\xi}|^2 + |f_{\xi\xi}|^2\right)
    \cdot \xi^2\cdot \mathfrak{w}(x,\xi)
    \,\mathrm{d}x \,\mathrm{d}\xi
\\
& {}
  + \int_{\mathbb{H}}
    \left( |f_{x}|^2 + |f_{\xi}|^2\right)
    \cdot (1 + \xi^2)\cdot \mathfrak{w}(x,\xi)
    \,\mathrm{d}x \,\mathrm{d}\xi
\\
& {}
  + \int_{\mathbb{H}} |f(x,\xi)|^2
    \cdot (1 + \xi)\cdot \mathfrak{w}(x,\xi)
    \,\mathrm{d}x \,\mathrm{d}\xi < \infty \,.
\end{aligned}
\end{equation*}
It easy to see that the Sobolev imbeddings
\begin{math}
  H^2(\mathbb{H};\mathfrak{w}) \hookrightarrow V =
  H^1(\mathbb{H};\mathfrak{w}) \hookrightarrow H =
  L^2(\mathbb{H};\mathfrak{w})
\end{math}
are bounded (i.e., continuous).

We will make use of the following {\em local\/} version of
the weighted Sobolev space
$H^2(\mathbb{H};\mathfrak{w})$:
Let $B_R(x_0,\xi_0)$ denote the open disc in $\RR^2$ with radius
$R > 0$ centered at the point $(x_0,\xi_0)\in \RR^2$.
If $\xi_0 = 0$, we define also the open upper half\--disc
\begin{equation*}
  B^{+}_R(x_0,0) =
  \{ (x,\xi)\in \RR^2\colon (x - x_0)^2 + \xi^2 < R^2 \,,
                     \hspace{7pt} \xi > 0
  \} \subset \HH \,.
\end{equation*}
Its closure in $\RR^2$ (hence, also in $\overline{\HH}$)
is denoted by
\begin{equation*}
  \overline{B}^{+}_R(x_0,0) =
  \{ (x,\xi)\in \RR^2\colon (x - x_0)^2 + \xi^2\leq R^2 \,,
                     \hspace{7pt} \xi\geq 0
  \} \subset \overline{\HH} \,.
\end{equation*}
We denote by $H^2(B^{+}_R(x_0,0); \mathfrak{w})$
the weighted Sobolev space of all functions
$f\in W^{2,2}_{\mathrm{loc}}(B^{+}_R(x_0,0))$
whose norm defined below is finite,
\begin{equation*}
%\label{loc:norm-H^2}
\begin{aligned}
  \left(
  \| f\|_{ H^2(B^{+}_R(x_0,0); \mathfrak{w}) }^{\sharp}
  \right)^2
& \eqdef \int_{ B^{+}_R(x_0,0) }
    \left( |f_{xx}|^2 + |f_{x\xi}|^2 + |f_{\xi\xi}|^2\right)
    \cdot \xi^2\cdot \mathfrak{w}(x,\xi)
    \,\mathrm{d}x \,\mathrm{d}\xi
\\
& {}
  + \int_{ B^{+}_R(x_0,0) }
    \left( |f_{x}|^2 + |f_{\xi}|^2\right)
    \cdot (1 + \xi^2)\cdot \mathfrak{w}(x,\xi)
    \,\mathrm{d}x \,\mathrm{d}\xi
\\
& {}
  + \int_{ B^{+}_R(x_0,0) } |f(x,\xi)|^2
    \cdot (1 + \xi)\cdot \mathfrak{w}(x,\xi)
    \,\mathrm{d}x \,\mathrm{d}\xi < \infty \,.
\end{aligned}
\end{equation*}
The half\--disc $B^{+}_R(x_0,0)$ being bounded in $\HH$,
this norm on $H^2(B^{+}_R(x_0,0); \mathfrak{w})$
is equivalent with the following simpler norm defined by
\begin{equation}
\label{e:norm-H^2}
\begin{aligned}
& \| f\|_{ H^2(B^{+}_R(x_0,0); \mathfrak{w}) }^2
  \eqdef \int_{ B^{+}_R(x_0,0) }
    \left( |f_{xx}|^2 + |f_{x\xi}|^2 + |f_{\xi\xi}|^2\right)
    \cdot \xi^{\beta + 1}\cdot \,\mathrm{d}x \,\mathrm{d}\xi
\\
& {}
  + \int_{ B^{+}_R(x_0,0) }
    \left( |f_{x}|^2 + |f_{\xi}|^2\right)
    \cdot \xi^{\beta - 1}\cdot \,\mathrm{d}x \,\mathrm{d}\xi
  + \int_{ B^{+}_R(x_0,0) } |f(x,\xi)|^2
    \cdot \xi^{\beta - 1}\cdot \,\mathrm{d}x \,\mathrm{d}\xi
    < \infty \,.
\end{aligned}
\end{equation}
We will employ the weighted Sobolev space
$H^2(B^{+}_R(x_0,0); \mathfrak{w})$
in Section~\ref{s:smooth_Hoelder}.

The weighted Sobolev space
$H^2(B^{+}_R(x_0,0); \mathfrak{w})$
will be imbedded into the weighted $L^p$-Lebesgue space
$L^p(B^{+}_R(x_0,0); \mathfrak{w})$ ($1\leq p < \infty$)
of all complex\--valued Lebesgue\--measurable functions
$f\colon B^{+}_R(x_0,0)\to \CC$ with the finite norm
\begin{equation}
\label{e:norm_L^p}
  \| f\|_{ L^p(B^{+}_R(x_0,0); \mathfrak{w}) }\eqdef
  \left( \int_{ B^{+}_R(x_0,0) } |f(x,\xi)|^p
         \cdot \xi^{\beta - 1}\cdot \,\mathrm{d}x \,\mathrm{d}\xi
  \right)^{1/p} < \infty \,.
\end{equation}

Finally, the local Schauder\--type regularity results near the boundary
\begin{math}
  \partial\HH = \RR\times \{ 0\} = \overline{\HH}\setminus \HH
\end{math}
of the half\--plane $\HH$ established in Section~\ref{s:smooth_Hoelder}
will be stated in the H\"older spaces
$C_s^{\alpha}( \overline{B}^{+}_R(x_0,0) )$ and
$C_s^{2+\alpha}( \overline{B}^{+}_R(x_0,0) )$
over any compact half\--disc
$\overline{B}^{+}_R(x_0,0)$ with $x_0\in \RR$ and $R\in (0,\infty)$.
The H\"older norm in these spaces corresponds to the so\--called
{\em cycloidal\/} Riemannian metric $s$ on $\HH$ defined by
\begin{math}
  \mathrm{d}s^2 = \xi^{-1} (\mathrm{d}x^2 + \mathrm{d}\xi^2) .
\end{math}
The associated {\em cycloidal distance function\/} on $\overline{\HH}$,
denoted by $s_{\mathrm{cycl}}(P_1,P_2)$ for two different points
$P_i = (x_i,\xi_i)\in \overline{\HH}$; $i=1,2$, is given by
\begin{equation*}
%\label{def:cyclo}
  s_{\mathrm{cycl}}(P_1,P_2)\eqdef
  \frac{ |x_1 - x_2| + |\xi_1 - \xi_2| }%
       { \sqrt{\xi_1} + \sqrt{\xi_2} +
         \sqrt{ |(x_1,\xi_1) - (x_2,\xi_2)| } } \,.
\end{equation*}
Of course, the expression
$|P_1 - P_1| = |(x_1,\xi_1) - (x_2,\xi_2)|$
stands for the Euclidean distance on~$\RR^2$.
We will use the following equivalent metric on $\overline{\HH}$
introduced in {\sc H.\ Koch} \cite[p.~11]{Koch-1999},
\begin{equation}
\label{def:cyclo}
  s(P_1,P_2)\eqdef
  \frac{ |(x_1,\xi_1) - (x_2,\xi_2)| }%
       { \sqrt{ \xi_1 + \xi_2 + |(x_1,\xi_1) - (x_2,\xi_2)| } } \,.
\end{equation}

As usual,
$C( \overline{B}^{+}_R(x_0,0) )$
denotes the Banach space of all continuous functions
$f\colon \overline{B}^{+}_R(x_0,0)\hfil\break \to \CC$
endowed with the maximum norm
\begin{equation*}
  \| f\|_{ C( \overline{B}^{+}_R(x_0,0) ) }\eqdef
  \max_{ (x,\xi)\in \overline{B}^{+}_R(x_0,0) } |f(x,\xi)| < \infty \,.
\end{equation*}
Given $\alpha\in (0,1)$, we denote by
$C_s^{\alpha}( \overline{B}^{+}_R(x_0,0) )$
the H\"older space of all functions
$f\in C( \overline{B}^{+}_R(x_0,0) )$ that satisfy
\begin{equation*}
  {[f]}_{ C_s^{\alpha}( \overline{B}^{+}_R(x_0,0) ) }\eqdef
  \sup_{ \genfrac{}{}{0pt}1{ P_1, P_2\in \overline{B}^{+}_R(x_0,0) }%
                           { P_1\neq P_2 }
    }\,
  \frac{ | f(P_1) - f(P_2) | }{ s(P_1,P_2)^{\alpha} }
  < \infty \,.
\end{equation*}
Recall that
$P_i = (x_i,\xi_i)\in \overline{\HH}$ for $i=1,2$.
The norm on this vector space is defined by
\begin{equation}
\label{norm:Hoelder}
\begin{aligned}
    \| f\|_{ C_s^{\alpha}( \overline{B}^{+}_R(x_0,0) ) }
  \eqdef
    \| f\|_{ C( \overline{B}^{+}_R(x_0,0) ) }
  + {[f]}_{ C_s^{\alpha}( \overline{B}^{+}_R(x_0,0) ) }
  < \infty \,.
\end{aligned}
\end{equation}
We denote by
$C_s^{2+\alpha}( \overline{B}^{+}_R(x_0,0) )$
its vector subspace consisting of all functions
$f\in C_s^{\alpha}( \overline{B}^{+}_R(x_0,0) )$
that are twice continuously differentiable in the open half\--disc
$B^{+}_R(x_0,0)$ and satisfy
\begin{equation}
\label{norm:2+alpha}
\begin{aligned}
    \| f\|_{ C_s^{2+\alpha}( \overline{B}^{+}_R(x_0,0) ) }
& \eqdef
    \| f\|_{ C_s^{\alpha}( \overline{B}^{+}_R(x_0,0) ) }
  + \| f_x    \|_{ C_s^{\alpha}( \overline{B}^{+}_R(x_0,0) ) }
  + \| f_{\xi}\|_{ C_s^{\alpha}( \overline{B}^{+}_R(x_0,0) ) }
\\
& {}
  + \|\xi\cdot f_{xx}(x,\xi)
    \|_{ C_s^{\alpha}( \overline{B}^{+}_R(x_0,0) ) }
  + \|\xi\cdot f_{x\xi}(x,\xi)
    \|_{ C_s^{\alpha}( \overline{B}^{+}_R(x_0,0) ) }
\\
& {}
  + \|\xi\cdot f_{\xi\xi}(x,\xi)
    \|_{ C_s^{\alpha}( \overline{B}^{+}_R(x_0,0) ) }
  < \infty \,.
\end{aligned}
\end{equation}
We endow $C_s^{2+\alpha}( \overline{B}^{+}_R(x_0,0) )$
with the norm
$\|\cdot\|_{ C_s^{2+\alpha}( \overline{B}^{+}_R(x_0,0) ) }$
defined above.
It is proved in
{\sc P.~M.~N.\ Feehan} and {\sc C.~A.\ Pop}
\cite{Feehan-Pop-13}, Lemma 3.1, Eq.\ (3.1), on p.~4409
(see also
 {\sc P.\ Daskalopoulos} and {\sc R.\ Hamilton}
 \cite{Daska-Hamil-98}, Prop.\ I.12.1 on p.~940)
that at every point
$P^{\ast} = (x^{\ast},0)\in \partial\HH$ with
$x^{\ast}\in (x_0 - R, x_0 + R)$ we have the zero limit
\begin{equation}
\label{lim:2+alpha}
  \lim_{ \genfrac{}{}{0pt}1{ P\to P^{\ast} }{ P\in B^{+}_R(x_0,0) }
       }\,
    \xi\cdot D^2 f(x,\xi) = 0 \quad\mbox{ for every }\,
    f\in C_s^{2+\alpha}( \overline{B}^{+}_R(x_0,0) ) \,,
\end{equation}
where $P = (x,\xi)\in \HH$ and
\begin{math}
  D^2 f = \left(
    \begin{matrix} f_{xx} , f_{x\xi}\\ f_{x\xi} , f_{\xi\xi}
    \end{matrix}
          \right)\in \RR^{2\times 2}
\end{math}
stands for the Hessian matrix of $f$ in $B^{+}_R(x_0,0)$
that consists of all second\--order partial derivatives of $f$.
This means that for any function $f\in C^2(B^{+}_R(x_0,0))$
the {\em\bfseries weighted\/} H\"older norm
\begin{math}
  \| f\|_{ C_s^{2+\alpha}( \overline{B}^{+}_R(x_0,0) ) } < \infty
\end{math}
forces the zero limit \eqref{lim:2+alpha}
which thus may be regarded as
an imposed homogeneous {\em\bfseries boundary condition\/}.

%%%%%%%%%%%%%%%%%%%%%%%%%%%%%%%%%%%%%%%%%%%%%%%%%%%%%%%%%%%%%%%%%%%%%%%
%%%%%    Formulation of the mathematical problem    %%%%%%%%%%%%%%%%%%%
%%%%%%%%%%%%%%%%%%%%%%%%%%%%%%%%%%%%%%%%%%%%%%%%%%%%%%%%%%%%%%%%%%%%%%%

\section{Formulation of the mathematical problem}
\label{s:MathProblem}

In this section we briefly describe {\sc S.~L.\ Heston}'s model
\cite[Sect.~1, pp.\ 328--332]{Heston}
and formulate the associated Cauchy problem as
an evolutionary equation of (degenerate) parabolic type.
A brief description of the ``economic'' model is provided in
Appendix~\ref{s:econom-Heston}.
The reader is referred to our earlier work in
{\sc B.\ Alziary} and {\sc P.\ Tak\'a\v{c}}
\cite[Sect.~2, pp.\ 6--13]{AlziaryTak}
for a more detailed analytical treatment of Heston's model.

%%%%%%%%%%%%%%%%%%%%%%%%%%%%%%%%%%%%%%%%%%%%%%%%%%%%%%%%%%%%%%%%%%%%%%%
%%%%%    Heston's stochastic volatility model    %%%%%%%%%%%%%%%%%%%%%%
%%%%%%%%%%%%%%%%%%%%%%%%%%%%%%%%%%%%%%%%%%%%%%%%%%%%%%%%%%%%%%%%%%%%%%%

\subsection{Heston's stochastic volatility model}
\label{ss:Heston}

We consider the Heston model given under
a \emph{risk neutral measure\/} via
equations $(1)-(4)$ in \cite[pp.\ 328--329]{Heston}.
The model is defined on a filtered probability space
$(\Omega, \mathcal{F}, (\mathcal{F}_t)_{t\geqslant 0}, \mathbb{P})$,
where $\mathbb{P}$ is a risk neutral probability measure, and
the filtration $(\mathcal{F}_t)_{t\geqslant 0}$
satisfies the usual conditions.
After a series of standard arguments based on \^{I}to's formula,
a (terminal value) Cauchy problem for the price of
a \emph{European call\/} or \emph{put option\/} is obtained
(see \cite[Eq.\ (2.4), p.~6]{AlziaryTak}).
This Cauchy problem is then transformed into
an initial value problem in the parabolic domain
$\mathbb{H}\times (0,T)\subset \mathbb{R}^3$
(\cite[Eq.\ (2.7), p.~8]{AlziaryTak})
with the (autonomous linear elliptic) {\it\bfseries Heston operator\/},
$\mathcal{A}$, given by
\cite[Eq.\ (2.9), p.~8]{AlziaryTak},
\begin{align}
\label{e:Heston-oper}
&
\begin{aligned}
  (\mathcal{A}u)(x,\xi) =
& {} - \frac{1}{2}\, \sigma\xi\cdot
  \left[ \frac{\partial}{\partial x}
  \left(
    \frac{\partial u}{\partial x}(x,\xi)
  + 2\rho\, \frac{\partial u}{\partial\xi}(x,\xi)
  \right)
  + \frac{\partial^2 u}{\partial\xi^2}(x,\xi)
  \right]
\\
& {}
  + \left( q_r + \genfrac{}{}{}1{1}{2} \sigma\xi\right)
    \cdot \frac{\partial u}{\partial x}(x,\xi)
  - \kappa (\theta_{\sigma} - \xi)
    \cdot \frac{\partial u}{\partial\xi}(x,\xi)
\end{aligned}
\\
\nonumber
&
\begin{aligned}
  \equiv
{} - \frac{1}{2}\, \sigma\xi\cdot
  \left[
  \left( u_x + 2\rho\, u_{\xi}\right)_x + u_{\xi\xi}
  \right]
   + \left( q_r + \genfrac{}{}{}1{1}{2} \sigma\xi\right) \cdot u_x
  {} - \kappa (\theta_{\sigma} - \xi) \cdot u_{\xi}
  \quad\mbox{ for $(x,\xi)\in \HH$, }
\end{aligned}
\end{align}
the {\it\bfseries boundary operator\/}, $\mathcal{B}$
(\cite[Eq.\ (2.10), p.~8]{AlziaryTak}),
on the boundary $\partial\HH\times (0,T)$,
and the boundary conditions as $x\to \pm\infty$ or $\xi\to +\infty$
(\cite[Eq.\ (2.11), p.~8]{AlziaryTak}).
Here, by $r-q\equiv - q_r$ $\in \RR$ we have abbreviated
the {\em\bfseries instantaneous drift of the stock price returns\/}
with $-\infty < r\leq q < \infty$, and by
$\theta_{\sigma}\equiv \theta / \sigma > 0$ the re\--scaled
{\em\bfseries long term\/} (or {\em\bfseries long\--run\/})
{\em\bfseries variance\/}
with $\theta, \sigma\in (0,\infty)$.
The {\em\bfseries correlation coefficient\/} $\rho$ satisfies
$\rho\in (-1,1)$.
Finally, $\kappa > 0$ denotes
the {\em\bfseries rate of mean reversion\/}; see
Eq.~\eqref{e:SVmodel} (Appendix~\ref{s:econom-Heston})
for motivation.
We now give a rigorous mathematical formulation of
this initial value Cauchy problem which follows
\cite[{\S}2.2, pp.\ 9--11]{AlziaryTak}.
Earlier motivation for formulation in similar
weighted Lebesgue and Sobolev spaces appears in
{\sc P.\ Daskalopoulos} and {\sc P.~M.~N.\ Feehan} \cite{Daska-Feehan-14}
and \cite[Sect.~2, p.~5048]{Daska-Feehan-16} and
{\sc P.~M.~N.\ Feehan} and {\sc C.~A.\ Pop} \cite{Feehan-Pop-15}.

We make use of the {\em Gel'fand triple\/}
$V\hookrightarrow H = H'\hookrightarrow V'$, i.e.,
we first identify the Hilbert space $H$ with its dual space $H'$,
by the Riesz representation theorem,
then use the imbedding $V\hookrightarrow H$,
which is dense and continuous, to construct its adjoint mapping
$H'\hookrightarrow V'$,
a dense and continuous imbedding of $H'$ into the dual space $V'$
of $V$ as well.
The (complex) inner product on $H$ induces
a sesquilinear duality between $V$ and $V'$;
we keep the notation
\begin{math}
  ( \,\cdot\, , \,\cdot\, )_H
\end{math}
also for this duality.
Now we define the linear operator $\mathcal{A}\colon V\to V'$
by the sesquilinear form
(cf.\ \cite[Eq.\ (2.21), p.~11]{AlziaryTak}),
for all $u,w\in V$,
\begin{align}
\label{e:Heston-diss:u=w}
&
\begin{aligned}
  (\mathcal{A}u, w)_H
& {}= \frac{\sigma}{2}\int_{\HH}
  \left(
    u_x\cdot \bar{w}_x + 2\rho\, u_{\xi}\cdot \bar{w}_x
  + u_{\xi}\cdot \bar{w}_{\xi}
  \right) \cdot \xi\cdot \mathfrak{w}(x,\xi)
    \,\mathrm{d}x \,\mathrm{d}\xi
\\
& {}+ \frac{\sigma}{2}\int_{\HH}
  (1 - \gamma\, \Sgn x)\, u_x\cdot \bar{w}
    \cdot \xi\cdot \mathfrak{w}(x,\xi)
    \,\mathrm{d}x \,\mathrm{d}\xi
\\
& {}+ \int_{\HH}
  \left( \kappa - \gamma\rho\sigma\, \Sgn x
       - \genfrac{}{}{}1{1}{2} \mu\sigma
  \right) u_{\xi}\cdot \bar{w}
    \cdot \xi\cdot \mathfrak{w}(x,\xi)
    \,\mathrm{d}x \,\mathrm{d}\xi
\end{aligned}
\\
\nonumber
&
\begin{aligned}
& {}+ q_r
    \int_{\HH} u_x\cdot \bar{w}\cdot \mathfrak{w}(x,\xi)
    \,\mathrm{d}x \,\mathrm{d}\xi
    + \left( \genfrac{}{}{}1{1}{2} \beta\sigma - \kappa\theta_{\sigma}
      \right)
    \int_{\HH} u_{\xi}\cdot \bar{w}\cdot \mathfrak{w}(x,\xi)
    \,\mathrm{d}x \,\mathrm{d}\xi \,.
\end{aligned}
\end{align}
All integrals on the right\--hand side converge absolutely for any pair
$u,w\in V$
(by the proof of Prop.\ 6.1 in \cite[pp.\ 21--23]{AlziaryTak}).

In order to derive the right\--hand side of Eq.~\eqref{e:Heston-diss:u=w}
from the left\--hand side which contains
the formal expression \eqref{e:Heston-oper} for $\mathcal{A}$
(see \cite[Eq.\ (2.20), p.~10]{AlziaryTak}),
the following {\it vanishing boundary conditions\/} are employed
(\cite[Eqs.\ (2.18), p.~9, and (2.19), p.~10]{AlziaryTak}):
\begin{align}
\label{bc_xi:Heston-bilin}
\left\{
\begin{aligned}
  \xi^{\beta}\cdot \int_{-\infty}^{+\infty}
    u_{\xi}(x,\xi)\cdot \bar{w}(x,\xi)
    \cdot \ee^{- \gamma |x|} \,\mathrm{d}x
& \,\longrightarrow\, 0 \quad\mbox{ as }\, \xi\to 0+ \,;
\\
  \xi^{\beta}\, \ee^{- \mu\xi}\cdot \int_{-\infty}^{+\infty}
    u_{\xi}(x,\xi)\cdot \bar{w}(x,\xi)
    \cdot \ee^{- \gamma |x|} \,\mathrm{d}x
& \,\longrightarrow\, 0 \quad\mbox{ as }\, \xi\to \infty \,,
\end{aligned}
\right.
\\
\label{bc_x:Heston-bilin}
    \ee^{- \gamma |x|}\cdot \int_0^{\infty}
  ( u_x + 2\rho\, u_{\xi} )\, \bar{w}(x,\xi)\cdot
    \xi^{\beta}\, \ee^{- \mu\xi} \,\mathrm{d}\xi
  \,\longrightarrow\, 0 \quad\mbox{ as }\, x\to \pm\infty \,,
\end{align}
for every function $w\in V$.
They are used in
{\sc B.\ Alziary} and {\sc P.\ Tak\'a\v{c}}
\cite[Eq.\ (2.20), p.~10]{AlziaryTak}
in order to perform integration by parts on
all second\--order partial derivatives of $u$ that appear in
the formal expression \eqref{e:Heston-oper} for $\mathcal{A}$
inserted into the inner product $(\mathcal{A}u, w)_H$
on the left\--hand side of Eq.~\eqref{e:Heston-diss:u=w}.
The boundary conditions in 
\eqref{bc_xi:Heston-bilin} and \eqref{bc_x:Heston-bilin}
are guaranteed by
the following (natural) {\it zero boundary conditions\/}
valid for every function
$w\in V = H^1(\mathbb{H};\mathfrak{w})$
(see \cite[Lemmas 10.2 and 10.3, pp.\ 44--45]{AlziaryTak}),
\begin{align}
\label{e:trace:v=0,infty}
\left\{
\begin{aligned}
  \xi^{\beta}\cdot \int_{-\infty}^{+\infty}
    |w(x,\xi)|^2\cdot \ee^{- \gamma |x|} \,\mathrm{d}x
& \,\longrightarrow\, 0 \quad\mbox{ as }\, \xi\to 0+ \,,
\\
  \xi^{\beta}\, \ee^{- \mu\xi}\cdot \int_{-\infty}^{+\infty}
    |w(x,\xi)|^2\cdot \ee^{- \gamma |x|} \,\mathrm{d}x
& \,\longrightarrow\, 0 \quad\mbox{ as }\, \xi\to \infty \,,
\end{aligned}
\right.
\\
\label{e:trace:x=+-infty}
    \mbox{ and }\quad
  \ee^{- \gamma |x|}\cdot \int_0^{\infty}
    |w(x,\xi)|^2\cdot \xi^{\beta}\, \ee^{- \mu\xi} \,\mathrm{d}\xi
  \,\longrightarrow\, 0 \quad\mbox{ as }\, x\to \pm\infty \,,
\end{align}
which are combined with the following additional boundary conditions
that we have to {\em impose\/}
(cf.\ \cite[Eqs.\ (2.23) and (2.24), p.~12]{AlziaryTak}):
\begin{align}
\label{bc_xi:bound_u}
\left\{
\begin{aligned}
  \xi^{\beta}\cdot \int_{-\infty}^{+\infty}
    |u_{\xi}(x,\xi)|^2\cdot \ee^{- \gamma |x|} \,\mathrm{d}x
& \leq \mathrm{const} < \infty \quad\mbox{ as }\, \xi\to 0+ \,;
\\
  \xi^{\beta}\, \ee^{- \mu\xi}\cdot \int_{-\infty}^{+\infty}
    |u_{\xi}(x,\xi)|^2\cdot \ee^{- \gamma |x|} \,\mathrm{d}x
& \leq \mathrm{const} < \infty \quad\mbox{ as }\, \xi\to \infty+ \,,
\end{aligned}
\right.
\\
\label{bc_x:bound_u}
    \ee^{- \gamma |x|}\cdot \int_0^{\infty}
  | u_x + 2\rho\, u_{\xi} |^2\cdot
    \xi^{\beta}\, \ee^{- \mu\xi} \,\mathrm{d}\xi
  \leq \mathrm{const} < \infty \quad\mbox{ as }\, x\to \pm\infty \,.
\end{align}
Indeed, we can apply the Cauchy\--Schwarz inequality to the integrals in
\eqref{bc_xi:Heston-bilin} and \eqref{bc_x:Heston-bilin}
to derive the zero limits from
\eqref{e:trace:v=0,infty}, \eqref{e:trace:x=+-infty},
\eqref{bc_xi:bound_u}, and \eqref{bc_x:bound_u}.

As we have just chosen a particular realization
$\mathcal{A}\colon V\to V'$
of the formal differential expression \eqref{e:Heston-oper}
defined by Eq.~\eqref{e:Heston-diss:u=w},
we no longer need to impose the boundary conditions
\eqref{bc_xi:bound_u} and \eqref{bc_x:bound_u}.

%%%%%%%%%%%%%%%%%%%%%%%%%%%%%%%%%%%%%%%%%%%%%%%%%%%%%%%%%%%%%%%%%%%%%%%
%%%%%    Cauchy problem in the weighted L^2 space H    %%%%%%%%%%%%%%%%
%%%%%%%%%%%%%%%%%%%%%%%%%%%%%%%%%%%%%%%%%%%%%%%%%%%%%%%%%%%%%%%%%%%%%%%

\subsection{The Cauchy problem in the weighted $L^2$-space $H$}
\label{ss:Cauchy-L^2}

The initial value Cauchy problem for the Heston model
mentioned in the previous paragraph ({\S}\ref{ss:Heston})
takes the following abstract form in the Hilbert space
$H = L^2(\mathbb{H};\mathfrak{w})$:
\begin{equation}
\label{e:Cauchy}
\left\{
\begin{alignedat}{2}
  \frac{\partial u}{\partial t} + \mathcal{A} u &= f(x,\xi,t)
  &&\quad\mbox{ in }\, \HH\times (0,T) \,;
\\
  u(x,\xi,0) &= u_0(x,\xi)
  &&\quad\mbox{ for }\, (x,\xi)\in \HH \,,
\end{alignedat}
\right.
\end{equation}
with the function $f(x,\xi,t)\equiv 0$ on the right\--hand side
and the initial data $u_0\in H$ at $t=0$.
The letter $T$ ($0 < T\leq +\infty$)
stands for an arbitrary (finite or infinite) upper bound on time $t$.
The (autonomous linear) {\it\bfseries Heston operator\/}
$\mathcal{A}\colon V\to V'$,
defined by the sesquilinear form \eqref{e:Heston-diss:u=w}
is bounded, by the Lax\--Milgram theorem.
Namely, the {\it\bfseries boundedness\/} and {\it\bfseries coercivity\/}
of this sesquilinear form are established in
\cite{AlziaryTak}, Prop.\ 6.1 on p.~21 and Prop.\ 6.2 on p.~23,
respectively, under certain restrictions on the constants which appear
in the weight $\mathfrak{w}$ and the operator $\mathcal{A}$
(see Eqs.\ \eqref{def:w} and \eqref{e:Heston-oper}).
We will discuss these rather fundamental restrictions in
Remark~\ref{rem-prop-Feller} at the end of this paragraph.

%%%%%%%%%%%%%%%%%%%%%%%%%%%%%%%%%%%%%%%%%%%%%%%%%%%%%%%%%%%%%%%%%%%%%%%
%%%%%    Weak solution (Definition)    %%%%%%%%%%%%%%%%%%%%%%%%%%%%%%%%
%%%%%%%%%%%%%%%%%%%%%%%%%%%%%%%%%%%%%%%%%%%%%%%%%%%%%%%%%%%%%%%%%%%%%%%
\begin{definition}\label{def-weak_sol}\nopagebreak
\begingroup\rm
\underline{\em Case\/} $0 < T < \infty$.
Let\/ $f\in L^2((0,T)\to V')$ and\/ $u_0\in H$.
A function $u\colon \HH\times [0,T]$ $\to \RR$
is called a {\it weak solution\/} to the initial value problem
\eqref{e:Cauchy}
if it has the following properties:
\begin{itemize}
\item[{\rm (i)}]
the mapping
$t\mapsto u(t)\equiv u(\,\cdot\,, \,\cdot\,, t)\colon [0,T]\to H$
is a continuous function, i.e.,
$u\in C([0,T]\to H)$;
\item[{\rm (ii)}]
the initial value $u(0) = u_0$ in $H$;
\item[{\rm (iii)}]
the mapping
$t\mapsto u(t)\colon (0,T)\to V$
is a B\^ochner square\--integrable function, i.e.,
$u\in L^2((0,T)\to V)$; and
\item[{\rm (iv)}]
for every function
\begin{equation*}
  \phi\in L^2((0,T)\to V)\cap W^{1,2}((0,T)\to V')
  \hookrightarrow C([0,T]\to H) \,,
\end{equation*}
the following equation holds,
\begin{equation*}
%\label{def:weak_sol}
\begin{aligned}
&   ( u(T), \phi(T) )_H
  - \int_0^T
    \left( u(t), \genfrac{}{}{}1{\partial\phi}{\partial t}(t)
    \right)_H \,\mathrm{d}t
  + \int_0^T (\mathcal{A}u(t), \phi(t))_H \,\mathrm{d}t
\\
& = ( u_0, \phi(0) )_H
  + \int_0^T (f(t), \phi(t))_H \,\mathrm{d}t \,.
\end{aligned}
\end{equation*}
\end{itemize}

\underline{\em Case\/} $T = +\infty$.
Let\/ $f\in L^2_{\mathrm{loc}}((0,\infty)\to V')$
(i.e., $f\in L^2((0,T_0)\to V')$ for every\/ $0 < T_0$ $< \infty$)
and let\/ $u_0\in H$.
A function $u\colon \HH\times [0,\infty)$ $\to \RR$
is called a {\it weak solution\/} to the initial value problem
\eqref{e:Cauchy} with $T = +\infty$,
if it is a {\it weak solution\/} to the initial value problem
\eqref{e:Cauchy} on every bounded time subinterval\/
$[0,T_0)\subset \RR_+$ with $0 < T_0 < T = +\infty$,
according to {\em Case\/} $0 < T < \infty$ above.
\endgroup
\end{definition}
%%%%%%%%%%%%%%%%%%%%%%%%%%%%%%%%%%%%%%%%%%%%%%%%%%%%%%%%%%%%%%%%%%%%%%%
\par\vskip 10pt

The following remarks are in order:

First, our definition of a weak solution is equivalent with that given in
{\sc L.~C.\ Evans} \cite[{\S}7.1]{Evans-98}, p.~352.
Here, for $0 < T < \infty$,
$W^{1,2}((0,T)\to V')$ denotes the Sobolev space of all functions
$\phi\in L^2((0,T)\to V')$ that possess a distributional time\--derivative
$\phi'\in L^2((0,T)\to V')$.
The norm is defined in the usual way; cf.\
{\sc L.~C.\ Evans} \cite[{\S}5.9]{Evans-98}.
The properties of
$V\equiv H^1(\mathbb{H};\mathfrak{w})$
justify the notation
$V'= H^{-1}(\mathbb{H};\mathfrak{w})$.
The continuity of the imbedding
\begin{equation*}
  L^2((0,T)\to V)\cap W^{1,2}((0,T)\to V')
  \hookrightarrow C([0,T]\to H)
\end{equation*}
is proved, e.g.,
in {\sc L.~C.\ Evans} \cite[{\S}5.9]{Evans-98}, Theorem~3 on p.~287.
We will see in Section~\ref{s:Main} that the initial value problem
\eqref{e:Cauchy} has a unique weak solution
$u\colon \HH\times [0,T]\to \RR$.

From now on, we use exclusively formula \eqref{e:Heston-diss:u=w}
to define the linear operator $\mathcal{A}\colon V\to V'$.
This means that we no longer need the boundary conditions in
\eqref{bc_xi:bound_u} and \eqref{bc_x:bound_u}
imposed on $u\in V$.

%%%%%%%%%%%%%%%%%%%%%%%%%%%%%%%%%%%%%%%%%%%%%%%%%%%%%%%%%%%%%%%%%%%%%%%
%%%%%    Feller's condition and another condition (Remark)    %%%%%%%%%
%%%%%%%%%%%%%%%%%%%%%%%%%%%%%%%%%%%%%%%%%%%%%%%%%%%%%%%%%%%%%%%%%%%%%%%
\begin{remark}\label{rem-prop-Feller}\nopagebreak
{\rm (Coercivity conditions.)}$\;$
\begingroup\rm
It is important to remark at this stage of our investigation of
the Heston operator $\mathcal{A}$ that,
in order to ensure the coercivity of $\mathcal{A} + c\, I$ on $V$,
one has to assume the well\--known {\it\bfseries Feller condition\/}
(\cite{Feller, Guo-Grzel-Ooster}),
\begin{equation}
\label{e:Feller}
  \genfrac{}{}{}1{1}{2} \sigma^2 - \kappa\theta < 0 \,.
\end{equation}

However, {\it Feller's condition\/} \eqref{e:Feller}
is {\it not sufficient\/} for obtaining the desired coercivity.
We need to guarantee also
\begin{equation*}
%\label{def:c_1'_max>0}
  c_{1,\mathrm{max}}'\eqdef \frac{1}{2}\, \sigma
    \left[
    \left( \frac{\kappa}{\sigma} - \gamma\, |\rho| \right)^2
  - \gamma (1+\gamma)
    \right] \geq 0 \,;
\end{equation*}
cf.\ Ineq.\ (6.15) in
{\sc B.\ Alziary} and {\sc P.\ Tak\'a\v{c}} \cite{AlziaryTak},
proof of Prop.\ 6.2, pp.\ 23--27.
That is, we need to assume the following
{\it\bfseries coercivity condition\/}:
\begin{equation}
\label{ineq:c_1'>0}
  \kappa\geq \sigma
    \left( \gamma\, |\rho| + \sqrt{ \gamma (1+\gamma) } \right)
  \quad \left(\, > \sigma\gamma (|\rho| + 1) \,\right) \,.
\end{equation}

The last inequality is an additional condition to
{\it Feller's condition\/},
$\genfrac{}{}{}1{1}{2} \sigma^2 - \kappa\theta < 0$,
both of them requiring the {\it rate of mean reversion\/}
$\kappa > 0$ of the stochastic volatility in Heston's model
to be sufficiently large.
This additional condition is caused by the fact that
{\sc W.\ Feller} \cite{Feller} considers only
an analogous problem in one space dimension ($\xi\in \RR_+$),
so that the solution $u = u(\xi)$ is independent from $x\in \RR$.
In particular, if the initial value
$u_0 = u(\,\cdot\,, \,\cdot\,, 0)\in H$ for $u(x,\xi,t)$
permits us to take $\gamma > 0$ arbitrarily small, then
inequality \eqref{ineq:c_1'>0} is easily satisfied,
provided {\it Feller's condition\/}
$\genfrac{}{}{}1{1}{2} \sigma^2 - \kappa\theta < 0$
is satisfied.
This is the case for a European {\em\bfseries put\/} option
with the initial condition
$u_0(x,\xi) = K\, (1 - \mathrm{e}^x)^{+}$ (${}\leq K$)
for $(x,\xi)\in \HH$.
However, if we wish to accommodate also initial values of type
$u_0(x,\xi) = K\, (\mathrm{e}^x - 1)^{+}$ for $(x,\xi)\in \HH$,
attached to a European {\em\bfseries call\/} option,
then we are forced to take $\gamma > 2$ to ensure that $u_0\in H$.

We refer the reader to the recent monograph by
{\sc G.~H. Meyer} \cite{Meyer-2015}
for a discussion of the role of Feller's condition
in the boundary conditions in Heston's model.
\hfill\Square
\endgroup
\end{remark}
%%%%%%%%%%%%%%%%%%%%%%%%%%%%%%%%%%%%%%%%%%%%%%%%%%%%%%%%%%%%%%%%%%%%%%%
\par\vskip 10pt

%%%%%%%%%%%%%%%%%%%%%%%%%%%%%%%%%%%%%%%%%%%%%%%%%%%%%%%%%%%%%%%%%%%%%%%
%%%%%    Main results    %%%%%%%%%%%%%%%%%%%%%%%%%%%%%%%%%%%%%%%%%%%%%%
%%%%%%%%%%%%%%%%%%%%%%%%%%%%%%%%%%%%%%%%%%%%%%%%%%%%%%%%%%%%%%%%%%%%%%%

\section{Main results}
\label{s:Main}

As our main results,
Theorems \ref{thm-Regular} and~\ref{thm-MaxPrinciple},
are only a~priori results for existing weak and strong solutions,
we state the following existence and uniqueness result taken from
our previous work
\cite[Prop.\ 4.1, p.~16]{AlziaryTak}.

%%%%%%%%%%%%%%%%%%%%%%%%%%%%%%%%%%%%%%%%%%%%%%%%%%%%%%%%%%%%%%%%%%%%%%%
%%%%%    J.-L. Lions' Proposition (Proposition)    %%%%%%%%%%%%%%%%%%%%
%%%%%%%%%%%%%%%%%%%%%%%%%%%%%%%%%%%%%%%%%%%%%%%%%%%%%%%%%%%%%%%%%%%%%%%
\begin{proposition}\label{prop-Lions}
Let\/
$\rho$, $\sigma$, $\theta$, $q_r$, and\/ $\gamma$
be given constants in $\RR$,
$\rho\in (-1,1)$, $\sigma > 0$, $\theta > 0$, and\/ $\gamma > 0$.
Assume that\/ $\kappa\in \RR$ is sufficiently large,
such that both inequalities,
\eqref{e:Feller} {\rm ({\it Feller's condition\/}) }
and \eqref{ineq:c_1'>0} are satisfied.
Set $\mu = \mu_{\mathrm{max}}$ where
\begin{equation}
\label{def:mu_max}
  \mu_{\mathrm{max}}\eqdef \frac{\kappa}{\sigma} - \gamma\, |\rho|
  \hspace{7pt} ({} > 0) \,; \quad\mbox{ hence, }\;
%\label{def:c_1'_max>0}
  c_{1,\mathrm{max}}'= \frac{1}{2}\, \sigma
    \left( \mu_{\mathrm{max}}^2 - \gamma (1+\gamma)
    \right) \geq 0 \,.
\end{equation}
Next, let us choose $\beta\in \RR$ such that\/
\begin{equation}
\label{ineq:beta-1<mu}
  1 < \beta\leq {2\kappa\theta} / {\sigma^2} \,.
\end{equation}
Let\/ $0 < T < \infty$, $f\in L^2((0,T)\to V')$, and\/ $u_0\in H$
be arbitrary.
Then the initial value problem \eqref{e:Cauchy} (with $u_0\in H$)
possesses a unique weak solution
\begin{equation*}
  u\in C([0,T]\to H)\cap L^2((0,T)\to V)
\end{equation*}
in the sense of\/ {\rm Definition~\ref{def-weak_sol}}.
Moreover, this solution satisfies also
$u\in W^{1,2}((0,T)\to V')$ and there exists a constant\/
$C\equiv C(T)\in (0,\infty)$, independent from $f$ and $u_0$, such that
\begin{equation*}
%\label{est:weak_sol}
    \sup_{t\in [0,T]} \| u(t)\|_H^2
  + \int_0^T \| u(t)\|_V^2 \,\mathrm{d}t
  + \int_0^T 
    \left\| \genfrac{}{}{}1{\partial u}{\partial t}(t)
    \right\|_{V'}^2 \,\mathrm{d}t
  \leq C
    \left( \| u_0\|_H^2
  + \int_0^T \| f(t)\|_{V'}^2 \,\mathrm{d}t \right) \,.
\end{equation*}

If\/ $T = +\infty$, $f\in L^2_{\mathrm{loc}}((0,\infty)\to V')$,
and\/ $u_0\in H$, the same existence and uniqueness result
(in the sense of\/ {\rm Definition~\ref{def-weak_sol}})
is valid with\/
\begin{equation*}
  u\in C([0,\infty)\to H)\cap L^2_{\mathrm{loc}}((0,\infty)\to V) \,.
\end{equation*}

Finally, if\/ $u_0\colon \HH\to \RR$ defined by
$u_0(x,\xi) = K\, (\mathrm{e}^x - 1)^{+}$, for\/ $(x,\xi)\in \HH$,
should belong to $H$, one needs to take $\gamma > 2$.
\end{proposition}
%%%%%%%%%%%%%%%%%%%%%%%%%%%%%%%%%%%%%%%%%%%%%%%%%%%%%%%%%%%%%%%%%%%%%%%
\par\vskip 10pt

The proof follows from
the {\it boundedness\/} and {\it coercivity\/}
of the sesquilinear form \eqref{e:Heston-diss:u=w} in $V\times V$
which are assumed in
{\sc J.-L.\ Lions} \cite[Chapt.~IV, {\S}1]{Lions-61},
inequalities (1.1) (p.~43) and (1.9) (p.~46), respectively.
For alternative proofs, see also e.g.\
{\sc L.~C.\ Evans} \cite[Chapt.~7, {\S}1.2(c)]{Evans-98},
Theorems 3 and~4, pp.\ 356--358,
{\sc J.-L.\ Lions} \cite[Chapt.~III, {\S}1.2]{Lions-71},
Theorem 1.2 (p.~102) and remarks thereafter (p.~103), or
{\sc A.\ Friedman} \cite{Friedman-64},
Chapt.~10, Theorem~17, p.~316.

Our first theorem contains global and local regularity results
for the weak solution
$u\colon \HH\times (0,T)\to \RR$ obtained in
Proposition~\ref{prop-Lions} above
for the special case $f\equiv 0$ in $\HH\times (0,T)$.
We formulate these regularity results using
the $C^0$-semigroup representation of the (unique) weak solution
\begin{math}
  u(\,\cdot\,,\,\cdot\,,t)\equiv u(t) = \ee^{-t\mathcal{A}} u_0\in H ,\
  t\in \RR_+ ,
\end{math}
to the homogeneous initial value problem \eqref{e:Cauchy}
(with $f\equiv 0$),
where we allow any $0 < T\leq +\infty$ and an arbitrary initial value
$u_0\in H$.
By the well\--known properties of $C^0$-semigroups,
\begin{math}
  \lambda I + \mathcal{A}\colon \mathcal{D}(\mathcal{A})\subset H\to H
\end{math}
is a closed linear operator in $H$ with the domain
$\mathcal{D}(\mathcal{A})\subset H$ which is invertible for all
$\lambda\in (\lambda_0,+\infty)$, with the bounded inverse
$(\lambda I + \mathcal{A})^{-1}\colon H\to H$.
We denote by
$\mathcal{D}\left( (\lambda I + \mathcal{A})^k \right)\subset H$
the domain of the $k$-th power of $\lambda I + \mathcal{A}$;
$k=1,2,3,\dots$.
Here, $\lambda_0\in (0,\infty)$ is a sufficiently large number
(called the {\em growth bound\/})
determined by the well\--known inequality \eqref{norm:e^(-tA)}
(in Section~\ref{s:smooth_Heston}).

The new result in this theorem is
a local Schauder\--type regularity result near
the boundary
\begin{math}
  \partial\HH\times (0,T) = \RR\times \{ 0\}\times (0,T)
\end{math}
of the parabolic domain $\HH\times (0,T)\subset \RR^3$
stated in the H\"older space
$C_s^{2+\alpha}( \overline{B}^{+}_R(x_0,0) )$
for every time $t\in (0,T)$.

%%%%%%%%%%%%%%%%%%%%%%%%%%%%%%%%%%%%%%%%%%%%%%%%%%%%%%%%%%%%%%%%%%%%%%%
%%%%%    Main Theorem (Theorem)    %%%%%%%%%%%%%%%%%%%%%%%%%%%%%%%%%%%%
%%%%%%%%%%%%%%%%%%%%%%%%%%%%%%%%%%%%%%%%%%%%%%%%%%%%%%%%%%%%%%%%%%%%%%%
\begin{theorem}\label{thm-Regular}
{\rm (Local and global regularity.)}$\;$
Let\/
$\rho$, $\sigma$, $\theta$, $q_r$, and\/ $\gamma$
be given constants in $\RR$,
$\rho\in (-1,1)$, $\sigma > 0$, $\theta > 0$, and\/ $\gamma > 0$.
Assume that\/ $\gamma$, $\kappa$, and\/ $\mu$
are chosen as specified in\/
{\rm Proposition~\ref{prop-Lions}} above and\/
$u_0\in H$ is arbitrary.
Finally, in addition to\/ {\rm Ineq.}~\eqref{ineq:beta-1<mu},
choose $\beta$ such that also $\beta(\beta - 1) < 4$, i.e.,
\begin{equation}
\label{ineq:beta.(beta-1)<4}
  1 < \beta\leq {2\kappa\theta} / {\sigma^2}
    \quad\mbox{ and }\quad
  \beta < \frac{1 + \sqrt{17}}{2} = 2.56\dots \;,
\end{equation}
respectively.
Then we have the following four statements for the weak solution
$u\colon (0,\infty)\to H$ obtained in\/
{\rm Proposition~\ref{prop-Lions}}:
\begin{itemize}
\item[{\rm (i)}]$\;$
\begin{math}
  u(\,\cdot\,,\,\cdot\,,t)\equiv u(t) = \ee^{-t\mathcal{A}} u_0
  \in \mathcal{D}_{\infty} = \bigcap_{k=1}^{\infty}
      \mathcal{D}\left( (\lambda I + \mathcal{A})^k \right)
  \subset H
\end{math}
holds for every $t\in (0,\infty)$.
\item[{\rm (ii)}]$\;$
\begin{math}
  u\in C^{\infty}( \HH\times (0,\infty) ) ,
\end{math}
i.e., $u$ is of class $C^{\infty}$ in $\HH\times (0,\infty)$.
Moreover, $u$ is a (local) classical solution of the parabolic equation
\begin{math}
%\label{e:Cauchy}
  \frac{\partial u}{\partial t} + \mathcal{A} u = 0
\end{math}
in the strong sense (pointwise) in $\HH\times (0,\infty)$.
\item[{\rm (iii)}]$\;$
Given $0 < T\leq +\infty$ and any\/ $x_0\in \mathbb{R}$,
there are a radius $R\in (0,\infty)$ and constants\/
$c_0, c_0'\in (0,\infty)$
such that, for every\/ $t\in (0,T)$, we have
\begin{math}
  u(t)\vert_{ \overline{B}^{+}_R(x_0,0) }
  \in C_s^{2+\alpha}( \overline{B}^{+}_R(x_0,0) )
\end{math}
and
\begin{equation*}
  U_{ \overline{B}^{+}_R(x_0,0) }(t)\colon u_0\longmapsto
  u(t)\vert_{ \overline{B}^{+}_R(x_0,0) }\colon
  H\longrightarrow C_s^{2+\alpha}( \overline{B}^{+}_R(x_0,0) )
\end{equation*}
is a bounded linear operator with the operator norm
\begin{math}
  \| U_{ \overline{B}^{+}_R(x_0,0) }(t) \|_{\mathrm{oper}}
  \leq (c_0' t^{-3} + c_0) \ee^{\lambda_0 t} .
\end{math}
\item[{\rm (iv)}]$\;$
Moreover, in the situation of\/ {\rm Part~(iii)} above, the mapping
\begin{equation*}
  t\,\longmapsto\,
    u(t)\vert_{ \overline{B}^{+}_R(x_0,0) }
  = U_{ \overline{B}^{+}_R(x_0,0) }(t) u_0
  \colon (0,T) \,\longrightarrow\,
    C_s^{2+\alpha}( \overline{B}^{+}_R(x_0,0) )
\end{equation*}
is continuous and differentiable, with
\begin{equation*}
  \| u(t+\tau) - u(t)\|_{ C_s^{2+\alpha}( \overline{B}^{+}_R(x_0,0) ) }
  \leq (c_0' t^{-3} + c_0) \ee^{\lambda_0 (t+\tau)}\cdot
       \| u(\tau) - u_0\|_H
\end{equation*}
and
\begin{equation*}
  \left\| \genfrac{}{}{}0{\partial u}{\partial t} (t)
  \right\|_{ C_s^{2+\alpha}( \overline{B}^{+}_R(x_0,0) ) }
  \leq (c_1' t^{-4} + c_1 t^{-1})
       \ee^{\lambda_0 t}\cdot \| u_0\|_H \,,
\end{equation*}
respectively, for all\/ $t\in (0,T)$ and for all\/
$\tau\in (0,\infty)$ such that\/ $t + \tau < T$.
Here,
$c_1, c_1'\in (0,\infty)$
are some other constants independent from $t$ and $u_0\in H$.
\end{itemize}
\end{theorem}
%%%%%%%%%%%%%%%%%%%%%%%%%%%%%%%%%%%%%%%%%%%%%%%%%%%%%%%%%%%%%%%%%%%%%%%
\par\vskip 10pt

Our {\em proof\/} of this theorem will be built up gradually
in the next two sections
(Sections~\ref{s:smooth_Heston} and~\ref{s:smooth_Hoelder})
and completed in Section~\ref{s:proofMain}.

We stress that the constants $c_0, c_0'\in (0,\infty)$
in {\rm Part~(iii)} do not depend on the choice of
$u_0\in H$ or $t\in (0,T)$.
However, the weighted norm on $H$ depends on the weight function
$\mathfrak{w}(x,\xi)$ which is not translation invariant
with respect to $x\in \RR$.
This property of $\mathfrak{w}$ means that the constants
$c_0, c_0'\in (0,\infty)$ may depend on $x_0\in \RR$.
We will see in the course of the proof of {\rm Part~(iii)}
(in Section~\ref{s:smooth_Hoelder})
that these constants are rendered independent from the length of
the time interval, $(0,T)$, $0 < T\leq +\infty$,
thanks to the multiplicative exponential factor $\ee^{\lambda_0 t}$.
The constant $\lambda_0\in (0,\infty)$
is determined solely by Ineq.~\eqref{norm:e^(-tA)}
(in Section~\ref{s:smooth_Heston}).
In particular, we obtain
\begin{math}
  \| U_{ \overline{B}^{+}_R(x_0,0) }(t) \|_{\mathrm{oper}}
  \leq (c_0' t^{-3} + c_0) \ee^{\lambda_0 t}
\end{math}
for all $t\in (0,T)$ (even if $T = +\infty$).

Concerning the behavior of the weak solution $u(x,\xi,t)$
to the Cauchy problem \eqref{e:Cauchy} in $\HH\times (0,T)$
near the boundary $\partial\HH\times (0,T)$,
{\rm Part~(iii)} of Theorem~\ref{thm-Regular}
has the following important consequence.

%%%%%%%%%%%%%%%%%%%%%%%%%%%%%%%%%%%%%%%%%%%%%%%%%%%%%%%%%%%%%%%%%%%%%%%
%%%%%    Main Theorem (Corollary)    %%%%%%%%%%%%%%%%%%%%%%%%%%%%%%%%%%
%%%%%%%%%%%%%%%%%%%%%%%%%%%%%%%%%%%%%%%%%%%%%%%%%%%%%%%%%%%%%%%%%%%%%%%
\begin{corollary}\label{cor-Regular}
{\rm (Boundary behavior.)}$\;$
Let\/ $0 < T\leq +\infty$ and\/ $t_0\in (0,T)$.
Under the hypotheses of\/ {\rm Theorem~\ref{thm-Regular}\/}, we have
\begin{math}
  u(\,\cdot\,,\,\cdot\,,t)\equiv u(t)\in C^1(\overline{\HH})
\end{math}
for every $t\in (0,T)$.
Furthermore, the function $u(x,\xi,t)$
verifies the following initial value Cauchy problem on
$\partial\HH\times (t_0,T)$,
\begin{equation}
\label{e:Cauchy:xi=0}
\left\{
\begin{alignedat}{2}
& \frac{\partial u}{\partial t}(x,0,t)
  + q_r
    \cdot \frac{\partial u}{\partial x}(x,0,t)
  - \kappa\theta_{\sigma}
    \cdot \frac{\partial u}{\partial\xi}(x,0,t)
  = 0
  &&\quad\mbox{ in }\, \RR\times (t_0,T) \,;
\\
& u(x,0,t_0) = u_{t_0}(x,0)
  &&\quad\mbox{ for }\, x\in \RR \,.
\end{alignedat}
\right.
\end{equation}
Here, we have denoted
$u_{t_0}(x,\xi)\eqdef u(x,\xi,t_0)$ 
for all $(x,\xi)\in \overline{\HH}$; hence,
\begin{math}
  u_{t_0}\equiv u(\,\cdot\,,\,\cdot\,,t_0)\in C^1(\overline{\HH}) .
\end{math}
This transport equation for the unknown function
$u(x,0,t)$ has a unique classical solution given by
\begin{equation}
\label{sol:Cauchy:xi=0}
\begin{aligned}
  u(x,0,t) = u(x - q_r (t-t_0),\, 0,t_0)
  + \kappa\theta_{\sigma}
    \int_{t_0}^t \frac{\partial u}{\partial\xi}
    \left( x - q_r (t-s),\, 0,s\right) \,\mathrm{d}s
\\
  \quad\mbox{ for }\, (x,t)\in \RR\times (t_0,T) \,.
\end{aligned}
\end{equation}
If, in addition, $u\in C^0(\overline{\HH}\times [0,T))$,
then we may take $t_0 = 0$ above, in Eqs.\
\eqref{e:Cauchy:xi=0} and \eqref{sol:Cauchy:xi=0}.
\end{corollary}
%%%%%%%%%%%%%%%%%%%%%%%%%%%%%%%%%%%%%%%%%%%%%%%%%%%%%%%%%%%%%%%%%%%%%%%
\par\vskip 10pt

This corollary will be {\em proved\/} in Section~\ref{s:proofMain}.

Our second theorem is a weak maximum principle which, in turn,
implies a pointwise bound on the weak solution
$u\colon \HH\times (0,T)\to \RR$ obtained in
Proposition~\ref{prop-Lions} above.
We begin with some auxiliary notation:

First, whenever $0 < T\leq +\infty$, let us denote by
$C^0(\HH\times [0,T))$ the vector space of all continuous functions
$u\colon \HH\times [0,T)\to \RR$ and by
$C^{2,1}(\HH\times (0,T))$
the vector space consisting of all continuous functions
$u\colon \HH\times (0,T)\to \RR$
that are continuously differentiable in $\HH\times (0,T)$ and also
twice continuously differentiable with respect to the space variables
$(x,\xi)\in \HH = \RR\times (0,\infty)$, i.e., all
\begin{math}
  u, u_t, u_x, u_{\xi}, u_{xx}, u_{x\xi}, u_{\xi\xi}
  \in C^0(\HH\times (0,T)) .
\end{math}

Second, let $\gamma_0\in (0,\infty)$
be an arbitrary constant, as large as needed.
Assuming {\it\bfseries Feller's condition\/} \eqref{e:Feller}, i.e.,
$\sigma^2 < 2\kappa\theta$, we allow any constants
$\beta_0, \mu_0\in (0,\infty)$ such that
\begin{equation}
\label{ineq:beta,mu}
  1\leq \beta_0 < {2\kappa\theta} / {\sigma^2}
  \quad\mbox{ and }\quad
  (0\leq {})\hspace{7pt} \beta_0 - 1 < \mu_0 < \infty \,.
\end{equation}
These two inequalities are motivated by
conditions \eqref{ineq:beta_1>1} and \eqref{ineq:mu>beta-1},
respectively, in the proof of the theorem below.
Notice that there is
{\em no upper bound\/} on the constant $\mu_0$.

Third, define a ``majorizing'' function
$\mathfrak{h}_0\colon \HH\to (0,\infty)$ by
\begin{equation}
\label{def:h_0}
\begin{aligned}
  \mathfrak{h}_0(x,\xi)
& {}\eqdef
    \exp\left[
        \gamma_0 (1 + x^2)^{1/2} + \mu_0\xi - (\beta_0 - 1)\, \ln\xi
        \right]
\\
& {}
  = \xi^{- (\beta_0 - 1)}\,
    \exp\left[ \gamma_0 (1 + x^2)^{1/2} + \mu_0\xi\right]
  \quad\mbox{ for $(x,\xi)\in \HH$. }
\end{aligned}
\end{equation}
A classical result on the weak maximum principle for
a parabolic Cauchy problem in $\RR^N\times (0,T)$
is valid under certain restrictions
on the growth of a strong solution $u(x,t)$ as $|x|\to \infty$,
$(x,t)\in \RR^N\times (0,T)$; see e.g.\
{\sc A.\ Friedman}
\cite[Chapt.~2, Sect.~4, Theorem~9, p.~43]{Friedman-64}.
Such restrictions in our case are reflected in the function
$\mathfrak{h}_0(x,\xi)$ introduced above.

Now we are ready to state our weak maximum principle.
This is an {\it a~priori\/} result for any strong subsolution $u$
to the parabolic Cauchy problem
\eqref{ineq:Cauchy}, \eqref{cond:Cauchy}, and \eqref{cond:Cauchy:t=0}
as described below.
As a consequence, we do {\it not\/} need to assume
hypothesis \eqref{ineq:c_1'>0} or \eqref{ineq:beta-1<mu}
(cf.\ Proposition~\ref{prop-Lions}).

%%%%%%%%%%%%%%%%%%%%%%%%%%%%%%%%%%%%%%%%%%%%%%%%%%%%%%%%%%%%%%%%%%%%%%%
%%%%%    Weak Maximum Principle (Theorem)    %%%%%%%%%%%%%%%%%%%%%%%%%%
%%%%%%%%%%%%%%%%%%%%%%%%%%%%%%%%%%%%%%%%%%%%%%%%%%%%%%%%%%%%%%%%%%%%%%%
\begin{theorem}\label{thm-MaxPrinciple}
{\rm (Weak maximum principle.)}$\;$
Let\/ $0 < T\leq +\infty$.
Assume that the constants\/
$\sigma, \kappa, \theta\in (0,\infty)$
satisfy the {\em\bfseries Feller condition\/} \eqref{e:Feller}.
Let\/ $\gamma_0\in (0,\infty)$ be arbitrary and assume that\/
$\beta_0, \mu_0\in (0,\infty)$ satisfy inequalities \eqref{ineq:beta,mu}.
Finally, assume that\/
$u\colon \HH\times [0,T)\to \RR$ is a function that satisfies\/
\begin{math}
  u\in C^0(\HH\times [0,T))\cap C^{2,1}(\HH\times (0,T))
\end{math}
together with\/
\begin{alignat}{2}
\label{ineq:Cauchy}
  \frac{\partial u}{\partial t} + \mathcal{A} u &\leq 0
&&  \quad\mbox{ in }\, \HH\times (0,T) \,,
\\
\label{cond:Cauchy}
  u(x,\xi,t) &\leq C\cdot \mathfrak{h}_0(x,\xi)
&&  \quad\mbox{ for }\, (x,\xi,t)\in \HH\times (0,T) \,;
\\
\label{cond:Cauchy:t=0}
  u(x,\xi,0) &\leq 0
&&  \quad\mbox{ for }\, (x,\xi)\in \HH \,,
\end{alignat}
where $C\in (0,\infty)$ is a positive constant independent from
$(x,\xi,t)\in \HH\times (0,T)$.

Then $u(x,\xi,t)\leq 0$ holds for all\/
$(x,\xi,t)\in \HH\times (0,T)$.
In particular, the Cauchy problem \eqref{e:Cauchy}
possesses at most one strong solution
\begin{math}
  u\in C^0(\HH\times [0,T))\cap C^{2,1}(\HH\times (0,T))
\end{math}
that satisfies the growth restriction
\begin{equation}
\label{cond:|u|<h_0}
  |u(x,\xi,t)|\leq C\cdot \mathfrak{h}_0(x,\xi)
  \quad\mbox{ for }\, (x,\xi,t)\in \HH\times [0,T) \,,
\end{equation}
where $C\in (0,\infty)$ is a positive constant.
\end{theorem}
%%%%%%%%%%%%%%%%%%%%%%%%%%%%%%%%%%%%%%%%%%%%%%%%%%%%%%%%%%%%%%%%%%%%%%%
\par\vskip 10pt

An important feature of this theorem is that there are
{\em no upper bounds\/} on the choice of the constants
$\gamma_0, \mu_0\in (0,\infty)$.
Once they have been chosen, the constant $\beta_0\in (0,\infty)$
must satisfy inequalities \eqref{ineq:beta,mu}.
Thus, any ``fast'' growth of the function $u(x,\xi,t)$,
as $x\to \pm\infty$ and/or $\xi\to +\infty$, of type
\begin{math}
  {}\leq \mathrm{const}\cdot \ee^{\gamma_0 |x| + \mu_0\xi}
\end{math}
is allowed in Ineq.~\eqref{cond:Cauchy}.
In contrast, as $x\to \pm\infty$ and $\xi\to 0+$,
the growth of $u(x,\xi,t)$ is limited to
\begin{math}
  {}\leq \mathrm{const}\cdot
    \xi^{- (\beta - 1)}\, \ee^{\gamma_0 |x|} .
\end{math}
A similar idea is offered by
Corollary~\ref{cor-MaxPrinciple} to
Theorem~\ref{thm-MaxPrinciple} below.
As we will infer from our proof of Corollary~\ref{cor-MaxPrinciple}
in Section~\ref{s:MaxPrinciple},
the case of $T < +\infty$ in Ineq.~\eqref{cond:Cauchy}
is of special importance.

Our {\it weak maximum principle\/} in Theorem~\ref{thm-MaxPrinciple}
differs from that in
{\sc P.~M.~N.\ Feehan} and {\sc C.~A.\ Pop} \cite{Feehan-Pop-13},
Lemma~3.4 on p.~4416.
Their conditions \cite[Eq.\ (3.29)]{Feehan-Pop-13}
imposed on the Heston operator $\mathcal{A}$ are weaker than ours.
We assume that the constants
$\sigma, \kappa, \theta\in (0,\infty)$
satisfy the {\em Feller condition\/} \eqref{e:Feller}.
On the other hand, we do not need that the functions
$u$, $u_t$, $u_x$, $u_{\xi}$, and
$\xi\, u_{xx}$, $\xi\, u_{x\xi}$, $\xi\, u_{\xi\xi}$
be continuous up to the boundary $\partial\mathbb{H}$
of the half\--plane $\HH = \RR\times (0,\infty)$; cf.\
\cite[Eq.\ (3.30)]{Feehan-Pop-13}.
Neither do we need the boundary condition in
\cite[Eq.\ (3.31)]{Feehan-Pop-13}.
In fact, we will show that this boundary condition is satisfied
also by our solutions to the Heston problem by combining our
growth hypothesis \eqref{cond:Cauchy} with Lemma~3.1 in
\cite[Eq.\ (3.1), p.~4409]{Feehan-Pop-13}.

%%%%%%%%%%%%%%%%%%%%%%%%%%%%%%%%%%%%%%%%%%%%%%%%%%%%%%%%%%%%%%%%%%%%%%%
%%%%%    Weak Maximum Principle (Corollary)    %%%%%%%%%%%%%%%%%%%%%%%%
%%%%%%%%%%%%%%%%%%%%%%%%%%%%%%%%%%%%%%%%%%%%%%%%%%%%%%%%%%%%%%%%%%%%%%%
\begin{corollary}\label{cor-MaxPrinciple}
Let\/ $0 < T\leq +\infty$, $\kappa\geq \sigma\rho$, and let\/
$r_0\in \RR_+$ satisfy\/
$r_0 + q_r = (r_0 - r) + q\geq 0$.
Assume that the constants\/
$\sigma, \kappa, \theta\in (0,\infty)$
satisfy the {\em Feller condition\/} \eqref{e:Feller}.
Let\/ $\gamma_0\in [1,\infty)$ be arbitrary and assume that\/
$\beta_0, \mu_0\in (0,\infty)$ satisfy inequalities \eqref{ineq:beta,mu}.
Finally, assume that\/
$u\colon \HH\times [0,T)\to \RR$ is a strong solution to
the homogeneous Cauchy problem \eqref{e:Cauchy} with $f\equiv 0$,
\begin{math}
  u\in C^0(\HH\times [0,T))\cap C^{2,1}(\HH\times (0,T)) ,
\end{math}
such that $u$ verifies the growth restriction \eqref{cond:|u|<h_0}
together with the following restriction at time\/ $t=0$,
\begin{equation}
\label{ineq:Cauchy:u_0}
  |u(x,\xi,0)|\leq U(x,\xi,0)
  \eqdef K_1\, \ee^{x + \varpi\xi} + K_0
  \quad\mbox{ for all }\, (x,\xi)\in \HH \,.
\end{equation}
Here, $K_0, K_1\in \RR_+$ are arbitrary constants, and\/
$\varpi\in \RR$ is another constant restricted by
\begin{equation}
\label{ineq:varpi}
  0\leq \varpi < \mu_0 \quad\mbox{ and }\quad
  (0\leq {})\hspace{7pt} \varpi\leq
    \min\left\{ \frac{r_0 + q_r}{ \kappa\theta_{\sigma} } ,\,
                \frac{ 2 (\kappa - \sigma\rho) }{\sigma}
        \right\} \,.
\end{equation}
Then
\begin{math}
  |u(x,\xi,t)|\leq U(x,\xi,t)\eqdef \ee^{r_0 t}\, U(x,\xi,0)
\end{math}
is valid in all of\/
$\HH\times [0,T)$, i.e., at all times $t\in [0,T)$.
\end{corollary}
%%%%%%%%%%%%%%%%%%%%%%%%%%%%%%%%%%%%%%%%%%%%%%%%%%%%%%%%%%%%%%%%%%%%%%%
\par\vskip 10pt

This corollary will be {\em proved\/}
in Section~\ref{s:MaxPrinciple}.

We remark that the condition in \eqref{ineq:Cauchy:u_0}
is satisfied for the initial value
$u_0\colon \HH\to \RR$ defined by
$u_0(x,\xi) = K\, (\mathrm{e}^x - 1)^{+}$, for\/ $(x,\xi)\in \HH$
(the European call option).
One may set $\varpi = 0$ together with $K_0 = 0$ and $K_1 = 1$.

We recall from Theorem~\ref{thm-Regular}, {\rm Part~(ii)},
that the (unique) weak solution
\begin{math}
  u(\,\cdot\,,\,\cdot\,,t)\equiv u(t) = \ee^{-t\mathcal{A}} u_0\in H ,\
  t\in \RR_+ ,
\end{math}
to the homogeneous initial value problem \eqref{e:Cauchy}
is of class $C^{\infty}$ in $\HH\times (0,\infty)$, i.e.,
\begin{math}
  u\in C^{\infty}( \HH\times (0,\infty) ) .
\end{math}
Thus, we conclude that $u$ verifies the parabolic equation
\begin{math}
%\label{e:Cauchy}
  \frac{\partial u}{\partial t} + \mathcal{A} u = 0
\end{math}
in the strong sense (pointwise) in $\HH\times (0,\infty)$, thanks to
$u\in C^{2,1}(\HH\times (0,\infty))$.
However, in order that $u$ be a strong (classical) solution of
problem \eqref{e:Cauchy} with $f\equiv 0$,
the additional continuity hypothesis
\begin{math}
  u\in C^0(\HH\times [0,T))
\end{math}
has to be made.

%%%%%%%%%%%%%%%%%%%%%%%%%%%%%%%%%%%%%%%%%%%%%%%%%%%%%%%%%%%%%%%%%%%%%%%
%%%%%    Some smoothing properties of the Heston semigroup    %%%%%%%%%
%%%%%%%%%%%%%%%%%%%%%%%%%%%%%%%%%%%%%%%%%%%%%%%%%%%%%%%%%%%%%%%%%%%%%%%

\section{Some smoothing properties of the Heston semigroup}
\label{s:smooth_Heston}

This section is concerned with some standard properties of
the $C^0$-semigroup $\ee^{-t\mathcal{A}}$ ($t\in \RR_+$)
of bounded linear operators $\ee^{-t\mathcal{A}}\colon H\to H$
on the complex Hilbert space $H$.
This semigroup has been already mentioned in
Section~\ref{s:Main}, Proposition~\ref{prop-Lions},
in connection with Theorem~\ref{thm-Regular}.
It is shown in
{\sc B.\ Alziary} and {\sc P.\ Tak\'a\v{c}} \cite{AlziaryTak},
Prop.\ 6.1 (p.~21) and Prop.\ 6.2 (p.~23), respectively,
that under conditions \eqref{e:Feller} and \eqref{ineq:c_1'>0}
the sesquilinear form
\begin{equation*}
  (u,w) \,\longmapsto\, \left( (\lambda I + \mathcal{A}) u, w\right)_H
  = (\mathcal{A}u, w)_H + \lambda\, (u,w)_H
    \colon V\times V\to \RR
\end{equation*}
defined in \eqref{e:Heston-diss:u=w}
(cf.\ \cite[Eq.\ (2.21), p.~11]{AlziaryTak})
is {\it bounded\/} and {\it coercive\/} on $V\hookrightarrow H$
(cf.\ also
 {\sc J.-L.\ Lions} \cite[Chapt.~IV, {\S}1]{Lions-61},
 inequalities (1.1) (p.~43) and (1.9) (p.~46), respectively),
provided
$\lambda\in (\lambda_0,\infty)$ where
$\lambda_0\in (0,\infty)$ is a sufficiently large constant, such that
\begin{equation*}
  \left( (\lambda_0 I + \mathcal{A}) u, u\right)_H
  = (\mathcal{A}u, u)_H + \lambda_0\, \| u\|_H^2 \geq 0
  \quad\mbox{ holds for all }\, u\in V \,.
\end{equation*}

We recall from Section~\ref{s:MathProblem}, {\S}\ref{ss:Cauchy-L^2},
that the (autonomous linear) operator
$\mathcal{A}\colon V\to V'$, defined by this sesquilinear form,
is bounded, by the Lax\--Milgram theorem.
The (unique) weak solution $u(x,\xi,t)$ to the Cauchy problem
\eqref{e:Cauchy} with $f\equiv 0$ in $\HH\times (0,T)$
and an arbitrary initial value $u_0\in H$
defines the $C^0$-semigroup representation of the solution
\begin{math}
  u(\,\cdot\,,\,\cdot\,,t)\equiv u(t) = \ee^{-t\mathcal{A}} u_0\in H ,\
  t\in \RR_+ ,
\end{math}
where $T\in (0,\infty)$ may be chosen arbitrarily large.

To be more precise, we denote by
\begin{math}
  -\mathcal{A}\colon \mathcal{D}(\mathcal{A})\subset H\to H
\end{math}
the infinitesimal generator of this semigroup which is the restriction
of the bounded linear operator $-\mathcal{A}\colon V\to V'$ to the domain
\begin{math}
  \mathcal{D}(\mathcal{A}) = \{ w\in V\colon \mathcal{A}w\in H\} .
\end{math}
In what follows, we keep the notation $\pm\mathcal{A}$
for this restriction.
It is verified in
\cite[Sect.~7, {\S}7.1, p.~29]{AlziaryTak}, that 
$\ee^{-t\mathcal{A}}$ ($t\in \RR_+$)
is a {\em holomorphic semigroup\/} of bounded linear operators on $H$
with the operator norm
\begin{equation}
\label{norm:e^(-tA)}
  \| \ee^{-t\mathcal{A}} \|_{ \mathcal{L}(H\to H) }
  \leq M_0\, \ee^{\lambda_0 t}
  \quad\mbox{ for all }\, t\in \RR_+ \,,
\end{equation}
by \cite[Ineq.\ (7.4), p.~29]{AlziaryTak}.
Here, $M_0\geq 1$ and $\lambda_0 > 0$ are some constants.
By a well\--known smoothing property of a holomorphic semigroup
({\sc A.\ Pazy} \cite[Eqs.\ (6.5)--(6.7), p.~70]{Pazy}),
we have
\begin{math}
  \ee^{-t\mathcal{A}} u_0\in
  \mathcal{D}\left( \mathcal{A}^k \right)\subset H
\end{math}
for all $u_0\in H$ and $t > 0$; $k=1,2,3,\dots$,
together with the bound on the operator norm
\begin{equation}
\label{norm:A^k.e^(-tA)}
  \| (\lambda I + \mathcal{A})^k\,
     \ee^{-t\mathcal{A}} \|_{ \mathcal{L}(H\to H) }
  \leq M_k t^{-k}\, \ee^{\lambda_0 t}
  \quad\mbox{ for all }\, t > 0 \,,
\end{equation}
where $M_k = (k M_1)^k > 0$ is a constant independent from time $t > 0$.
Indeed, employing the factorization
\begin{math}
    (\lambda I + \mathcal{A})^k\, \ee^{-t\mathcal{A}}
  = \left[
    (\lambda I + \mathcal{A})\, \ee^{- (t/k) \mathcal{A}}
    \right]^k
\end{math}
for $k=1,2,3,\dots$, we deduce the value
$M_k = (k M_1)^k$ in Ineq.~\eqref{norm:A^k.e^(-tA)}
for every $k=1,2,3,\dots$ from the case $k=1$.

Finally, the factorization
\begin{equation}
\label{e:A^(-k).A^k.e^(-tA)}
    \ee^{-t\mathcal{A}}
  = (\lambda I + \mathcal{A})^{-k}
  \left[
    (\lambda I + \mathcal{A})^k\, \ee^{-t\mathcal{A}}
  \right]
    \quad\mbox{ for $t > 0$ and $k=1,2,3,\dots$ }
\end{equation}
renders the smoothing property of
the holomorphic $C^0$-semigroup $\ee^{-t\mathcal{A}}$ ($t\in \RR_+$)
stated in the next lemma.
As usual, we endow the domain
\begin{math}
  \mathcal{D}(\mathcal{A}^k) =
  \mathcal{D}\left( (\lambda I + \mathcal{A})^k \right)
\end{math}
of the $k$-th power of $\mathcal{A}$ with its graph norm
for $k=1,2,3,\dots$
(see {\sc A.\ Pazy} \cite[Def.\ 6.7 and Thm.\ 6.8, p.~72]{Pazy}).
Hence, $\mathcal{D}(\mathcal{A}^k)$ is a Banach space
continuously imbedded into $H$, thanks to the graph of each
$\mathcal{A}^k$ being closed in $H\times H$.
Keeping the meaning of $M_k = (k M_1)^k$ from above, we get

%%%%%%%%%%%%%%%%%%%%%%%%%%%%%%%%%%%%%%%%%%%%%%%%%%%%%%%%%%%%%%%%%%%%%%%
%%%%%    Smoothing properties of Heston's semigroup (Lemma)    %%%%%%%%
%%%%%%%%%%%%%%%%%%%%%%%%%%%%%%%%%%%%%%%%%%%%%%%%%%%%%%%%%%%%%%%%%%%%%%%
\begin{lemma}\label{lem-smoothing}
{\rm (Smoothing property.)}$\;$
Under the hypotheses of\/ {\rm Theorem~\ref{thm-Regular}\/}
(cf.\ {\rm Proposition~\ref{prop-Lions}}), for any\/ $t > 0$
and every $k=1,2,3,\dots$, the bounded linear operator
$\ee^{-t\mathcal{A}}\colon H\to H$ maps $H$ into
$\mathcal{D}(\mathcal{A}^k)$ with the operator norm satisfying
\begin{equation}
\label{norm:A^(-k).e^(-tA)}
  \| \ee^{-t\mathcal{A}}
  \|_{ \mathcal{L}( H\to \mathcal{D}(\mathcal{A}^k) ) }
  \leq M_k t^{-k}\, \ee^{\lambda_0 t} \cdot
  \| (\lambda I + \mathcal{A})^{-k}
  \|_{ \mathcal{L}( H\to \mathcal{D}(\mathcal{A}^k) ) }
    \quad\mbox{ for all }\, t > 0 \,.
\end{equation}
\end{lemma}
%%%%%%%%%%%%%%%%%%%%%%%%%%%%%%%%%%%%%%%%%%%%%%%%%%%%%%%%%%%%%%%%%%%%%%%
\par\vskip 10pt

%\par\vskip 10pt
%%%%%%%%%%%%%%%%%
\proof
The estimate in Ineq.~\eqref{norm:A^(-k).e^(-tA)}
is obtained by applying \eqref{norm:A^k.e^(-tA)}
to the right\--hand side of \eqref{e:A^(-k).A^k.e^(-tA)}.
%\null\hfill\qed
\qed
%%%%%%%%%%%%%%%%%
\par\vskip 10pt

%%%%%%%%%%%%%%%%%%%%%%%%%%%%%%%%%%%%%%%%%%%%%%%%%%%%%%%%%%%%%%%%%%%%%%%
%%%%%    Smoothing properties in H\"older spaces    %%%%%%%%%%%%%%%%%%%
%%%%%%%%%%%%%%%%%%%%%%%%%%%%%%%%%%%%%%%%%%%%%%%%%%%%%%%%%%%%%%%%%%%%%%%

\section{Smoothing properties in H\"older spaces}
\label{s:smooth_Hoelder}

We apply Lemma~\ref{lem-smoothing} step by step for $k=1,2,3$.
We define the auxiliary functions $f_{j,k}(x,\xi,t)$
for $0\leq j\leq k\leq 3$ as follows:
First, for any time $t > 0$ we set
\begin{equation}
\label{def:f_0,k(t)}
    f_{0,k}(\,\cdot\,,\,\cdot\,,t)\equiv f_{0,k}(t)\eqdef
    (\lambda I + \mathcal{A})^k \ee^{-t\mathcal{A}} u_0\in H \,,\quad
    k=1,2,3 \,.
\end{equation}
Next, for $t > 0$ we introduce
\begin{equation}
\label{def:f_j,k(t)}
    f_{j,k}(\,\cdot\,,\,\cdot\,,t)\equiv f_{j,k}(t)\eqdef
    (\lambda I + \mathcal{A})^{-j} f_{0,k}(t)
  \quad\mbox{ for }\, 1\leq j\leq k\leq 3 \,.
\end{equation}
Clearly, for $j=k$; $k=1,2,3$, and $t > 0$ we obtain
$f_{k,k}(t) = u(t)$.

%%%%%%%%%%%%%%%%%%%%%%%%%%%%%%%%%%%%%%%%%%%%%%%%%%%%%%%%%%%%%%%%%%%%%%%
%%%%%    Smoothing with $(\lambda I + \mathcal{A})^{-1}$    %%%%%%%%%%%
%%%%%%%%%%%%%%%%%%%%%%%%%%%%%%%%%%%%%%%%%%%%%%%%%%%%%%%%%%%%%%%%%%%%%%%

\subsection{Smoothing with the factor $(\lambda I + \mathcal{A})^{-1}$}
\label{ss:smooth_Hoelder-1}

For $j=1$ and $k=1,2,3$ we get
\begin{math}
  (\lambda I + \mathcal{A}) f_{1,k}(t) = f_{0,k}(t)\in H ,\
  t > 0 \,.
\end{math}
We apply an interior (local) $H^2$-type regularity result due to
{\sc P.~M.~N.\ Feehan} and {\sc C.~A.\ Pop}
\cite{Feehan-Pop-15}, Theorem 3.16, Eq.\ (3.12), on p.~385
(stated in Lemma~\ref{lem-smooth_Hoelder-1},
 Appendix~\ref{s:ellipt_regul}),
to conclude that
\begin{math}
  f_{1,k}(t)\in H^2\left( B_{R_1}^{+}(x_0,0); \mathfrak{w}\right)
\end{math}
holds with any radius $R_1\in (0,\infty)$.
More precisely, there is a constant $C_1 > 0$ depending only on
the center point $x_0\in \mathbb{R}$ and the radii
$0 < R_1 < R_0 < \infty$, but independent from
$u_0\in H$ and $t > 0$, such that
\begin{equation}
\label{est:f_1,k(t)}
\begin{aligned}
  \| f_{1,k}(t)\|_{ H^2\left( B_{R_1}^{+}(x_0,0); \mathfrak{w}\right) }
& {}
  \leq C_1
  \left(
    \| f_{0,k}(t)\|_{ L^2\left( B_{R_0}^{+}(x_0,0); \mathfrak{w}\right) }
  + \| f_{1,k}(t)\|_{ L^2\left( B_{R_0}^{+}(x_0,0); \mathfrak{w}\right) }
  \right) \,.
\end{aligned}
\end{equation}
(The weighted Sobolev norm on the left\--hand side has been introduced
 in Eq.~\eqref{e:norm-H^2}.)

For $k=1$ we take advantage of the well\--known fact that
the operator norms of the family of bounded linear operators
\begin{math}
  t (\lambda I + \mathcal{A})\, \ee^{-t\mathcal{A}}\colon H\to H
\end{math}
are bounded above by
$M_1\, \ee^{\lambda_0 t}$ for all $t > 0$,
by Ineq.~\eqref{norm:A^k.e^(-tA)}.
Consequently, we get the estimate
\begin{equation}
\label{est:f_0,1(t)}
  \| f_{0,1}(t)\|_{ L^2\left( B_{R_0}^{+}(x_0,0); \mathfrak{w}\right) }
  \leq \| f_{0,1}(t)\|_H
  \leq M_1 t^{-1} \ee^{\lambda_0 t}\, \| u_0\|_H
    \quad\mbox{ for }\, t > 0 \,.
\end{equation}
Recalling
$u(t) = f_{1,1}(t) = \ee^{-t\mathcal{A}} u_0$
with the operator norms
\begin{math}
  \| \ee^{-t\mathcal{A}} \|_{ \mathcal{L}(H\to H) }
  \leq M_0\, \ee^{\lambda_0 t}
\end{math}
for $t > 0$, by Ineq.~\eqref{norm:e^(-tA)}, and applying
\eqref{est:f_0,1(t)} to \eqref{est:f_1,k(t)} to deduce
\begin{equation*}
%\label{est:f_1,1(t)}
\begin{aligned}
  \| u(t)\|_{ H^2\left( B_{R_1}^{+}(x_0,0); \mathfrak{w}\right) }
& \leq C_1
  \left( M_1 t^{-1} \ee^{\lambda_0 t}\, \| u_0\|_H + \| u(t)\|_H
  \right)
\\
& \leq \left( C_{1,1} t^{-1} + C_{1,0}\right)
       \ee^{\lambda_0 t}\, \| u_0\|_H
    \quad\mbox{ for all }\, t\in (0,\infty) \,.
\end{aligned}
\end{equation*}
The constants $C_{1,1}, C_{1,0} > 0$ are given by
$C_{1,1} = C_1 M_1$ and $C_{1,0} = C_1 M_0$.

We conclude that, for every $t > 0$,
\begin{math}
  u_0 \,\longmapsto\, u(t)\vert_{ B_{R_1}^{+}(x_0,0) }\colon
  H\hookrightarrow
  H^2\left( B_{R_1}^{+}(x_0,0); \mathfrak{w}\right)
\end{math}
is a bounded linear operator with the operator norm bounded above by
\begin{math}
  \left( C_{1,1} t^{-1} + C_{1,0}\right) \ee^{\lambda_0 t} .
\end{math}
%

%%%%%%%%%%%%%%%%%%%%%%%%%%%%%%%%%%%%%%%%%%%%%%%%%%%%%%%%%%%%%%%%%%%%%%%
%%%%%    Smoothing with $(\lambda I + \mathcal{A})^{-2}$    %%%%%%%%%%%
%%%%%%%%%%%%%%%%%%%%%%%%%%%%%%%%%%%%%%%%%%%%%%%%%%%%%%%%%%%%%%%%%%%%%%%

\subsection{Smoothing with the factor $(\lambda I + \mathcal{A})^{-2}$}
\label{ss:smooth_Hoelder-2}

Now we take $j=2$ and $k=2,3$.
Hence, we get
\begin{math}
  (\lambda I + \mathcal{A}) f_{2,k}(t) = f_{1,k}(t)
  \in H^2\left( B_{R_1}^{+}(x_0,0); \mathfrak{w}\right) ,
\end{math}
$t > 0$.
Thanks to our hypotheses $\beta > 1$ and
$\beta (\beta - 1) < 4$ in Ineq.~\eqref{ineq:beta.(beta-1)<4},
there is a number $p > 4$ such that
$2+\beta < p < 2 + \frac{4}{\beta - 1}$.
In particular, Ineq.~\eqref{cond:beta} is valid.
For instance, if $1 < \beta < 2$, we may choose $p=6$.
By Lemma~\ref{lem-Sobol-Morrey}, Ineq.~\eqref{ineq:H^2->Castro_p}
(Appendix~\ref{s:Sobolev_trace}),
the {\em restricted imbedding\/}
\begin{equation}
\label{e:H^2->L^p}
  u\vert_{ B_{R_1}^{+}(x_0,0) }\longmapsto
  u\vert_{ B_{R_1'}^{+}(x_0,0) }\colon
  H^2\left( B_{R_1}^{+}(x_0,0); \mathfrak{w}\right) \hookrightarrow
  L^p\left( B_{R_1'}^{+}(x_0,0); \mathfrak{w}\right)
\end{equation}
is continuous, whenever
$R_1'= R_1 / 2$ and $0 < R_1 < R_0 < \infty$.
(We refer to {\sc R.~A.\ Adams} and {\sc J.~J.~F. Fournier}
 \cite[Chapt.~6, {\S}6.1, p.~167]{AdamsFournier}
 for the definition of a {\em restricted imbedding\/}
 concerning Sobolev and Lebesgue function spaces.
 Typically, a restricted imbedding is {\em not\/} injective.)
Consequently, we have also
\begin{math}
  (\lambda I + \mathcal{A}) f_{2,k}(t) = f_{1,k}(t)
  \in L^p\left( B_{R_1'}^{+}(x_0,0); \mathfrak{w}\right) ,\
  t > 0 \,.
\end{math}
This is an elliptic equation for the unknown function
$f_{2,k}(t)\in V\hookrightarrow H$.
This observation allows us to apply
a local H\"older regularity result from
{\sc P.~M.~N.\ Feehan} and {\sc C.~A.\ Pop}
\cite{Feehan-Pop-17}, Theorem 1.11, Eq.\ (1.31), on p.~1083
(stated in Lemma~\ref{lem-smooth_Hoelder-2},
 Appendix~\ref{s:ellipt_regul}; see also
 \cite{Feehan-Pop-15}, Theorem 2.5, Eq.\ (2.12), pp.\ 375--376)
in order to derive
\begin{math}
  f_{2,k}(t)\in C_s^{\alpha}\left( \overline{B}_{R_2}^{+}(x_0,0)\right)
\end{math}
for $t\in (0,T)$, together with the estimate
\begin{equation}
\label{est:f_2,k(t)}
\begin{aligned}
  \| f_{2,k}(t)
  \|_{ C_s^{\alpha}\left( \overline{B}_{R_2}^{+}(x_0,0)\right) }
& {}
  \leq C_2
  \left(
    \| f_{1,k}(t)\|_{ L^p\left( B_{R_1'}^{+}(x_0,0); \mathfrak{w}\right) }
  + \| f_{2,k}(t)\|_{ L^2\left( B_{R_1'}^{+}(x_0,0); \mathfrak{w}\right) }
  \right)
\end{aligned}
\end{equation}
with some $0 < R_2 < R_1'$.
In this local H\"older regularity result
(Lemma~\ref{lem-smooth_Hoelder-2}),
only the condition $p > \max\{ 4,\, 2+\beta\}$ is needed.
(The H\"older norm on the left\--hand side has been introduced
 in Eq.~\eqref{norm:Hoelder}.)
All constants $R_2\equiv R_2(R_1')$ ($0 < R_2 < R_1'$),
$\alpha\in (0,1)$, and $C_2 > 0$ depend on
the center point $x_0\in \mathbb{R}$ and the radius $R_1'$
with $R_1'= R_1 / 2$ (${} < R_1 < R_0 < \infty$),
but are independent from $u_0\in H$ and $t > 0$.
We now employ
Lemma~\ref{lem-Sobol-Morrey}, Ineq.~\eqref{ineq:H^2->Castro_p}
(Appendix~\ref{s:Sobolev_trace}),
again to estimate the norm of the restricted Sobolev imbedding
in \eqref{e:H^2->L^p},
\begin{equation*}
%\label{ineq:f_1,k(t)}
  \| f_{1,k}(t)\|_{ L^p\left( B_{R_1'}^{+}(x_0,0); \mathfrak{w}\right) }
  \leq C'(R_1)\,
    \| f_{1,k}(t)
    \|_{ H^2\left( B^{+}_{R_1}(x_0,0); \mathfrak{w}\right) } \,,
\end{equation*}
where $0 < C'(R_1) < \infty$ is a constant depending only on
the center point $x_0\in \mathbb{R}$ and the radius $R_1 > 0$,
but neither on $u_0\in H$ nor on $t > 0$.
We combine the last estimate with \eqref{est:f_1,k(t)}
in order to estimate the right\--hand side of
Ineq.~\eqref{est:f_2,k(t)} by
\begin{equation}
\label{ineq:f_2,k(t)}
\begin{aligned}
& \| f_{2,k}(t)
  \|_{ C_s^{\alpha}\left( \overline{B}_{R_2}^{+}(x_0,0)\right) }
\\
& \leq C_2\cdot C'(R_1)\,
    \| f_{1,k}(t)
    \|_{ H^2\left( B^{+}_{R_1}(x_0,0); \mathfrak{w}\right) }
  + C_2\,
    \| f_{2,k}(t)
    \|_{ L^2\left( B_{R_1'}^{+}(x_0,0); \mathfrak{w}\right) }
\\
& {}
  \leq C_1 C_2\cdot C'(R_1)\,
  \left(
    \| f_{0,k}(t)\|_{ L^2\left( B_{R_0}^{+}(x_0,0); \mathfrak{w}\right) }
  + \| f_{1,k}(t)\|_{ L^2\left( B_{R_0}^{+}(x_0,0); \mathfrak{w}\right) }
  \right)
\\
& {}
  + C_2\,
    \| f_{2,k}(t)
    \|_{ L^2\left( B_{R_1'}^{+}(x_0,0); \mathfrak{w}\right) }
\end{aligned}
\end{equation}
with some $0 < R_2 < R_1'= R_1 / 2$ and $0 < R_1 < R_0 < \infty$.

For $k=2$ we now employ the fact that
the operator norms of both families of bounded linear operators,
\begin{math}
  t (\lambda I + \mathcal{A})\, \ee^{-t\mathcal{A}}\colon H\to H
\end{math}
and
\begin{math}
  t^2 (\lambda I + \mathcal{A})^2\, \ee^{-t\mathcal{A}}\colon H\to H ,
\end{math}
are bounded above by $M_1\, \ee^{\lambda_0 t}$ and
$(2 M_1)^2\, \ee^{\lambda_0 t}$, respectively,
by Ineq.~\eqref{norm:A^k.e^(-tA)}, i.e., by
$\mathrm{const}\cdot \ee^{\lambda_0 t}$ for all $t > 0$.
We thus estimate
\begin{align}
\label{est:f_0,2(t)}
& \| f_{0,2}(t)\|_{ L^2\left( B_{R_0}^{+}(x_0,0); \mathfrak{w}\right) }
  \leq \| f_{0,2}(t)\|_H
  \leq (2 M_1)^2 t^{-2}\, \ee^{\lambda_0 t}\, \| u_0\|_H
    \quad\mbox{ and }
\\
\label{est:f_1,2(t)}
&
\begin{aligned}
&   \| f_{1,2}(t)\|_{ L^2\left( B_{R_0}^{+}(x_0,0); \mathfrak{w}\right) }
  = \| f_{0,1}(t)\|_{ L^2\left( B_{R_0}^{+}(x_0,0); \mathfrak{w}\right) }
\\
& \leq \| f_{0,1}(t)\|_H
  \leq M_1 t^{-1}\, \ee^{\lambda_0 t}\, \| u_0\|_H
    \quad\mbox{ for }\, t\in (0,\infty) \,.
\end{aligned}
\end{align}
The latter estimate follows directly from
$f_{1,2}(t) = f_{0,1}(t)$ and Ineq.~\eqref{est:f_0,1(t)}.
Consequently, recalling
$u(t) = f_{2,2}(t) = \ee^{-t\mathcal{A}} u_0$
with the operator norms
\begin{math}
  \| \ee^{-t\mathcal{A}} \|_{ \mathcal{L}(H\to H) }
  \leq M_0\, \ee^{\lambda_0 t}
\end{math}
for $t > 0$, we apply the estimates in
\eqref{est:f_0,2(t)} and \eqref{est:f_1,2(t)}
to \eqref{ineq:f_2,k(t)}, thus arriving at
\begin{equation}
\label{est:f_2,2(t)}
\begin{aligned}
& \| u(t)\|_{ C_s^{\alpha}\left( \overline{B}_{R_2}^{+}(x_0,0)\right) }
\\
& {}
  \leq C_1 C_2\cdot C'(R_1)\left( (2 M_1)^2 t^{-2} + M_1 t^{-1}\right)
       \ee^{\lambda_0 t}\, \| u_0\|_H
  + C_2\, \| u(t)\|_H
\\
& {}
  \leq C_1 C_2\cdot C'(R_1)\left( (2 M_1)^2 t^{-2} + M_1 t^{-1}\right)
       \ee^{\lambda_0 t}\, \| u_0\|_H
  + C_2 M_0\, \ee^{\lambda_0 t}\, \| u_0\|_H
\\
& {}
  = \left( C_{2,2} t^{-2} + C_{2,1} t^{-1} + C_{2,0}\right)
    \ee^{\lambda_0 t}\, \| u_0\|_H
    \quad\mbox{ for all }\, t\in (0,\infty) \,.
\end{aligned}
\end{equation}
The constants $C_{2,j} > 0$; $j=0,1,2$, are given by
$C_{2,2} = C_1 C_2\cdot C'(R_1) (2 M_1)^2$,
$C_{2,1} = C_1 C_2\cdot C'(R_1) M_1$, and
$C_{2,0} = C_2 M_0$.

We have shown that, for every $t > 0$,
\begin{math}
  u_0 \,\longmapsto\, u(t)\vert_{ \overline{B}_{R_2}^{+}(x_0,0) }
  \colon H\hookrightarrow
  C_s^{\alpha}\left( \overline{B}_{R_2}^{+}(x_0,0)\right)
\end{math}
is a bounded linear operator with the operator norm bounded above by
\hfil\break
\begin{math}
    \left( C_{2,2} t^{-2} + C_{2,1} t^{-1} + C_{2,0}\right)
    \ee^{\lambda_0 t} .
\end{math}
%

%%%%%%%%%%%%%%%%%%%%%%%%%%%%%%%%%%%%%%%%%%%%%%%%%%%%%%%%%%%%%%%%%%%%%%%
%%%%%    Smoothing with $(\lambda I + \mathcal{A})^{-3}$    %%%%%%%%%%%
%%%%%%%%%%%%%%%%%%%%%%%%%%%%%%%%%%%%%%%%%%%%%%%%%%%%%%%%%%%%%%%%%%%%%%%

\subsection{Smoothing with the factor $(\lambda I + \mathcal{A})^{-3}$}
\label{ss:smooth_Hoelder-3}

Here, we take $j=k=3$, that is, we factorize
\begin{math}
  u(t) = f_{3,3}(t) =
  (\lambda I + \mathcal{A})^{-1}\, f_{2,3}(t)
\end{math}
with
\begin{math}
  f_{2,3}(t) = f_{0,1}(t) =
  (\lambda I + \mathcal{A})\, \ee^{-t\mathcal{A}} \in H .
\end{math}
In Paragraph {\S}\ref{ss:smooth_Hoelder-2}
above we have obtained the local H\"older regularity
\begin{math}
  f_{2,3}(t)\in
  C_s^{\alpha}\left( \overline{B}_{R_2}^{+}(x_0,0)\right)
\end{math}
for $t\in (0,\infty)$, together with
the estimate \eqref{ineq:f_2,k(t)} ($k=3$).
Applying Ineq.~\eqref{norm:A^k.e^(-tA)} to
\eqref{ineq:f_2,k(t)} with $k=3$, where
$f_{1,3} = f_{0,2}$ and $f_{2,3} = f_{0,1}$, we obtain further
\begin{equation}
\label{est:f_2,3(t)}
\begin{aligned}
& \| f_{2,3}(t)
  \|_{ C_s^{\alpha}\left( \overline{B}_{R_2}^{+}(x_0,0)\right) }
\\
& \leq C_1 C_2\cdot C'(R_1)
       \left( (3 M_1)^3 t^{-3} + (2 M_1)^2 t^{-2}\right)
       \ee^{\lambda_0 t}\, \| u_0\|_H
  + C_2 M_1 t^{-1}\, \ee^{\lambda_0 t}\, \| u_0\|_H
\\
& {}
  = \left( c_{3,3} t^{-3} + c_{3,2} t^{-2} + c_{3,1} t^{-1}\right)
       \ee^{\lambda_0 t}\, \| u_0\|_H
    \quad\mbox{ for all }\, t\in (0,\infty) \,,
\end{aligned}
\end{equation}
with some constant $R_2\in \RR$ satisfying
$0 < R_2 < R_1'= R_1 / 2$ and $0 < R_1 < R_0 < \infty$.
The constants $c_{3,j} > 0$; $j=1,2,3$, are given by
$c_{3,3} = C_1 C_2\cdot C'(R_1) (3 M_1)^3$,
$c_{3,2} = C_1 C_2\cdot C'(R_1) (2 M_1)^2 = C_{2,2}$, and
$c_{3,1} = C_2 M_1$.

The function $u(t) = f_{3,3}(t)\in V$ verifies the elliptic equation
\begin{math}
  (\lambda I + \mathcal{A}) f_{3,3}(t) = f_{2,3}(t)
  \in C_s^{\alpha}\left( \overline{B}_{R_2}^{+}(x_0,0)\right) ,\
  t > 0 .
\end{math}
In Paragraph {\S}\ref{ss:smooth_Hoelder-2} we have shown also
\begin{math}
  u(t) = f_{2,2}(t)
  \in C_s^{\alpha}\left( \overline{B}_{R_2}^{+}(x_0,0)\right) ,
\end{math}
$t > 0$,
together with the norm estimate \eqref{est:f_2,2(t)}.
We apply another local H\"older regularity result from
{\sc P.~M.~N.\ Feehan} and {\sc C.~A.\ Pop}
\cite{Feehan-Pop-14}, Theorem 8.1, Eq.\ (8.4), pp.\ 937--938
(stated in Lemma~\ref{lem-smooth_Hoelder-3},
 Appendix~\ref{s:ellipt_regul}; see also
 \cite{Feehan-17}, Theorem 1.1, Part~2, on pp.\ 2487--2488)
in order to derive
\begin{math}
  u(t) = f_{3,3}(t)\in
  C_s^{2+\alpha}\left( \overline{B}_{R_2'}^{+}(x_0,0)\right) ,\
  t > 0 ,
\end{math}
together with the estimate
\begin{equation*}
%\label{est:f_3,k(t)}
\begin{aligned}
&   \| u(t)\|_{ C_s^{2+\alpha}\left( \overline{B}_{R_2'}^{+}(x_0,0)\right) }
  = \| f_{3,3}(t)
    \|_{ C_s^{2+\alpha}\left( \overline{B}_{R_2'}^{+}(x_0,0)\right) }
\\
& \leq C_3
  \left(
    \| f_{2,3}(t)
    \|_{ C_s^{\alpha}\left( \overline{B}_{R_2}^{+}(x_0,0)\right) }
  + \| f_{3,3}(t)\|_{ C\left( \overline{B}_{R_2}^{+}(x_0,0)\right) }
  \right) ,
\end{aligned}
\end{equation*}
provided $0 < R_2' < R_2 < R_1'= R_1 / 2$
and $0 < R_1 < R_0 < \infty$.
We estimate the right\--hand side by a combination of inequalities
\eqref{est:f_2,2(t)} and \eqref{est:f_2,3(t)}, thus arriving at
\begin{equation}
\label{est:f_3,k(t)}
\begin{aligned}
&   \| u(t)\|_{ C_s^{2+\alpha}\left( \overline{B}_{R_2'}^{+}(x_0,0)\right) }
  = \| f_{3,3}(t)
    \|_{ C_s^{2+\alpha}\left( \overline{B}_{R_2'}^{+}(x_0,0)\right) }
\\
&
\begin{aligned}
  \leq C_3
& \left[
    \left( c_{3,3} t^{-3} + c_{3,2} t^{-2} + c_{3,1} t^{-1}\right)
  + \left( C_{2,2} t^{-2} + C_{2,1} t^{-1} + C_{2,0}\right)
  \right]
    \ee^{\lambda_0 t}\, \| u_0\|_H
\end{aligned}
\\
&
\begin{aligned}
  = \left( C_{3,3} t^{-3} + C_{3,2} t^{-2} + C_{3,1} t^{-1} + C_{3,0}
    \right)
    \ee^{\lambda_0 t}\, \| u_0\|_H
    \quad\mbox{ for all }\, t\in (0,\infty) \,.
\end{aligned}
\\
\end{aligned}
\end{equation}
We have abbreviated the constants $C_{3,j} > 0$; $j=0,1,2,3$, given by
\begin{align*}
  C_{3,3}&= C_3 c_{3,3} = C_1 C_2 C_3\cdot C'(R_1) (3 M_1)^3 \,,
\\
  C_{3,2}&= C_3 (c_{3,2} + C_{2,2})
          = 2 C_1 C_2 C_3\cdot C'(R_1) (2 M_1)^2 \,,
\\
  C_{3,1}&= C_3 (c_{3,1} + C_{2,1})
          = (1 + C_1\cdot C'(R_1)) C_2 C_3 M_1 \,,\quad
  C_{3,0} = C_3 C_{2,0} = C_2 C_3 M_0 \,.
\end{align*}

In particular, we have shown that
\begin{equation*}
  U_{ \overline{B}^{+}_{R_2'}(x_0) }(t)\colon
    u_0\longmapsto         u(t)\vert_{ \overline{B}_{R_2'}^{+}(x_0,0) }
  = ( \ee^{-t\mathcal{A}} u_0 )\vert_{ \overline{B}_{R_2'}^{+}(x_0,0) }
    \colon H\longrightarrow
  C_s^{2+\alpha}\left( \overline{B}_{R_2'}^{+}(x_0,0)\right)
\end{equation*}
is a bounded linear operator with the operator norm
\begin{equation*}
  \| U_{ \overline{B}^{+}_{R_2'}(x_0,0) }(t) \|_{\mathrm{oper}}
  \leq
    \left( C_{3,3} t^{-3} + C_{3,2} t^{-2} + C_{3,1} t^{-1} + C_{3,0}
    \right)
    \ee^{\lambda_0 t} \quad\mbox{ for all }\, t\in (0,\infty) \,.
\end{equation*}
%

%%%%%%%%%%%%%%%%%%%%%%%%%%%%%%%%%%%%%%%%%%%%%%%%%%%%%%%%%%%%%%%%%%%%%%%
%%%%%    Completion of the proof of the main regularity result    %%%%%
%%%%%%%%%%%%%%%%%%%%%%%%%%%%%%%%%%%%%%%%%%%%%%%%%%%%%%%%%%%%%%%%%%%%%%%

\section{Completion of the proof of the main regularity result}
\label{s:proofMain}

In this section we finish the proof of our main regularity result,
Theorem~\ref{thm-Regular}, started in the two previous sections,
Section~\ref{s:smooth_Heston}
and
Section~\ref{s:smooth_Hoelder},
and prove also its Corollary~\ref{cor-Regular}.

\par\vskip 10pt
%%%%%%%%%%%%%%%%%
%\proof
{\it Proof of\/} {\bf Theorem~\ref{thm-Regular}.}$\;$
The regularity statement in {\rm Part~(i)} follows directly from
the results in Section~\ref{s:smooth_Heston},
Ineq.~\eqref{norm:A^k.e^(-tA)}.
The $C^{\infty}$-regularity in {\rm Part~(ii)}
is a (local) interior regularity result for
(local) weak solutions to a locally strictly parabolic equation
established (in a more general setting) in
{\sc A.\ Friedman} \cite[Chapt.~10, Sect.~4]{Friedman-64},
Theorem~11 (p.~302) and its Corollary (p.~303).
The complete proof of {\rm Part~(iii)} has been given in
Section~\ref{s:smooth_Hoelder}.
The radius $R\in (0,\infty)$ stands for the radius
$R_2'\in (0,\infty)$ that appears in Eq.~\eqref{est:f_3,k(t)}.

Finally, we derive {\rm Part~(iv)} from {\rm Part~(iii)} as follows.
The continuity and differentiability of the mapping
\begin{math}
  t\mapsto u(t)\vert_{ \overline{B}^{+}_R(x_0,0) }
\end{math}
from $(0,T)$ to the H\"older space
$C_s^{2+\alpha}( \overline{B}^{+}_R(x_0,0) )$
follow from the respective formulas
\begin{align}
\label{ineq:u(t+tau)-u(t)}
&   \left(
    u(t+\tau) - u(t)
    \right)\vert_{ \overline{B}^{+}_R(x_0,0) }
  = U_{ \overline{B}^{+}_R(x_0,0) }(t) (u(\tau) - u_0)
    \quad\mbox{ and }\quad
\\
\label{ineq:du/dt(t+tau)}
&
\begin{aligned}
&   \left(
    \frac{\partial u}{\partial t} (t+\tau)
    \right)\Big\vert_{ \overline{B}^{+}_R(x_0,0) }
  = \left(
    {}- \mathcal{A} u(t+\tau)
    \right)\big\vert_{ \overline{B}^{+}_R(x_0,0) }
  = U_{ \overline{B}^{+}_R(x_0,0) }(t)
    \left( \frac{\partial u}{\partial t} (\tau) \right)
\\
& = U_{ \overline{B}^{+}_R(x_0,0) }(t)
    \left( {}- \mathcal{A} u(\tau)\right)
  = U_{ \overline{B}^{+}_R(x_0,0) }(t)
    \left( ({}- \mathcal{A})\, \ee^{- \tau\mathcal{A}} u_0\right)
\end{aligned}
\end{align}
for all $t\in (0,T)$ and for all $\tau\in (0,\infty)$ such that
$t + \tau < T$, combined with the locally uniform upper bound
on the operator norm of the bounded linear operator
$U_{ \overline{B}^{+}_R(x_0,0) }(t)$.
Whereas the norm in the H\"older space
$C_s^{2+\alpha}( \overline{B}^{+}_R(x_0,0) )$
of the expression in Eq.~\eqref{ineq:u(t+tau)-u(t)}
above is estimated easily by
the operator norm of $U_{ \overline{B}^{+}_R(x_0,0) }(t)$
from {\rm Part~(iii)}, estimating
the expression in Eq.~\eqref{ineq:du/dt(t+tau)}
requires also the following estimate which follows from inequalities
\eqref{norm:e^(-tA)} and \eqref{norm:A^k.e^(-tA)} ($k=1$),
\begin{align*}
%\label{norm:A^k.e^(-tA)}
  \| ({}- \mathcal{A})\, \ee^{- \tau\mathcal{A}}
  \|_{ \mathcal{L}(H\to H) }
& \leq \| (\lambda I + \mathcal{A})\, \ee^{- \tau\mathcal{A}}
       \|_{ \mathcal{L}(H\to H) }
     + |\lambda |\cdot
       \| \ee^{- \tau\mathcal{A}} \|_{ \mathcal{L}(H\to H) }
\\
& \leq (M_1\tau^{-1} + \lambda_0 M_0)\, \ee^{\lambda_0\tau}
  \quad\mbox{ for all }\, \tau > 0 \,.
\end{align*}
The desired estimate for the norm of
$\frac{\partial u}{\partial t} (t)$
is obtained from that for
$\frac{\partial u}{\partial t} (t+\tau)$
in Eq.~\eqref{ineq:du/dt(t+tau)}
by replacing both, $t$ and $\tau$, by the common value of $t/2$
which means that the sum $t+\tau$ is replaced by $t\in (0,T)$.
%\null\hfill\qed
\qed
%%%%%%%%%%%%%%%%%
\par\vskip 10pt

%\par\vskip 10pt
%%%%%%%%%%%%%%%%%
%\proof
{\it Proof of\/} {\bf Corollary~\ref{cor-Regular}.}$\;$
Let $0 < T\leq +\infty$ and an arbitrary $x_0\in \mathbb{R}$
be fixed.
By the (local) boundary regularity result obtained in
{\rm Part~(iii)} of Theorem~\ref{thm-Regular},
there is a radius $R\in (0,\infty)$ such that,
for every $t\in (0,T)$, we have
\begin{math}
  u(t)\vert_{ \overline{B}^{+}_R(x_0,0) }
  \in C_s^{2+\alpha}( \overline{B}^{+}_R(x_0,0) )
\end{math}
with the norm
\begin{math}
  \| u(t)\|_{ C_s^{2+\alpha}( \overline{B}^{+}_R(x_0,0) ) }
  \leq (c_0' t^{-3} + c_0) \ee^{\lambda_0 t}\, \| u_0\|_H .
\end{math}
Similarly, by {\rm Part~(iv)} of Theorem~\ref{thm-Regular},
for every $t\in (0,T)$, we have
\begin{math}
    \frac{\partial u}{\partial t} (t)\vert_{ \overline{B}^{+}_R(x_0,0) }
  \in C_s^{2+\alpha}( \overline{B}^{+}_R(x_0,0) )
\end{math}
with the norm
\begin{math}
    \left\|
    \frac{\partial u}{\partial t} (t)
    \right\|_{ C_s^{2+\alpha}( \overline{B}^{+}_R(x_0,0) ) }
  \leq (c_1' t^{-4} + c_1 t^{-1})
       \ee^{\lambda_0 t}\cdot \| u_0\|_H
\end{math}
for all $t\in (0,T)$.
We recall that the constants
$c_0, c_0', c_1, c_1'\in (0,\infty)$
do not depend on the choice of $u_0\in H$ or $t\in (0,T)$,
although they may depend on $x_0\in \RR$.

Now let $t_0, T_0\in (0,\infty)$ be arbitrary, but fixed, such that
$0 < t_0 < T_0 < T$ ($\leq +\infty$).
We combine the (local) boundary regularity result from above
with the (local) interior regularity result,
$u\in C^{\infty}( \HH\times (0,\infty) )$
from {\rm Part~(ii)} of Theorem~\ref{thm-Regular},
to conclude that
\begin{math}
  u(t)\vert_{ \overline{B}^{+}_{R_0}(x_0,0) } ,\,
\hfil\break
    \frac{\partial u}{\partial t} (t)\vert_{ \overline{B}^{+}_R(x_0,0) }
  \in C_s^{2+\alpha}( \overline{B}^{+}_{R_0}(x_0,0) )
\end{math}
holds with any finite radius $R_0\in (0,\infty)$ and at any time
$t\in [t_0,T_0]$.
The H\"older norms of $u(t)$ and
\begin{math}
  \frac{\partial u}{\partial t} (t)
\end{math}
satisfy
\begin{equation}
\label{norm:u,du/dt}
    \| u(t)\|_{ C_s^{2+\alpha}( \overline{B}^{+}_{R_0}(x_0,0) ) }
  + \left\|
    \frac{\partial u}{\partial t} (t)
    \right\|_{ C_s^{2+\alpha}( \overline{B}^{+}_R(x_0,0) ) }
 \leq \Gamma
    \quad\mbox{ for every }\, t\in [t_0,T_0] \,.
\end{equation}
The constant
$\Gamma\equiv \Gamma(R_0,t_0,T_0)\in (0,\infty)$
does {\it not\/} depend on the choice of $t\in [t_0,T_0]$.
Moreover, $u$ is a (local) classical solution of the parabolic equation
\begin{math}
%\label{e:Cauchy}
  \frac{\partial u}{\partial t} + \mathcal{A} u = 0
\end{math}
in the strong sense (pointwise) in $\HH\times (0,T)$.
Our next step is to take the limit (as $\xi\to 0+$)
of the function $u(x,\xi,t)$,
its first\--order partial derivatives, and the expressions
$\xi\cdot u_{xx}(x,\xi,t)$,
$\xi\cdot u_{x\xi}(x,\xi,t)$,
$\xi\cdot u_{\xi\xi}(x,\xi,t)$,
for an arbitrary, but fixed pair $(x,t)\in (-R_0,R_0)\times [t_0,T_0]$.
More generally, we fix any pair
$(x^{\ast},t)\in (-R_0,R_0)\times [t_0,T_0]$
which means that
\begin{math}
  P^{\ast} = (x^{\ast},0)\in \partial\HH\cap \overline{B}^{+}_R(x_0,0) .
\end{math}
We will take any point $P = (x,\xi)\in B^{+}_R(x_0,0)$
and calculate the limit (as $P\to P^{\ast}$)
of the functions $u(x,\xi,t)$, $u_t$, etc.\ (as indicated above).

To this end, let us abbreviate the function
\begin{align*}
%\label{e:Cauchy=g(xi)}
&
\begin{aligned}
  g(x,\xi,t)\eqdef
& \;
  \frac{1}{2}\, \sigma\xi\cdot
  \left(
    \frac{\partial^2 u}{\partial x^2}(x,\xi,t)
  + 2\rho\, \frac{\partial^2 u}{\partial x\;\partial\xi}(x,\xi,t)
  + \frac{\partial^2 u}{\partial\xi^2}(x,\xi,t)
  \right)
\\
& {}
  - \xi\cdot
    \left( \frac{1}{2}\, \sigma\cdot
           \frac{\partial u}{\partial x}(x,\xi,t)
         + \kappa\cdot
           \frac{\partial u}{\partial\xi}(x,\xi,t)
    \right)
\end{aligned}
\end{align*}
of $(x,\xi,t)\in \HH\times (0,T)$
and the {\it\bfseries boundary operator\/}, $\mathcal{B}$
(cf.\ \cite[Eq.\ (2.10), p.~8]{AlziaryTak}),
near the boundary $\partial\HH\times (0,T)$,
\begin{align*}
%\label{e:Heston_bc}
&
\begin{aligned}
  (\mathcal{B}u) (x,\xi,t)\eqdef
& \;
    q_r\cdot \frac{\partial u}{\partial x}(x,\xi,t)
  - \kappa\theta_{\sigma}\cdot \frac{\partial u}{\partial\xi}(x,\xi,t)
\end{aligned}
\end{align*}
for $(x,\xi,t)\in \HH\times (0,T)$.
Notice that
\begin{math}
  \mathcal{A}u = \mathcal{B}u - g
\end{math}
holds in $\HH\times (0,T)$.
Hence, the parabolic equation
\begin{math}
%\label{e:Cauchy}
  \frac{\partial u}{\partial t} + \mathcal{A} u = 0
\end{math}
for a (local) classical solution
$u\in C^{2,1}(\HH\times (0,T))$ is equivalent with
\begin{equation}
\label{e:Heston+B}
  \frac{\partial u}{\partial t} + (\mathcal{B}u)(x,\xi,t)
  = g(x,\xi,t)
  \equiv (\mathcal{B} - \mathcal{A}) u
  \quad\mbox{ for }\, (x,\xi,t)\in \HH\times (0,T) \,.
\end{equation}
Fixing any $\xi\in (0,\infty)$
(arbitrarily small for our purpose),
we can easily solve Eq.~\eqref{e:Heston+B}
as a first\--order transport equation for the unknown function
\begin{math}
  (x,t)\mapsto u^{(\xi)}(x,t)\eqdef
               u(x,\xi,t)\colon \RR\times (0,T)\to \RR ,
\end{math}
thus obtaining the following formula,
valid for any $(x,t)\in \RR\times [t_0,T)$:
\begin{align}
\nonumber
& {}
  u^{(\xi)}(x,t)
  = u(x,\xi,t) = u(x - q_r (t-t_0),\, \xi,t_0)
\\
\label{sol:Cauchy:xi>0}
& {}
  + \kappa\theta_{\sigma}
    \int_{t_0}^t \frac{\partial u}{\partial\xi}
    \left( x - q_r (t-s),\, \xi,s\right) \,\mathrm{d}s
  + \int_{t_0}^t g
    \left( x - q_r (t-s),\, \xi,s\right) \,\mathrm{d}s \,.
\end{align}

To complete our proof, let us recall
the (local) boundary regularity results obtained above, in addition to
$u\in C^{\infty}( \HH\times (0,\infty) )$, namely,
\begin{math}
  u(t)\vert_{ \overline{B}^{+}_{R_0}(x_0,0) } ,\,
    \frac{\partial u}{\partial t} (t)\vert_{ \overline{B}^{+}_R(x_0,0) }
  \in C_s^{2+\alpha}( \overline{B}^{+}_{R_0}(x_0,0) )
\end{math}
with any finite radius $R_0\in (0,\infty)$ and at any time
$t\in [t_0,T_0]$.
Moreover, Ineq.~\eqref{norm:u,du/dt}
holds for every $t\in [t_0,T_0]$, with a constant
$\Gamma\equiv \Gamma(R_0,t_0,T_0)\in (0,\infty)$.
Let $x^{*}\in \RR$ be given.
We choose $x_0\in \RR$ arbitrary and $R_0\in (0,\infty)$
large enough, such that
$x_0 - R_0 < x^{*} - q_r T_0 < x^{*} < x_0 + R_0$.
All these inequalities are guaranteed by choosing
$R_0 > |x_0 - x^{*}| + q_r T_0$.
We apply the H\"older regularity from \eqref{norm:u,du/dt}
to all expressions in Eq.~\eqref{e:Heston+B}
in order to conlude that all these expressions belong to
the H\"older space
$C_s^{\alpha}( \overline{B}^{+}_{R_0}(x_0,0) )$,
at any fixed time $t\in [t_0,T_0]$.
In particular, we may take the limit (as $\xi\to 0+$)
of all these expressions in order to conclude that
$g(x,\xi,t)\to g(x^{*},0,t) = 0$ owing to
$P = (x,\xi)\to P^{\ast} = (x^{*},0)\in \partial\HH$.
Here, the limits of both first\--order partial derivatives
$u_x$ and $u_{\xi}$ as $\xi\to 0+$ exist and are bounded by
\eqref{norm:u,du/dt} and the definition of the H\"older space
$C_s^{2+\alpha}( \overline{B}^{+}_{R_0}(x_0,0) )$,
cf.\ Eq.~\eqref{norm:2+alpha},
whereas the limits of all expressions containing
the second\--order partial derivatives,
$\xi\cdot u_{xx}$, $\xi\cdot u_{x\xi}$, and $\xi\cdot u_{\xi\xi}$,
vanish as $\xi\to 0+$, by
{\sc P.~M.~N.\ Feehan} and {\sc C.~A.\ Pop}
\cite{Feehan-Pop-13}, Lemma 3.1, Eq.\ (3.1), on p.~4409
(see also
 {\sc P.\ Daskalopoulos} and {\sc R.\ Hamilton}
 \cite{Daska-Hamil-98}, Prop.\ I.12.1 on p.~940).
We complete our proof by applying these limits to
Eqs.\ \eqref{e:Heston+B} and \eqref{sol:Cauchy:xi>0},
thus arriving at
Eqs.\ \eqref{e:Cauchy:xi=0} and \eqref{sol:Cauchy:xi=0},
as desired.
%\null\hfill\qed
\qed
%%%%%%%%%%%%%%%%%
\par\vskip 10pt

%%%%%%%%%%%%%%%%%%%%%%%%%%%%%%%%%%%%%%%%%%%%%%%%%%%%%%%%%%%%%%%%%%%%%%%
%%%%%    A maximum principle and growth at low and high volatil.    %%%
%%%%%%%%%%%%%%%%%%%%%%%%%%%%%%%%%%%%%%%%%%%%%%%%%%%%%%%%%%%%%%%%%%%%%%%

\section{A maximum principle and growth at low and high volatilities}
\label{s:MaxPrinciple}

According to a classical result on the weak maximum principle for
a uniformly parabolic Cauchy problem in $\RR^N\times (0,T)$,
see e.g.\ {\sc A.\ Friedman}
\cite[Chapt.~2, Sect.~4, Theorem~9, p.~43]{Friedman-64},
the weak maximum principle is valid under ``very weak'' restrictions
on the growth of a strong solution $u(x,t)$ as $|x|\to \infty$,
$(x,t)\in \RR^N\times (0,T)$.
Consequently, one may speak of practicaly
{\it\bfseries no boundary conditions\/} being imposed on
the strong solution $u(x,t)$ as $|x|\to \infty$,
at least in contrast with classical boundary conditions
of Dirichlet, Neumann, or oblique derivative (Robin) types.
Nevertheless, thanks to the weak maximum principle,
the uniqueness of any strong solution to the Cauchy problem
with prescribed initial data is still guaranteed.

Now we are ready to prove our Theorem~\ref{thm-MaxPrinciple}.

\par\vskip 10pt
%%%%%%%%%%%%%%%%%
%\proof
{\it Proof of\/} {\bf Theorem~\ref{thm-MaxPrinciple}.}$\;$
Let us recall that $\gamma_0\in (0,\infty)$
is an arbitrary constant, as large as needed, the constants
$\beta_0, \mu_0\in (0,\infty)$ satisfy inequalities
\eqref{ineq:beta,mu},
and the function $\mathfrak{h}_0$ is defined in \eqref{def:h_0}.

We will compare the function $u\colon \HH\times (0,T)\to \RR$
to the smooth function $h$ defined as follows:
\begin{equation}
\label{def:h}
\begin{aligned}
  h(x,\xi,t) \eqdef
    \exp\left(
    \genfrac{}{}{}0{ \gamma_1 (1 + x^2)^{1/2} + \mu_1\xi
                   - (\beta_1 - 1)\, \ln\xi }{1 - \omega t}
  + \nu t
        \right)
\\
\end{aligned}
\end{equation}
for $(x,\xi,t)\in \HH\times (0,T)$,
where $\beta_1\geq 1$, $\gamma_1 > 0$, $\mu_1 > 0$, $\nu\geq 0$, and
$\omega > 0$ are suitable positive constants to be specified later
in the proof.
Clearly, $h(x,\xi,t)^{-1}$ replaces the weight function
$\mathfrak{w}\colon \mathbb{H}\to (0,\infty)$ defined in Eq.~\eqref{def:w}.

We calculate the partial derivatives of $h(x,\xi,t)$
at $(x,\xi)\in \HH$ and $0 < t < T$:
\begin{align*}
%\label{e:dh/dt}
    h^{-1}\, \frac{\partial h}{\partial t}
& = \frac{\omega}{ (1 - \omega t)^2 }
    \left[ \gamma_1 (1 + x^2)^{1/2} + \mu_1\xi - (\beta_1 - 1)\, \ln\xi
    \right]
  + \nu \,,
\\
%\label{e:dh/dx}
    h^{-1}\, \frac{\partial h}{\partial x}
& = \frac{\gamma_1}{1 - \omega t}\, \frac{x}{ (1 + x^2)^{1/2} } \,,
  \quad
%\label{e:dh/d_xi}
    h^{-1}\, \frac{\partial h}{\partial\xi}
  = \frac{1}{1 - \omega t}
    \left( \mu_1 - \frac{\beta_1 - 1}{\xi} \right) \,.
\end{align*}
Similarly, we calculate the second\--order partial derivatives:
\begin{align*}
%\label{e:d^2_h/dx^2}
&
\begin{aligned}
    h^{-1}\, \frac{\partial^2 h}{\partial x^2}
& = \genfrac{(}{)}{}0{\gamma_1}{1 - \omega t}^2\, \frac{x^2}{1 + x^2}
  + \frac{\gamma_1}{1 - \omega t}
    \left[ \frac{1}{ (1 + x^2)^{1/2} }
         - \frac{x^2}{ (1 + x^2)^{3/2} }
    \right]
\\
& = \frac{\gamma_1^2}{ (1 - \omega t)^2 }
    \left( 1 - \frac{1}{1 + x^2}\right)
  + \frac{\gamma_1}{1 - \omega t}\, \frac{1}{ (1 + x^2)^{3/2} } \,,
\end{aligned}
\\
%\label{e:d^2_h/d_xi^2}
&
\begin{aligned}
    h^{-1}\, \frac{\partial^2 h}{\partial\xi^2}
& = \frac{1}{ (1 - \omega t)^2 }
    \left( \mu_1 - \frac{\beta_1 - 1}{\xi} \right)^2
  + \frac{\beta_1 - 1}{1 - \omega t}\, \frac{1}{\xi^2} \,,
\end{aligned}
\\
%\label{e:d^2_h/dx.d_xi}
&
\begin{aligned}
    h^{-1}\, \frac{\partial^2 h}{\partial x\, \partial\xi}
& = \frac{\gamma_1}{ (1 - \omega t)^2 }
    \left( \mu_1 - \frac{\beta_1 - 1}{\xi} \right)
    \frac{x}{ (1 + x^2)^{1/2} } \,.
\end{aligned}
\end{align*}
We plug these partial derivatives of $h$ into formula
\eqref{e:Heston-oper} to calculate
\begin{align}
\label{e:dh/dt+Heston-oper}
& {}
  - h^{-1}
  \left( \frac{\partial h}{\partial t} + \mathcal{A} h\right)
  = \frac{1}{2}\, \sigma\xi
  \Biggl[
    \frac{\gamma_1^2}{ (1 - \omega t)^2 }
    \left( 1 - \frac{1}{1 + x^2}\right)
  + \frac{\gamma_1}{1 - \omega t}\, \frac{1}{ (1 + x^2)^{3/2} }
\\
\nonumber
& {}
  + \frac{2\rho\gamma_1}{ (1 - \omega t)^2 }
    \left( \mu_1 - \frac{\beta_1 - 1}{\xi} \right)
    \frac{x}{ (1 + x^2)^{1/2} }
  + \frac{1}{ (1 - \omega t)^2 }
    \left( \mu_1 - \frac{\beta_1 - 1}{\xi} \right)^2
  + \frac{\beta_1 - 1}{1 - \omega t}\, \frac{1}{\xi^2}
  \Biggr]
\\
\nonumber
& {}
  - \left( q_r + \genfrac{}{}{}1{1}{2} \sigma\xi\right)
    \frac{\gamma_1}{1 - \omega t}\, \frac{x}{ (1 + x^2)^{1/2} }
  + \kappa (\theta_{\sigma} - \xi)
    \frac{1}{1 - \omega t}
    \left( \mu_1 - \frac{\beta_1 - 1}{\xi} \right)
\\
\nonumber
& {}
  - \frac{\omega}{ (1 - \omega t)^2 }
    \left[ \gamma_1 (1 + x^2)^{1/2} + \mu_1\xi - (\beta_1 - 1)\, \ln\xi
    \right]
  - \nu
  \equiv J_1\xi + J_0 + J_{-1}\xi^{-1} \,,
\end{align}
where we recall
$\theta_{\sigma} = \theta / \sigma$ and abbreviate
\begin{align}
\label{e:J_1.xi}
&
\begin{aligned}
  J_1\eqdef \frac{1}{2}\, \frac{\sigma}{1 - \omega t}
  \Biggl[
& {}
    \frac{\gamma_1^2}{1 - \omega t}
    \left( 1 - \frac{1}{1 + x^2}\right)
  + \frac{\gamma_1}{ (1 + x^2)^{3/2} }
\\
& {}
  + \frac{2\rho\gamma_1\mu_1}{1 - \omega t}\,
    \frac{x}{ (1 + x^2)^{1/2} }
  + \frac{\mu_1^2}{1 - \omega t}
  - \gamma_1\, \frac{x}{ (1 + x^2)^{1/2} }
  \Biggr]
\end{aligned}
\\
\nonumber
& {}
    \phantom{J_1\eqdef}
  - \frac{\kappa\mu_1}{1 - \omega t}
  - \frac{\omega}{ (1 - \omega t)^2 }
    \left( \mu_1 - (\beta_1 - 1)\frac{\ln\xi}{\xi} \right) \,,
\\
%\label{e:J_0.xi}
\nonumber
&
\begin{aligned}
  J_0\eqdef
& {}
  - \frac{ \sigma\rho\gamma_1 (\beta_1 - 1) }{ (1 - \omega t)^2 }\,
    \frac{x}{ (1 + x^2)^{1/2} }
  - \frac{ \sigma\mu_1 (\beta_1 - 1) }{ (1 - \omega t)^2 }
  - \frac{ q_r\gamma_1 }{1 - \omega t}\,
    \frac{x}{ (1 + x^2)^{1/2} }
\\
& {}
  + \frac{ \kappa [ \theta_{\sigma}\mu_1 + (\beta_1 - 1) ] }{1 - \omega t}
  - \frac{\omega\gamma_1}{ (1 - \omega t)^2 }\, (1 + x^2)^{1/2}
  - \nu \,, \quad\mbox{ and }
\end{aligned}
\\
%\label{e:J_(-1).xi}
\nonumber
&
\begin{aligned}
  J_{-1}
& {}
  \eqdef \frac{1}{2}\, \frac{\sigma (\beta_1 - 1)}{1 - \omega t}
  \left( \frac{\beta_1 - 1}{1 - \omega t} + 1\right)
  - \frac{ \kappa\theta_{\sigma} (\beta_1 - 1) }{1 - \omega t}
\\
& {}
  = \frac{\beta_1 - 1}{1 - \omega t}
    \left( \frac{1}{2}\, \frac{\sigma (\beta_1 - 1)}{1 - \omega t}
         + \frac{1}{2}\, \sigma - \kappa\theta_{\sigma}
    \right) \,.
\end{aligned}
\end{align}
Our assumption on $\beta_0\in (0,\infty)$ in \eqref{ineq:beta,mu}
allows us to find a constant $\tau\in (0,1)$ small enough, such that
\begin{equation*}
%\label{ineq:tau}
  (1\leq {})\hspace{7pt}
  \beta_0 < 1 + (1 - \tau)
            \left( \frac{2\kappa\theta}{\sigma^2} - 1\right)
          < \frac{2\kappa\theta}{\sigma^2}
\end{equation*}
or, equivalently,
\begin{equation*}
%\label{ineq:tau}
  0 < \tau < \tau_0\eqdef
      \left( \frac{2\kappa\theta}{\sigma^2} - \beta_0\right)
    \bigg\slash
      \left( \frac{2\kappa\theta}{\sigma^2} - 1\right)
    \hspace{7pt} ({}\leq 1) \,.
\end{equation*}
From now on we restrict ourselves to the time interval
$0 < t\leq T_{\omega}\eqdef \tau / \omega$ with $\omega > 0$
to be determined as follows.

We begin with estimating the last expression, $J_{-1}$.
We fix any $\gamma_1\in (0,\infty)$ such that $\gamma_1 > \gamma_0$.
Recalling {\it\bfseries Feller's condition\/} \eqref{e:Feller}
and the first inequality in \eqref{ineq:beta,mu},
let us choose $\beta_1\in [1,\infty)$ such that
\begin{equation}
\label{ineq:beta_1>1}
  (1\leq {})\hspace{7pt}
    \beta_0 < \beta_1\leq
        1 + (1 - \tau)
            \left( \frac{2\kappa\theta}{\sigma^2} - 1\right)
          < \frac{2\kappa\theta}{\sigma^2} \,.
\end{equation}
This choice guarantees the following inequality, whenever
$0 < t\leq T_{\omega}$,
\begin{equation}
\label{ineq:J_(-1)}
\begin{aligned}
  J_{-1}
& {}\leq
    \frac{\beta_1 - 1}{1 - \omega t}
    \left( \frac{1}{2}\, \frac{\sigma (\beta_1 - 1)}{1 - \tau}
         + \frac{1}{2}\, \sigma - \kappa\theta_{\sigma}
    \right)
\\
& {}
  = \frac{ \sigma (\beta_1 - 1) }{ 2(1 - \omega t) (1 - \tau) }
    \left[ \beta_1 - 1
    - (1 - \tau)\left( \frac{2\kappa\theta}{\sigma^2} - 1\right)
    \right]\leq 0 \,.
\end{aligned}
\end{equation}

We fix a suitable constant $\mu_1 > 0$ in the first expression, $J_1$,
as follows:
\begin{equation}
\label{ineq:mu>beta-1}
  (0\leq \beta_0 - 1 < {})\hspace{7pt}
    \max\{ \beta_1 - 1,\, \mu_0\} < \mu_1 < \infty \,.
\end{equation}

This choice, combined with the standard inequality
$\ln\xi\leq \xi - 1$ for all $\xi > 0$, guarantees
the following estimate for the expression in the parentheses
of the last summand in $J_1$, Eq.~\eqref{e:J_1.xi},
\begin{equation*}
%\label{e:mu-(beta-1)}
\begin{aligned}
& {}
  - \frac{\omega}{ (1 - \omega t)^2 }
    \left( \mu_1 - (\beta_1 - 1)\, \frac{\ln\xi}{\xi} \right)
  \leq
{}- \frac{\omega}{ (1 - \omega t)^2 }
    \left( \mu_1 - (\beta_1 - 1)\, \frac{\xi - 1}{\xi} \right)
\\
  =
& {}
  - \frac{\omega}{ (1 - \omega t)^2 }
    \left( \mu_1 - (\beta_1 - 1) + \frac{\beta_1 - 1}{\xi} \right)
  \leq
{}- \frac{\omega}{ (1 - \omega t)^2 }\, [ \mu_1 - (\beta_1 - 1) ] \,.
\end{aligned}
\end{equation*}
We apply this inequality and the trivial relation
$|x|\leq (1 + x^2)^{1/2}$ for all $x\in \RR$ to estimate $J_1$,
whenever $0 < t\leq T_{\omega}$:
\begin{equation*}
%\label{est:J_1}
\begin{aligned}
  J_1
& {}\leq \frac{1}{2}\, \frac{\sigma}{1 - \omega t}
  \Biggl(
    \frac{\gamma_1^2}{1 - \tau} + \gamma_1
  + \frac{ 2\, |\rho|\, \gamma_1\mu_1}{1 - \tau}
  + \frac{\mu_1^2}{1 - \tau} + \gamma_1
  \Biggr)
  - \kappa\mu_1 - \omega\, [ \mu_1 - (\beta_1 - 1) ]
\\
& {}\leq \frac{\sigma}{1 - \tau}
  \Biggl(
    \frac{ \gamma_1^2 + 2\, |\rho|\, \gamma_1\mu_1 + \mu_1^2 }%
         {2(1 - \tau)}
  + \gamma_1
  \Biggr)
  - \kappa\mu_1 - \omega\, [ \mu_1 - (\beta_1 - 1) ] \,.
\end{aligned}
\end{equation*}
Recall that the correlation coefficient $\rho$ satisfies
$\rho\in (-1,1)$.
All constants $\beta_1\geq 1$, $\gamma_1 > 0$, and $\mu_1 > 0$
having been fixed, such that all inequalities
\eqref{e:Feller}, \eqref{ineq:beta_1>1}, and \eqref{ineq:mu>beta-1}
are valid, we now choose $\omega\in (0,\infty)$ large enough to guarantee
$\omega\geq \tau / T$ and also
\begin{equation}
\label{ineq:J_1}
  J_1\leq \frac{\sigma}{1 - \tau}
  \Biggl( \frac{ (\gamma_1 + \mu_1)^2 }{2(1 - \tau)} + \gamma_1
  \Biggr)
  - \kappa\mu_1 - \omega\, [ \mu_1 - (\beta_1 - 1) ]
  \leq 0
\end{equation}
whenever $0 < t\leq \tau / \omega$ $({}= T_{\omega}\leq T)$.

The constant $\nu$ appears in the expression $J_0$ only;
we take $\nu\in \RR_+ = [0,\infty)$ arbitrary.
Since $|x|\leq (1 + x^2)^{1/2}$ holds for every $x\in \RR$,
we can choose $\omega\in (0,\infty)$ even greater than above to obtain also
\begin{equation}
\label{ineq:J_0}
\begin{aligned}
  J_0 + \nu
& {}
  \leq \frac{ \sigma\rho\gamma_1 (\beta_1 - 1) }{ (1 - \tau)^2 }
  - \sigma\mu_1 (\beta_1 - 1) + \frac{ q_r\gamma_1 }{1 - \tau}
  + \frac{ \kappa\, [ \theta_{\sigma}\mu_1 + (\beta_1 - 1) ] }{1 - \tau}
  - \omega\gamma_1
\\
& {}
  = \sigma (\beta_1 - 1)
    \left( \frac{\rho\gamma_1}{ (1 - \tau)^2 } - \mu_1\right)
  + \frac{ q_r\gamma_1
         + \kappa\, [ \theta_{\sigma}\mu_1 + (\beta_1 - 1) ] }{1 - \tau}
  - \omega\gamma_1
  \leq 0
\end{aligned}
\end{equation}
whenever $0 < t\leq T_{\omega} = \tau / \omega$ $({}\leq T)$.
In other words, the constant $\omega\in (0,\infty)$
must be large enough in order to obey all three inequalities,
$\omega\geq \tau / T$, \eqref{ineq:J_1}, and \eqref{ineq:J_0}.

We remark that, in the works by
{\sc P.\ Daskalopoulos} and {\sc P.~M.~N.\ Feehan} \cite{Daska-Feehan-14}
and \cite[Sect.~2, p.~5048]{Daska-Feehan-16},
the constants $\beta$ and $\mu$ are chosen to be
\begin{math}
%\label{ineq:beta_1>1}
  \beta = {2\kappa\theta}/{\sigma^2} > 1
\end{math}
and
\begin{math}
%\label{ineq:mu>beta-1}
  \mu = {2\kappa}/{\sigma^2} = {\beta}/{\theta} .
\end{math}

Finally, we apply inequalities
\eqref{ineq:J_(-1)}, \eqref{ineq:J_1}, and \eqref{ineq:J_0}
to Eq.~\eqref{e:dh/dt+Heston-oper}
to infer that
\begin{equation}
\label{ineq:dh/dt+Heston-oper}
\begin{aligned}
& {}- h^{-1}
  \left( \frac{\partial h}{\partial t} + \mathcal{A} h\right)
  \equiv J_1\xi + J_0 + J_{-1}\xi^{-1}
  \leq -\nu\leq 0
\\
& \quad\mbox{ is valid for all }\, (x,\xi)\in \HH
  \;\mbox{ and for all }\, 0 < t\leq T_{\omega} \,.
\end{aligned}
\end{equation}

In order to obtain a weak maximum principle for a strong solution
$u\colon \HH\times (0,T)\to \RR$
of the initial value problem \eqref{e:Cauchy}, such that
$u(x,\xi,t)\leq \mathrm{const}\cdot h(x,\xi,t)$ for all
$(x,\xi)\in \HH$ and $t\in (0,T)$,
from the parabolic equation in the Cauchy problem \eqref{e:Cauchy}
we derive an analogous parabolic equation for the ratio
\begin{math}
  w(x,\xi,t)\eqdef u(x,\xi,t) / h(x,\xi,t)\leq \mathrm{const}
  < \infty .
\end{math}
Using $u = wh$ we have
\begin{align*}
%\label{e:d_wh/dt}
    \frac{\partial u}{\partial t}
& = \frac{\partial w}{\partial t}\, h
  + w\, \frac{\partial h}{\partial t} \,,
  \quad
%\label{e:d_wh/dx}
    \frac{\partial u}{\partial x}
  = \frac{\partial w}{\partial x}\, h
  + w\, \frac{\partial h}{\partial x} \,,
  \quad
%\label{e:d_wh/d_xi}
    \frac{\partial u}{\partial\xi}
  = \frac{\partial w}{\partial\xi}\, h
  + w\, \frac{\partial h}{\partial\xi} \,,
\end{align*}
and similarly
\begin{align*}
%\label{e:d^2_wh/dx^2}
    \frac{\partial^2 u}{\partial x^2}
& {}
  = \frac{\partial^2 w}{\partial x^2}\, h
  + 2\frac{\partial w}{\partial x}\, \frac{\partial h}{\partial x}
  + w\, \frac{\partial^2 h}{\partial x^2} \,,
   \quad
%\label{e:d^2_wh/d_xi^2}
    \frac{\partial^2 u}{\partial\xi^2}
  = \frac{\partial^2 w}{\partial\xi^2}\, h
  + 2\frac{\partial w}{\partial\xi}\, \frac{\partial h}{\partial\xi}
  + w\, \frac{\partial^2 h}{\partial\xi^2} \,,
\\
%\label{e:d^2_wh/dx.d_xi}
    \frac{\partial^2 u}{\partial x\, \partial\xi}
& = \frac{\partial^2 w}{\partial x\, \partial\xi}\, h
  + \frac{\partial w}{\partial x}\, \frac{\partial h}{\partial\xi}
  + \frac{\partial w}{\partial\xi}\, \frac{\partial h}{\partial x}
  + w\, \frac{\partial^2 h}{\partial x\, \partial\xi} \,.
\end{align*}
We plug these partial derivatives of $u$ into formula
\eqref{e:Heston-oper} to calculate
\begin{align*}
%\label{e:du/dt+Heston-oper}
& {}- h^{-1}
  \left( \frac{\partial u}{\partial t} + \mathcal{A} u\right)
  = {}- \frac{\partial w}{\partial t} - \mathcal{A} w
  - h^{-1}
  \left( \frac{\partial h}{\partial t} + \mathcal{A} h\right)
  \cdot w(x,\xi,t)
\\
& {}
  + \sigma\xi
  \left[
    \frac{\partial w}{\partial x}\, \frac{\partial h}{\partial x}
  + \rho
    \left(
    \frac{\partial w}{\partial x}\, \frac{\partial h}{\partial\xi}
  + \frac{\partial w}{\partial\xi}\, \frac{\partial h}{\partial x}
    \right)
  + \frac{\partial w}{\partial\xi}\, \frac{\partial h}{\partial\xi}
  \right] \,,
\end{align*}
or equivalently, for all $(x,\xi,t)\in \HH\times (0,T)$,
\begin{align}
\label{e:du/dt+Heston-oper}
&
\begin{aligned}
  {}- h^{-1}
  \left( \frac{\partial u}{\partial t} + \mathcal{A} u\right)
& {}
  = {}- \frac{\partial w}{\partial t} - \mathcal{A} w
  - h^{-1}
  \left( \frac{\partial h}{\partial t} + \mathcal{A} h\right)
  \cdot w(x,\xi,t)
\\
& {}
  + \sigma\xi
  \left( \frac{\partial h}{\partial x}
  + \rho\, \frac{\partial h}{\partial\xi}
  \right)
  \cdot \frac{\partial w}{\partial x}
  + \sigma\xi
  \left( \frac{\partial h}{\partial\xi}
  + \rho\, \frac{\partial h}{\partial x}
  \right)
  \cdot \frac{\partial w}{\partial\xi} \,.
\end{aligned}
\end{align}
We recall that the multiplicative coefficient at $w(x,\xi,t)$
is ${}\leq -\nu\leq 0$, by Ineq.~\eqref{ineq:dh/dt+Heston-oper}.

Recalling formula \eqref{def:h} for $h$, the ratio
\begin{align*}
\nonumber
  \widetilde{h}(x,\xi,t)\eqdef
    \frac{ \mathfrak{h}_0(x,\xi) }{h(x,\xi,t)}
%\label{def:tilde_h}
& {}
  = \exp\left(
    {}- (\gamma_1 - \gamma_0) (1 + x^2)^{1/2} - (\mu_1 - \mu_0)\xi
      + (\beta_1 - \beta_0)\, \ln\xi
        \right)
\\
\nonumber
& {}
  \times \exp
        \left( {}- \genfrac{}{}{}0{\omega t}{1 - \omega t}
        \left[
        \gamma_1\, (1 + x^2)^{1/2} + \mu_1\xi - (\beta_1 - 1)\, \ln\xi
        \right] - \nu t
        \right)
\\
\nonumber
& {}
  \leq \exp\left(
    {}- (\gamma_1 - \gamma_0) (1 + x^2)^{1/2} - (\mu_1 - \mu_0)\xi
      + (\beta_1 - \beta_0)\, \ln\xi 
        \right)
\end{align*}
has the following asymptotic behavior,
for $(x,\xi)\in \HH$ and
$0 < t < T_{\omega} = \tau / \omega$ $({}\leq T)$:
\begin{align*}
%\label{lim:tilde_h_x}
& \lim_{|x|\to \infty}
  \sup_{ (\xi,t)\in (0,\infty)\times (0,T_{\omega}) }
    \widetilde{h}(x,\xi,t) {}= 0 \,,
\\
%\label{lim:tilde_h_xi=0,infty}
& \lim_{\xi\to 0+}
  \sup_{ (x,t)\in \RR\times (0,T_{\omega}) }
    \widetilde{h}(x,\xi,t) {}= 0 \,,\quad
  \lim_{\xi\to +\infty}
  \sup_{ (x,t)\in \RR\times (0,T_{\omega}) }
    \widetilde{h}(x,\xi,t) {}= 0 \,.
\end{align*}
These limits follow from inequalities
\eqref{ineq:beta_1>1} and \eqref{ineq:mu>beta-1} combined with
$\omega\geq \tau / T$.
From the limits above we derive analogous results for the ratio
\begin{equation*}
%\label{e:w=u/h}
\begin{aligned}
  w(x,\xi,t)\eqdef \frac{u(x,\xi,t)}{h(x,\xi,t)}
  = \frac{u(x,\xi,t)}{ \mathfrak{h}_0(x,\xi) }\,
    \frac{ \mathfrak{h}_0(x,\xi) }{h(x,\xi,t)}
  \leq C\cdot \widetilde{h}(x,\xi,t)
\\
  \quad\mbox{ for }\, (x,\xi)\in \HH \,\mbox{ and }\,
    0 < t < T_{\omega}\; ({}\leq T) ,
\end{aligned}
\end{equation*}
namely,
\begin{align}
\label{lim:w_x}
& \limsup_{|x|\to \infty}
  \sup_{ (\xi,t)\in (0,\infty)\times (0,T_{\omega}) } w(x,\xi,t)
  {}\leq 0 \,,
\\
\label{lim:w_xi=0,infty}
& \limsup_{\xi\to 0+}
  \sup_{ (x,t)\in \RR\times (0,T_{\omega}) } w(x,\xi,t)
  {}\leq 0 \,,\quad
  \limsup_{\xi\to +\infty}
  \sup_{ (x,t)\in \RR\times (0,T_{\omega}) } w(x,\xi,t)
  {}\leq 0 \,.
\end{align}

In order to complete our proof, we recall the parabolic equation
\eqref{e:du/dt+Heston-oper} for $w$ with the right\--hand side
${}\geq 0$, by Ineq.~\eqref{ineq:Cauchy},
and the multiplicative coefficient ${}\leq -\nu\leq 0$ at $w(x,\xi,t)$,
by Ineq.~\eqref{ineq:dh/dt+Heston-oper},
or equivalently, for all $(x,\xi,t)\in \HH\times (0,T)$,
\begin{align*}
%\label{e:du/dt+Heston-oper}
  \frac{\partial w}{\partial t} + \mathcal{A} w
& {}
  - \sigma\xi
  \left( \frac{\partial h}{\partial x}
  + \rho\, \frac{\partial h}{\partial\xi}
  \right)
  \cdot \frac{\partial w}{\partial x}
  - \sigma\xi
  \left( \frac{\partial h}{\partial\xi}
  + \rho\, \frac{\partial h}{\partial x}
  \right)
  \cdot \frac{\partial w}{\partial\xi}
\\
& = h^{-1}
  \left( \frac{\partial u}{\partial t} + \mathcal{A} u\right)
  - h^{-1}
  \left( \frac{\partial h}{\partial t} + \mathcal{A} h\right)
  \cdot w(x,\xi,t)
  \leq 0 \,.
\end{align*}
Taking advantage of the initial condition \eqref{cond:Cauchy:t=0}
at time $t=0$, which is equivalent with
$w(x,\xi,0)$ $\leq 0$ for all $(x,\xi)\in \HH$, in addition to
the boundary behavior \eqref{lim:w_x} and \eqref{lim:w_xi=0,infty},
we may apply the weak maximum principle from {\sc A.~Friedman}
\cite[Chapt.~2, Sect.~4, Lemma~5, p.~43]{Friedman-64}
to conlude that $w(x,\xi,t)\leq 0$ holds for all
$(x,\xi)\in \HH$ at all times $t\in [0,T_{\omega})$.
We may apply this result in any subinterval
$[t_0,\, t_0 + T_{\omega}) \subset [0,T)$
of length $0 < T_{\omega} = \tau / \omega\leq T$
to extend the weak maximum principle in $\HH\times (0,T_{\omega})$
to the entire domain $\HH\times (0,T)$.
The corresponding result for $u = wh$ now follows exactly as in
\cite[Chapt.~2, Sect.~4, Theorem~9, p.~43]{Friedman-64}.

The uniqueness for the Cauchy problem \eqref{e:Cauchy}
follows from the weak maximum principle exactly as in
\cite[Chapt.~2, Sect.~4, Theorem~10, p.~44]{Friedman-64}.

Theorem~\ref{thm-MaxPrinciple} is proved.
%\null\hfill\qed
\qed
%%%%%%%%%%%%%%%%%
\par\vskip 10pt

It remains to give

\par\vskip 10pt
%%%%%%%%%%%%%%%%%
%\proof
{\it Proof of\/} {\bf Corollary~\ref{cor-MaxPrinciple}.}$\;$
Our strategy of the proof is to verify that
the weak maximum principle in Theorem~\ref{thm-MaxPrinciple}
can be applied to both functions
$W_{\pm}(x,\xi,t) = {}- U(x,\xi,t)\pm u(x,\xi,t)$
for $(x,\xi,t)\in \HH\times [0,T)$.
In Theorem~\ref{thm-MaxPrinciple},
we take $T\in (0,\infty)$ arbitrarily large, but finite.

We begin with the growth restriction \eqref{cond:|u|<h_0}.
The strong solution to
the homogeneous Cauchy problem \eqref{e:Cauchy} with $f\equiv 0$,
\begin{math}
  u\in C^0(\HH\times [0,T))\cap C^{2,1}(\HH\times (0,T)) ,
\end{math}
obeys this restriction by hypothesis.
Hence, it remains to verify that so does the function
$U\colon \HH\times (0,\infty)\to \RR$, that is to say,
\begin{align*}
  (0\leq {})\hspace{7pt}
  U(x,\xi,t)
& {}
  = \ee^{r_0 t}\left( K_1\, \ee^{x + \varpi\xi} + K_0\right)
  \leq C\, \ee^{r_0 T}\cdot \mathfrak{h}_0(x,\xi)
\\
& {}
  = C\, \ee^{r_0 T}\cdot \xi^{- (\beta_0 - 1)}\,
    \exp\left[ \gamma_0 (1 + x^2)^{1/2} + \mu_0\xi\right]
\end{align*}
holds for all $(x,\xi,t)\in \HH\times (0,\infty)$,
with some constant $C\in (0,\infty)$; see
Eqs.\ \eqref{def:h_0} and \eqref{ineq:Cauchy:u_0}.
We recall from the hypotheses in Theorem~\ref{thm-MaxPrinciple}
(the weak maximum principle)
that the constant $\gamma_0\in [1,\infty)$ is arbitrary and
$\beta_0, \mu_0\in (0,\infty)$ satisfy inequalities \eqref{ineq:beta,mu}.
Consequently, we have to choose $\gamma_0\geq 1$,
as we have already done in the hypotheses, and
$0\leq \varpi < \mu_0$ assumed in Ineq.~\eqref{ineq:varpi}, as well.
We conclude that the functions
$W_{\pm}\colon \HH\times [0,T)\to \RR$
obey the growth restriction \eqref{cond:|u|<h_0}.

Furthermore, the restriction at the initial time $t=0$
in Ineq.~\eqref{ineq:Cauchy:u_0} guarantees
$W_{\pm}(x,\xi,0)$ $\leq 0$ for all $(x,\xi)\in \HH$.

Thus, conditions
\eqref{cond:Cauchy} and \eqref{cond:Cauchy:t=0}
having been verified above, only
Ineq.~\eqref{ineq:Cauchy} for $W_{\pm}$ in place of~$u$
remains to be proved.
Notice that $\frac{\partial U}{\partial t} = r_0 U$ and
\begin{math}
%\label{e:Cauchy}
  \frac{\partial u}{\partial t} + \mathcal{A} u = 0
\end{math}
in the strong sense (pointwise) in $\HH\times (0,T)$, thanks to
$u\in C^{2,1}(\HH\times (0,T))$.

The first and second partial derivatives of $U$ are
\begin{align*}
& U_x = K_1\, \ee^{r_0 t}\cdot \ee^{x + \varpi\xi} \,,\quad
  U_{\xi} = \varpi K_1\, \ee^{r_0 t}\cdot \ee^{x + \varpi\xi} \,,\quad
\\
& U_{xx} = K_1\, \ee^{r_0 t}\cdot \ee^{x + \varpi\xi} \,,\quad
  U_{x\xi} = \varpi K_1\, \ee^{r_0 t}\cdot \ee^{x + \varpi\xi} \,,\quad
  U_{\xi\xi} = \varpi^2 K_1\, \ee^{r_0 t}\cdot \ee^{x + \varpi\xi} \,.
\end{align*}
We insert them into the Heston operator \eqref{e:Heston-oper},
$\mathcal{A}$,
\begin{align*}
%\label{e:Heston-oper}
&
\begin{aligned}
  \ee^{- r_0 t}
  \left( \frac{\partial U}{\partial t} + (\mathcal{A}U)(x,\xi)
  \right)
& {}
  = r_0 K_0
  + K_1\cdot \ee^{x + \varpi\xi}
    \left[ r_0 - \frac{1}{2}\, \sigma\xi\cdot
  \left( 1 + 2\rho\varpi + \varpi^2\right)
    \right]
\\
& {}
  + \left( q_r + \genfrac{}{}{}1{1}{2} \sigma\xi\right)
    \cdot K_1\cdot \ee^{x + \varpi\xi}
  - \kappa (\theta_{\sigma} - \xi)
    \cdot \varpi K_1\cdot \ee^{x + \varpi\xi}
\end{aligned}
\\
& {}
  = r_0 K_0
  + K_1\cdot \ee^{x + \varpi\xi}
  \left\{
  \xi\left[ - \genfrac{}{}{}1{1}{2} \sigma
    \left( 1 + 2\rho\varpi + \varpi^2\right)
         + \genfrac{}{}{}1{1}{2} \sigma + \kappa\varpi
    \right]
  + r_0 + q_r - \kappa\theta_{\sigma} \varpi
  \right\}
\\
& {}
  = r_0 K_0
  + K_1\cdot \ee^{x + \varpi\xi}
  \left\{
  \xi\varpi
      \left[ - \genfrac{}{}{}1{1}{2} \sigma \varpi
             + (\kappa - \sigma\rho)
      \right]
  + r_0 + q_r - \kappa\theta_{\sigma} \varpi
  \right\} \geq 0 \,.
\end{align*}
The last inequality follows from $r_0, K_0\in \RR_+$
combined with our conditions on $\varpi$
in \eqref{ineq:varpi}.
Finally, we combine this inequality with
\begin{math}
%\label{e:Cauchy}
  \frac{\partial u}{\partial t} + \mathcal{A} u = 0
\end{math}
in $\HH\times (0,T)$ to derive the desired inequality
\eqref{ineq:Cauchy} for $W_{\pm} = {}- U\pm u$ in place of~$u$.

We finish our proof by applying the weak maximum principle
(Theorem~\ref{thm-MaxPrinciple})
to the functions $W_{\pm}\colon \HH\times [0,T)\to \RR$
which guarantees $W_{\pm}\leq 0$ throughout $\HH\times [0,T)$.
%\null\hfill\qed
\qed
%%%%%%%%%%%%%%%%%
\par\vskip 10pt

%%%%%%%%%%%%%%%%%%%%%%%%%%%%%%%%%%%%%%%%%%%%%%%%%%%%%%%%%%%%%%%%%%%%%%%
%%%%%    Discussion of the boundary conditions    %%%%%%%%%%%%%%%%%%%%%
%%%%%%%%%%%%%%%%%%%%%%%%%%%%%%%%%%%%%%%%%%%%%%%%%%%%%%%%%%%%%%%%%%%%%%%

\section{Discussion of the boundary conditions}
\label{s:BoundaryCond}

It is not difficult to see, as we will show below,
that at any time $t\in (0,T)$
the Cauchy problem for the Heston model
({\S}\ref{ss:Heston} and Appendix~\ref{s:econom-Heston})
imposes on the solution $u(\,\cdot\,,t)\colon \RR^N\to \RR$
the ``boundary'' behavior at infinity (as $|x|\to \infty$)
exhibited precisely by the initial value
$u(\,\cdot\,,0) = u_0\colon \RR^N\to \RR$.
More specifically, this is the case for the European call and put options,
$u_0(x,\xi) = K\, (\mathrm{e}^x - 1)^{+}$ and
$u_0(x,\xi) = K\, (1 - \mathrm{e}^x)^{+}$, respectively,
for $(x,\xi)\in \HH$; cf.\
Eq.~\eqref{log-bc:Ito-oper} (for the European call option) and
{\sc J.-P.\ Fouque}, {\sc G.\ Papanicolaou}, and {\sc K.~R.\ Sircar}
\cite[Fig.\ 1.2 (p.~17) and Fig.\ 1.3 (p.~18)]{FouqPapaSir}
(for both, European call and put options, respectively).
This means that, at least in the case of
European call and put options, the boundary conditions for
$u(x,t)$ as $|x|\to \infty$ are determined by the asymptotic behavior of
the initial data $u_0(x)$ as $|x|\to \infty$.
Hence, if any boundary conditions at infinity
(independent from time $t\in (0,T)$)
are to be imposed on the strong solution to the Cauchy problem,
they should be obeyed also by the initial data (at time $t=0$).
An apparent {\em open question\/} is
\begingroup\it
if those boundary conditions
\endgroup
(i.e., {\it boundary behavior\/}) {\it at infinity\/}
obeyed by the initial data $u_0(x)$ as $|x|\to \infty$
are {\it\bfseries inherited\/}
by the (unique) solution for all times $t\in (0,T)$ and
{\it in what sense\/}.

To illustrate this question, one may consider
the well\--known Black\--Scholes model as treated in
\cite[{\S}1.3, pp.\ 12--18]{FouqPapaSir}
with the closed\--form solution provided in
\cite[Eq.\ (1.37), p.~16]{FouqPapaSir}.
Of economic importance is the {\em Delta hedging\/} ratio,
$\ee^{-x}\cdot \frac{\partial u}{\partial x}$,
defined in
\cite[Eq.\ (1.32), p.~14]{FouqPapaSir}
and calculated in
\cite[{\S}1.3.3, p.~15]{FouqPapaSir}.
The limit of this ratio as $|x|\to \infty$
is obeyed by the closed\--form solution to the Black\--Scholes model
for the European call and put options.
In the analogous form it is imposed also in the Heston model
for the European call option; cf.\
Eq.~\eqref{log-bc:Ito-oper} (equivalent to Eq.~\eqref{bc:Ito-oper})
in the next section
(Appendix~\ref{s:econom-Heston}).
It is well\--known
(see, e.g., \cite[{\S}1.5.3, p.~26]{FouqPapaSir})
that the Black\--Scholes partial differential equation
\cite[Eq.\ (1.35), p.~14]{FouqPapaSir}
can be easily transformed
(by a few elementary substitutions of variables)
into the standard diffusion (i.e., heat) equation over
the space\--time domain $\RR\times (0,\infty)$.
The solution of this standard evolutionary equation is given by
the classical formula
\begin{align}
\label{e:heat}
  u(x,t) = \int_{-\infty}^{+\infty} G(|x-y|;t)\, u_0(y) \,\mathrm{d}y
    \quad\mbox{ for }\, (x,t)\in \RR\times (0,\infty) \,,
\\
\nonumber
\textstyle
    \quad\mbox{ where }\quad
  G(|x-y|;t)\eqdef \frac{1}{ \sqrt{4\pi t} }\cdot
  \exp\left( {}- \frac{|x-y|^2}{4t}\right) \,.
\end{align}
In order to obtain a classical solution
$u\colon \RR\times (0,\tau)\to \RR$ on a sufficiently short
time interval $(0,\tau)\subset (0,\infty)$,
any Lebesgue\--measurable initial data
$u_0\colon \RR\to \RR$ satisfying the growth restriction
$|u_0(x)|\leq M\, \ee^{c x^2}$ for a.e.\ $x\in \RR$ will do.
Here, $M,c\in (0,\infty)$ are some positive constants.
Applying this procedure in any time interval
$(t_0, t_0 + \tau)\subset (0,\infty)$
of length $\tau > 0$, one obtains a classical solution
$u\colon \RR\times (0,\infty)\to \RR$, global in time.
This solution is unique among all classical solutions satisfying
the growth restriction
$|u(x,t)|\leq M\, \ee^{c x^2}$ for all
$(x,t)\in \RR\times (0,\infty)$.
Let us rewrite Eq.~\eqref{e:heat} as
\begin{equation}
\label{e:heat-u_0}
\begin{aligned}
  u(x,t) - u_0(x) = \int_{-\infty}^{+\infty}
    G(|x-y|;t)\, [ u_0(y) - u_0(x) ] \,\mathrm{d}y
\\
    \quad\mbox{ for }\, (x,t)\in \RR\times (0,\infty) \,.
\end{aligned}
\end{equation}

Next, given any $\delta\in (0,1)$, we fix a number
$A_{\delta}\in (0,\infty)$
large enough, such that
\begin{equation}
\label{e:heat_kernel}
\textstyle
    \int_{ A_{\delta}\sqrt{t} }^{+\infty} G(s;t) \,\mathrm{d}s
  = \int_{ A_{\delta} }^{+\infty} G(s';1) \,\mathrm{d}s' < \delta / 2
    \quad\mbox{ for }\, t\in (0,\infty) \,.
\end{equation}
If we wish to impose (possibly inhomogenouos) Dirichlet boundary conditions
on the initial data $u_0(x)$ as $|x|\to \infty$,
for the sake of simplicity, let us assume that
$u_0\colon \RR\to \RR$ is a bounded continuous function,
$|u_0(x)|\leq M\equiv \mathrm{const} < \infty$
for all $x\in \RR$, with the limits
$\lim_{x\to -\infty} u_0(x) = u_0(-\infty)$ and
$\lim_{x\to +\infty} u_0(x) = u_0(+\infty)$.
Making use of Ineq.~\eqref{e:heat_kernel},
we now estimate the difference in
Eq.~\eqref{e:heat-u_0} as $|x|\to \infty$:
\begin{align*}
%\label{eq:heat-u_0}
  |u(x,t) - u_0(x)| \leq
    \int_{ - A_{\delta}\sqrt{t} }^{ + A_{\delta}\sqrt{t} }
    G(|z|;t)\, |u_0(x+z) - u_0(x)| \,\mathrm{d}z
\\
  + \int_{-\infty}^{ - A_{\delta}\sqrt{t} }
    G(|z|;t)\, |u_0(x+z) - u_0(x)| \,\mathrm{d}z
  + \int_{ + A_{\delta}\sqrt{t} }^{+\infty}
    G(|z|;t)\, |u_0(x+z) - u_0(x)| \,\mathrm{d}z
\\
  \leq \left( \int_{-\infty}^{+\infty} G(|z|;t) \,\mathrm{d}z \right)
  \cdot \sup_{x\in \RR}
        \sup_{ |z|\leq A_{\delta}\sqrt{t} } |u_0(x+z) - u_0(x)|
  + 2M\int_{ |z'|\geq A_{\delta} } G(|z'|;1) \,\mathrm{d}z'
\\
  \leq \sup_{x\in \RR}
       \sup_{ |z|\leq A_{\delta}\sqrt{t} } |u_0(x+z) - u_0(x)|
  + 2M\delta
    \quad\mbox{ for }\, (x,t)\in \RR\times (0,\infty) \,.
\end{align*}
Letting $x\to \pm\infty$ we arrive at
\begin{equation*}
%\label{lim:heat-u_0}
  \limsup_{ x\to \pm\infty } |u(x,t) - u_0(x)| \leq 2M\delta
    \quad\mbox{ for every }\, t\in (0,\infty) \,.
\end{equation*}
The number $\delta\in (0,1)$ being arbitrary, we conclude that
$\lim_{x\to -\infty} u(x,t) = u_0(-\infty)$ and
$\lim_{x\to +\infty} u(x,t) = u_0(+\infty)$ as desired.

Neumann boundary conditions can be treated in a similar manner
using the following formula derived from Eq.~\eqref{e:heat-u_0}
by simple differentiation:
\begin{equation*}
%\label{e_x:heat-u_0}
\begin{aligned}
    \frac{\partial u}{\partial x} (x,t)
  - \frac{\partial u_0}{\partial x} (x)
  = \int_{-\infty}^{+\infty} G(|x-y|;t)
  \left[ \frac{\partial u_0}{\partial x} (y)
       - \frac{\partial u_0}{\partial x} (x)
  \right] \,\mathrm{d}y
\\
    \quad\mbox{ for }\, (x,t)\in \RR\times (0,\infty) \,.
\end{aligned}
\end{equation*}
However, caution must be payed to the ``weighted''
Neuman boundary conditions suggested by Eq.~\eqref{log-bc:Ito-oper}
(Appendix~\ref{s:econom-Heston}):
The European call option prescribes the limits
\begin{equation*}
  L_{-\infty} =
  \lim_{x\to -\infty}
  \left( \ee^{-x}\cdot \frac{\partial u}{\partial x} (x,t) \right)
  = 0
  \quad\mbox{ and }\quad
  L_{+\infty} =
  \lim_{x\to +\infty}
  \left( \ee^{-x}\cdot \frac{\partial u}{\partial x} (x,t) \right)
  = 1
\end{equation*}
for the {\em Delta hedging\/} ratio
(cf.\ \cite[Eq.\ (1.37), p.~16]{FouqPapaSir}).
In the case of the Black\--Scholes model,
these limits can be verified in a manner similar to
the Dirichlet boundary conditions above.
We leave the details to an interested reader.

The boundary condition as $|x|\to \infty$, given by
Eq.~\eqref{bc_x:bound_u} with $\gamma > 2$,
that we have used in the definition of the Heston operator by
Eq.~\eqref{e:Heston-oper}
(cf.\ \cite[Eq.\ (2.24), p.~12]{AlziaryTak}),
is in fact weaker than the corresponding boundary condition
in Eq.~\eqref{log-bc:Ito-oper}.
Nevertheless, our condition \eqref{bc_x:bound_u}
is still sufficient for obtaining a unique solution to
the Heston model.
We recall that the choice of $\gamma > 2$ is necessary to ensure
$u_0\in H$ for the European call option with the initial condition
$u_0(x,\xi) = K\, (\mathrm{e}^x - 1)^{+}$ for $(x,\xi)\in \HH$.

\appendix
%%%%%%%%%%%%%%%%%%%%%%%%%%%%%%%%%%%%%%%%%%%%%%%%%%%%%%%%%%%%%%%%%%%%%%%
%%%%%    The Heston model in finance    %%%%%%%%%%%%%%%%%%%%%%%%%%%%%%%
%%%%%%%%%%%%%%%%%%%%%%%%%%%%%%%%%%%%%%%%%%%%%%%%%%%%%%%%%%%%%%%%%%%%%%%

\section{Appendix: The Heston model in finance}
\label{s:econom-Heston}

A number of {\it\bfseries stochastic volatility models\/} for
derivative pricing
(e.g., of {\it call\/} or {\it put options\/} on stocks)
are known in the literature; see
{\sc J.-P.\ Fouque}, {\sc G.\ Papanicolaou}, and {\sc K.~R.\ Sircar}
\cite[Table 2.1, p.~42]{FouqPapaSir}.
We focus our attention on
{\sc S.~L.\ Heston}'s model \cite{Heston}
which has attracted significant attention of
a broad community of researchers from Finance and Mathematics.
We consider this model under a \emph{risk neutral measure\/}
via equations $(1)-(4)$ in \cite[pp.\ 328--329]{Heston}.
The model is defined on a filtered probability space
$(\Omega, \mathcal{F}, (\mathcal{F}_t)_{t\geqslant 0}, \mathbb{P})$,
where $\mathbb{P}$ is a risk neutral probability measure, and
the filtration $(\mathcal{F}_t)_{t\geqslant 0}$
satisfies the usual conditions.
Denoting by $S_t$ the {\em stock price\/} and by
$V_t$ the (stochastic) {\em variance\/} of the stock market
at (the real) time $t\geq 0$,
the {\sc Heston\/} model requires that
the unknown pair $(S_t,V_t)_{t\geqslant 0}$ satisfies
the following system of {\em stochastic\/} differential equations,
\begin{equation}
\label{e:SVmodel}
\begin{aligned}
  \frac{\,\mathrm{d} S_t}{S_t} =
  {}- q_r \,\mathrm{d}t + \sqrt{V_t} \,\mathrm{d} W_t \,,\quad
  \,\mathrm{d} V_t =
    \kappa\, (\theta - V_t) \,\mathrm{d}t
  + \sigma\, \sqrt{V_t}\, \,\mathrm{d} Z_t \,.
\end{aligned}
\end{equation}
Here, $(W_t)_{t\geqslant 0}$ and $(Z_t)_{t\geqslant 0}$
are two Brownian motions with the {\it correlation coefficient\/}
$\rho\in (-1,1)$, a constant given by
$\,\mathrm{d}\langle W,Z\rangle_t = \rho\,\mathrm{d}t$.
Furthermore,
$q_r = q-r\in \RR$ and $\sigma, \kappa, \theta\in (0,\infty)$
are some given constants whose economic meaning is explained, e.g., in
{\sc B.\ Alziary} and {\sc P.\ Tak\'a\v{c}}
\cite[Sect.~1, pp.\ 3--4]{AlziaryTak}
or
{\sc C.\ Chiarella}, {\sc B.\ Kang}, and {\sc G.~H. Meyer}
\cite[Chapt.~2, pp.\ 3--5]{Chiarella-K-M-2015}.

If $X_t = \ln (S_t/K)$ denotes the (natural) logarithm of
the {\em scaled stock price\/} $S_t / K$ at time $t\geq 0$
relative to the strike price $K>0$ at maturity $T>0$,
then the pair $(X_t,V_t)_{t\geqslant 0}$ satisfies
the following system of stochastic differential equations,
\begin{equation}
\label{e:XVmodel}
\begin{aligned}
  \,\mathrm{d} X_t =
  {}- \left( q_r + \genfrac{}{}{}1{1}{2} V_t \right) \,\mathrm{d}t
                 + \sqrt{V_t} \,\mathrm{d} W_t \,,\quad
  \,\mathrm{d} V_t = \kappa\, (\theta - V_t) \,\mathrm{d}t
                 + \sigma\, \sqrt{V_t}\, \,\mathrm{d} Z_t \,.
\end{aligned}
\end{equation}

Following \cite[Sect.~4]{Davis-Obloj}, let us consider
a European call option written in this market with {\em payoff\/}
$\hat{h}(S_T,V_T)\equiv \hat{h}(S_T)\geq 0$ at maturity $T$, where
$\hat{h}(s) = (s-K)^{+}$ for all $s>0$.
Recalling {\sc Heston}'s notation in
\cite[Eq.~(11), p.~330]{Heston}, we denote by
$x = X_t(\omega)\in \RR$ the logarithm of
the {\it\bfseries spot price\/} of stock and by
$v = V_t(\omega)\in \RR$ the {\it\bfseries variance\/}
of stock market at time~$t$.
We set
$h(x,v)\equiv h(x) = K\, (\ee^x - 1)^{+}$
for all $x = \ln (s/K)\in \mathbb{R}$, so that
$h(x) = \hat{h}(s) = \hat{h}(K\ee^x)$ for $x\in \RR$.
Hence, if the instant values
$(X_t(\omega) ,\, V_t(\omega)) = (x,v)\in \HH$
are known at time $t\in (0,T)$, where
$\mathbb{H} = \mathbb{R}\times (0,\infty) \subset \RR^2$,
the \emph{arbitrage\--free price\/}
$P^h_t$ of the European call option at this time is given by
the following expectation formula
(with respect to the risk neutral probability measure $\mathbb{P}$)
which is justified in \cite{Davis-Obloj} and \cite{Takac-12}:
$P^h_t = p(X_t,V_t,t)$ where
\begin{equation}
\label{e:optionprice}
\begin{aligned}
  p(x,v,t)
& = \mathrm{e}^{-r(T-t)}\, \mathbb{E}_{\mathbb{P}}
    \left[ \hat{h}(S_T)\mid \mathcal{F}_t \right]
  = \mathrm{e}^{-r(T-t)}\, \mathbb{E}_{\mathbb{P}}
    \left[ h(X_T)\mid \mathcal{F}_t \right]
\\
& = \mathrm{e}^{-r(T-t)}\, \mathbb{E}_{\mathbb{P}}
    \left[ h(X_T)\mid X_t = x ,\ V_t = v \right] \,.
\end{aligned}
\end{equation}
Furthermore, $p$ solves the ({\it terminal\/} value) Cauchy problem
\begin{equation}
\label{e:IVP.p}
\left\{
\begin{alignedat}{2}
  \frac{\partial p}{\partial t} + \mathcal{G}_t\, p - rp &= 0 \,,\quad
&&(x,v,t)\in \HH\times (0,T) \,;
\\
  p(x,v,T) &= h(x) \,,\quad &&(x,v)\in \HH \,,
\end{alignedat}
\right.
\end{equation}
with $\mathcal{G}_t$ being
the (time\--independent) infinitesimal generator of
the time\--homogeneous Markov process $(X_t,V_t)$;
cf.\ {\sc A.\ Friedman} \cite[Chapt.~6]{Friedman-75}
or {\sc B.\ {\O}ksendal} \cite[Chapt.~8]{Oksendal}.
Indeed, to justify Eq.~\eqref{e:IVP.p},
we take advantage of It\^{o}'s formula to derive the following equation
from eqs.\ \eqref{e:XVmodel} and \eqref{e:optionprice}:
\begin{equation}
\label{eq:Heston}
  \left( \frac{\partial}{\partial t} + \mathbf{A} \right)
  U(s,v,t) = 0 \,,
\end{equation}
where we use the instant values
$(s,v) = (S_t(\omega) ,\, V_t(\omega))\in (0,\infty)^2$
and substitute
$U(s,v,t) = p(x,v,t)$ with
$s = K\mathrm{e}^x$, ${\mathrm{d}s} / {\mathrm{d}x} = s$, and
\begin{equation}
\label{eq:Ito-oper}
\begin{aligned}
  (\mathbf{A}U)(s,v,t)\eqdef
{}&
    \frac{1}{2}\, v\cdot
  \left(
    s^2\, \frac{\partial^2 U}{\partial s^2}(s,v,t)
  + 2\rho\sigma\, s\, \frac{\partial^2 U}{\partial s\;\partial v}(s,v,t)
  + \sigma^2\, \frac{\partial^2 U}{\partial v^2}(s,v,t)
  \right)
\\
& {}
  + (r-q)\, s\, \frac{\partial U}{\partial s}(s,v,t)
  + \left[ \kappa (\theta - v) - \lambda\cdot v\right]
    \frac{\partial U}{\partial v}(s,v,t)
\\
& {}
  - r\, U(s,v,t)
  \quad\mbox{ for $s>0$, $v>0$, and $0\leq t\leq T$, }
\end{aligned}
\end{equation}
denotes the (usual)
{\it\bfseries Black\--Scholes\-(-\^{I}to) operator\/}.
Eq.~\eqref{eq:Heston} entails
the desired diffusion equation \eqref{e:IVP.p} using
\begin{align*}
    \frac{\partial p}{\partial x}(x,v,t)
& {}
  = s\, \frac{\partial U}{\partial s}(s,v,t) \,,
\\
    \frac{\partial^2 p}{\partial x^2}(x,v,t)
& {}
  = s\, \frac{\partial U}{\partial s}(s,v,t)
  + s^2\, \frac{\partial^2 U}{\partial s^2}(s,v,t)
  = \frac{\partial p}{\partial x}(x,v,t)
  + s^2\, \frac{\partial^2 U}{\partial s^2}(s,v,t) \,.
\end{align*}
Hence, the function
$\overline{p}\colon (x,v,t)\mapsto p(x,v,T-t)$
verifies a linear ({\it initial\/} value) Cauchy problem derived from
\eqref{e:IVP.p}.
The functional\--analytic formulation of this problem is given in
Eq.~\eqref{e:Cauchy}
with the initial data corresponding to
$\overline{p}(x,v,0) = p(x,v,T) = h(x)$ at $t=0$.

We have replaced the meaning of the temporal variable $t$
as real time ($t\leq T$)
by the {\em time to maturity\/} $t$ ($t\geq 0$), so that
the real time (time to maturity) has become $\tau = T-t$.
According to
{\sc S.~L.\ Heston} \cite[Eq.~(6), p.~329]{Heston},
the unspecified term $\lambda(x,v,T-t)$ in the second drift term
on the right\--hand side of Eq.~\eqref{eq:Ito-oper}
(with ${\partial U} / {\partial v}$)
represents the {\it\bfseries price of volatility risk\/}
and is specifically chosen to be
$\lambda(x,v,T-t)\equiv \lambda v$ with a constant $\lambda\geq 0$.
As we have already pointed out in the Introduction
(Section~\ref{s:Intro}),
we can treat much more general initial conditions
$\overline{p}(x,v,0) = h(x,v)$
than just those given by
the terminal condition for the European call option,
$p(x,v,T) = h(x) =  K\, (\ee^x - 1)^{+}$ for $(x,v)\in \HH$,
which does not depend on the instant value
$v = V_T(\omega)$ of the variance at maturity~$T$;
see Section~\ref{s:Main}.

Next, we eliminate the constants $r\in \RR$ and $\lambda\geq 0$,
respectively, from Eq.~\eqref{e:IVP.p} by substituting
\begin{equation*}
%\label{e:U-->u}
  p^{*}(x,v,t)\eqdef \ee^{rt}\, \overline{p}(x,v,t)
                   = \ee^{rt}\, p(x,v,T-t) = \ee^{rt}\, U(x,v,T-t)
\end{equation*}
for $\overline{p} (x,v,t)$ and replacing
$\kappa$ by $\kappa^{*} = \kappa + \lambda > 0$ and
$\theta$ by
\begin{math}
  \theta^{*} = \frac{\kappa\theta}{\kappa + \lambda} > 0 .
\end{math}
Hence, we may set $r = \lambda$ $= 0$.
Finally, we introduce also the re\--scaled variance
$\xi = v / \sigma > 0$ for $v\in (0,\infty)$
and abbreviate
$\theta_{\sigma}\eqdef \theta / \sigma\in \RR$.
These substitutions have a simplifying effect
on our calculations.
Eq.~\eqref{e:IVP.p} then yields the initial value problem
\eqref{e:Cauchy} for the unknown function
$u(x,\xi,t) = p^{*}(x, \sigma\xi ,t)$,
with the initial data
$u_0(x,\xi) = p^{*}(x, \sigma\xi ,0)\equiv h(x)$ at $t=0$,
where the (autonomous linear) {\it\bfseries Heston operator\/}
$\mathcal{A}$, derived from Eq.~\eqref{e:IVP.p},
takes the standard elliptic form \cite[Eq.\ (2.8), p.~8]{AlziaryTak};
we prefer to use the asymmetric ``divergence'' form of $\mathcal{A}$
given by Eq.~\eqref{e:Heston-oper}
(cf.\ \cite[Eq.\ (2.9), p.~8]{AlziaryTak}).

The original work by
{\sc S.~L.\ Heston} \cite[Eq.~(9), p.~330]{Heston}
imposes the following {\it\bfseries boundary conditions\/}:
The {\it\bfseries boundary operator\/} as $v\to 0+$,
defined by
\begin{equation*}
%\label{eq:Ito-oper_BC}
\begin{aligned}
  (\mathbf{B}U)(s,0,t)\eqdef
    (r-q)\, s\, \frac{\partial U}{\partial s}(s,0,t)
  + \kappa\theta\, \frac{\partial U}{\partial v}(s,0,t)
  - r\, U(s,0,t)
\\
  \quad\mbox{ for $s>0$, $v=0$, and $0\leq t\leq T$, }
\end{aligned}
\end{equation*}
transforms the left\--hand side of
the boundary condition as $v\to 0+$,
\begin{equation}
\label{bc:Heston}
  \left( \frac{\partial}{\partial t} + \mathbf{B} \right)
  U(s,0,t) = 0 \,,
\end{equation}
into the following (logarithmic) form
on the boundary $\partial\HH = \RR\times \{ 0\}$ of $\HH$:
\begin{equation}
\label{log:Ito-oper_BC}
\begin{aligned}
& \ee^{-r\tau}\,
  \left( \frac{\partial}{\partial\tau} + \mathbf{B} \right)
  U(s,0,\tau) \Big\vert_{\tau = T-t}
  = {}- \left( \frac{\partial}{\partial t} + \mathcal{B} \right)
  u(x,0,t)
\\
& = {}
  - \frac{\partial u}{\partial t}(x,0,t)
  - q_r\, \frac{\partial u}{\partial x}(x,0,t)
  + \kappa\theta_{\sigma}\,
    \frac{\partial u}{\partial\xi}(x,0,t)
\\
& \quad\mbox{ for $x\in \RR$ and $0 < t < \infty$. }
\end{aligned}
\end{equation}
The remaining boundary conditions
(in addition to \eqref{bc:Heston}),
\begin{equation}
\label{bc:Ito-oper}
\left\{
\begin{aligned}
& U(0,v,t) = 0 \,;
\\
& \lim_{s\to \infty}\,
  \frac{\partial}{\partial s}\, ( U(s,v,t) - s ) = 0 \,;
\\
& \lim_{v\to \infty}\, ( U(s,v,t) - s ) = 0 \,,
\end{aligned}
\right.
\end{equation}
for $s>0$, $v>0$, and $0\leq t\leq T$, become
(in addition to \eqref{log:Ito-oper_BC}),
\begin{equation}
\label{log-bc:Ito-oper}
\left\{
\begin{alignedat}{2}
  u(-\infty,\xi,t)\eqdef
  \lim_{x\to -\infty}\,
  \left( u(x,\xi,t) - K\mathrm{e}^{x+rt} \right)
& = 0 \quad\mbox{ for $\xi > 0$; }
\\
  \lim_{x\to +\infty}\,
  \left[
  \ee^{-x}\cdot
  \frac{\partial}{\partial x}\,
  \left( u(x,\xi,t) - K\mathrm{e}^{x+rt} \right)
  \right]
& = 0 \quad\mbox{ for $\xi > 0$; }
\\
  \lim_{\xi\to \infty}\,
  \left( u(x,\xi,t) - K\mathrm{e}^{x+rt} \right)
& = 0 \quad\mbox{ for $x\in \RR$, }
\end{alignedat}
\right.
\end{equation}
at all times $t\in (0,\infty)$.
We remark that the first equation in \eqref{log-bc:Ito-oper}
is a consequence of the initial conditions
(for a European call option)
\begin{math}
  0\leq u_0(x,\xi) = K\, (\mathrm{e}^x - 1)^{+}
\end{math}
for $(x,\xi)\in \HH$,
the lower bound $u(x,\xi,t)\geq 0$ for all
$(x,\xi,t)\in \HH\times (0,\infty)$
obtained from the weak maximum principle in
Theorem~\ref{thm-MaxPrinciple}, and the upper bound
$u(x,\xi,t)\leq K\mathrm{e}^{x+rt}$ for all
$(x,\xi,t)\in \HH\times (0,\infty)$
obtained from the weak maximum principle in
Corollary~\ref{cor-MaxPrinciple} with
$K_1 = K$, $K_0 = 0$, and $\varpi = 0$.

%%%%%%%%%%%%%%%%%%%%%%%%%%%%%%%%%%%%%%%%%%%%%%%%%%%%%%%%%%%%%%%%%%%%%%%
%%%%%    Regularity for the elliptic Heston operator in $L^2$    %%%%%%
%%%%%%%%%%%%%%%%%%%%%%%%%%%%%%%%%%%%%%%%%%%%%%%%%%%%%%%%%%%%%%%%%%%%%%%

\section{Appendix: Weighted Sobolev spaces and boundary traces}
\label{s:Sobolev_trace}

We denote by
$H^1(B^{+}_R(x_0,0); \mathfrak{w})$
the weighted Sobolev space of all functions
$f\in W^{1,2}_{\mathrm{loc}}(B^{+}_R(x_0,0))$
whose norm defined below is finite,
\begin{equation}
\label{loc:norm-H^1}
\begin{aligned}
  \| f\|_{ H^1(B^{+}_R(x_0,0); \mathfrak{w}) }^2
& \eqdef \int_{ B^{+}_R(x_0,0) }
    \left( |f_{x}|^2 + |f_{\xi}|^2\right)
    \cdot \xi^{\beta}\cdot \,\mathrm{d}x \,\mathrm{d}\xi
\\
& {}
  + \int_{ B^{+}_R(x_0,0) } |f(x,\xi)|^2
    \cdot \xi^{\beta - 1}\cdot \,\mathrm{d}x \,\mathrm{d}\xi
    < \infty \,.
\end{aligned}
\end{equation}
Let us recall that the weighted Sobolev space
$H^2(B^{+}_R(x_0,0); \mathfrak{w})$
has been defined by its norm in Eq.~\eqref{e:norm-H^2}.
As we will see later, in Lemma~\ref{lem-Sobol-Morrey},
functions from the weighted Sobolev spaces
$H^j(B^{+}_R(x_0,0); \mathfrak{w})$; $j=1,2$,
must satisfy certain homogeneous boundary conditions as
$\xi\to 0+$, i.e., near the boundary
\begin{math}
  \partial\HH \cap B^{+}_R(x_0,0) =
  (x_0 - R,\, x_0 + R) \times \{ 0\} .
\end{math}
We will see in the proof of this result
(Lemma~\ref{lem-Sobol-Morrey}), as well,
that these boundary conditions, if satisfied by a function,
imply that this function belongs to a particular weighted Sobolev space.

A simple motivation for such a result is the classical Sobolev space
$W^{1,2}(0,1)$:
This space is continuously imbedded into the H\"older space
$C^{1/2}[0,1]$; hence, the limit
$\lim_{\xi\to 0+} f(\xi) = f(0)\in \RR$
is valid for every function $f\in W^{1,2}(0,1)$.
By Hardy's inequality
(proved in
 {\sc G.~H.\ Hardy}, {\sc J.~E.\ Littlewood}, and {\sc G.\ P\'olya}
 \cite[Theorem 330, pp.\ 245--246]{Hardy-Polya}; see also
 {\sc A.\ Kufner} \cite[Section~5]{Kufner}),
we have
\begin{equation*}
  \int_0^1 |f(\xi) - f(0)|^2\, \xi^{-2} \,\mathrm{d}\xi
  \leq \mathrm{const}\cdot
  \int_0^1 |f'(\xi)|^2\, \xi^{-2} \,\mathrm{d}\xi
\end{equation*}
with a positive constant independent from $f\in W^{1,2}(0,1)$,
whenever $f$ satisfies $f(1) = f(0)$.
Clearly, the homogeneous boundary condition $f(0) = 0$ is valid
if and only if
\begin{math}
  \int_0^1 |f(\xi)|^2\, \xi^{-2} \,\mathrm{d}\xi < \infty .
\end{math}
If this is the case, then even
\begin{math}
  \lim_{\xi\to 0+} ( f(\xi) / \xi^{1/2} ) = 0
\end{math}
holds.

A much less trivial example appears in our earlier work
\cite[Sect.~10 (Appendix)]{AlziaryTak}.
We now show this example only for the bounded half\--disc
$B^{+}_R(x_0,0)\subset \HH$ near the boundary
$\partial\HH \cap B^{+}_R(x_0,0)$.
Let us fix any
$r\in (0,R)$ and set
$\varrho\equiv \varrho(r) = \sqrt{R^2 - r^2}$.
Given any function
$f\in W^{1,2}_{\mathrm{loc}}( B^{+}_R(x_0,0) )$,
we begin with the identity
\begin{equation*}
%\label{ineq:trace:v}
  \frac{\partial}{\partial\xi}
  \left( \xi^{\beta - 1}\, f(x,\xi)^2\right)
  = 2\, f f_{\xi} \cdot \xi^{\beta - 1}
  + (\beta - 1)\, f(x,\xi)^2\cdot \xi^{\beta - 2} \,,
\end{equation*}
for $(x,\xi)\in B^{+}_R(x_0,0)$ satisfying
$x\in (x_0 - r,\, x_0 + r)$ and $0 < \xi < \varrho$.
We apply Cauchy's inequality
\begin{equation*}
    2\, | f f_{\xi} |
  \leq \frac{\beta - 1}{2}\, \xi^{-1}\, f^2
     + \frac{2}{\beta - 1}\, \xi\, f_{\xi}^2
\end{equation*}
to the equation above to estimate the partial derivative
\begin{equation}
\label{e:d(xi^(beta-1)f)/d_xi}
\begin{aligned}
    \frac{\beta - 1}{2}\, \xi^{\beta - 2}\, f^2
  - \frac{2}{\beta - 1}\, \xi^{\beta}\, f_{\xi}^2
& \leq
  \frac{\partial}{\partial\xi}
  \left( \xi^{\beta - 1}\, f(x,\xi)^2\right)
\\
& \leq
    \frac{3}{2}\, (\beta - 1)\, \xi^{\beta - 2}\, f^2
  + \frac{2}{\beta - 1}\, \xi^{\beta}\, f_{\xi}^2 \,.
\end{aligned}
\end{equation}
Assuming the integrability
\begin{equation*}
\textstyle
  \iint_{ B^{+}_R(x_0,0) }
    \left( f_{\xi}^2\cdot \xi^{\beta} + f^2\cdot \xi^{\beta - 1}
    \right) \,\mathrm{d}x \,\mathrm{d}\xi < \infty ,
\end{equation*}
we deduce from the inequalities in \eqref{e:d(xi^(beta-1)f)/d_xi}
that the function
\begin{equation*}
  \xi \,\longmapsto\, F_{\beta,r}(\xi)\eqdef \xi^{\beta - 1}
\textstyle
    \int_{x_0 - r}^{x_0 + r} f(x,\xi)^2 \,\mathrm{d}x
  \colon (0,\varrho) \,\longrightarrow\, \RR
\end{equation*}
is absolutely continuous over the compact interval
$[0,\varrho]$ with finite boundary limits
\begin{equation*}
%\label{lim:d(xi^(beta-1)f)/d_xi}.
  F_{\beta,r}(0)\eqdef \liminf_{\xi\to 0+} F_{\beta,r}(\xi)
    \quad\mbox{ and }\quad
  F_{\beta,r}(\varrho) =
  \varrho^{\beta - 1}
    \int_{x_0 - r}^{x_0 + r} f(x,\varrho)^2 \,\mathrm{d}x
  \quad\mbox{ (as $\xi\to {\varrho}-$) }
\end{equation*}
if and only if
\begin{equation*}
  \int_0^{\varrho} F_{\beta,r}(\xi)\cdot \xi^{-1} \,\mathrm{d}\xi
  = \int_0^{\varrho} \xi^{\beta - 2} \int_{x_0 - r}^{x_0 + r}
    f^2 \,\mathrm{d}x \,\mathrm{d}\xi
  < \infty \,.
\end{equation*}
However, the last integral is finite if and only if
$F_{\beta,r}(0) = 0$.
If this is the case, then also the limit
$\lim_{\xi\to 0+} F_{\beta,r}(\xi) = 0$.
We conclude that the homogeneous boundary condition given by
$F_{\beta,r}(0) = 0$
is equivalent with the convergence of the last integral.
Greater details can be found in
\cite[Sect.~10, pp.\ 43--48]{AlziaryTak}, Lemmas 10.1 through 10.5.

We will follow a similar procedure in treating the case
$f\in H^2(B^{+}_R(x_0,0); \mathfrak{w})$.
More precisely, we wish to show that if we take a weaker norm,
\begin{math}
  \|\cdot\|_{ H^2(B^{+}_R(x_0,0); \mathfrak{w}) }^{\flat}
\end{math}
on $H^2(B^{+}_R(x_0,0); \mathfrak{w})$
defined in Eq.~\eqref{eq:norm-H^2} below,
the restriction mapping
\begin{equation}
\label{eq:restr-H^2}
  f \,\longmapsto\, f\vert_{ \overline{B}^{+}_{ R/\sqrt{2} }(x_0,0) }
  \colon \tilde{H}^2(B^{+}_R(x_0,0); \mathfrak{w}) \,\longrightarrow\,
         H^2( B^{+}_{ R/\sqrt{2} }(x_0,0); \mathfrak{w} )
\end{equation}
is still continuous from
$\tilde{H}^2(B^{+}_R(x_0,0); \mathfrak{w})$ to
$H^2( B^{+}_{ R/\sqrt{2} }(x_0,0); \mathfrak{w} )$,
where
$\tilde{H}^2(B^{+}_R(x_0,0); \mathfrak{w})$
stands for the completion of the Sobolev space
$H^2(B^{+}_R(x_0,0); \mathfrak{w})$ under the new norm
\begin{math}
  \|\cdot\|_{ H^2(B^{+}_R(x_0,0); \mathfrak{w}) }^{\flat}
\end{math}
defined as follows:
\begin{align} 
\label{eq:norm-H^2}
&
\begin{aligned}
& \left(
  \| f\|_{ H^2(B^{+}_R(x_0,0); \mathfrak{w}) }^{\flat}
  \right)^2
  \eqdef \int_{ B^{+}_R(x_0,0) } 
    \left( |f_{xx}|^2 + |f_{x\xi}|^2 + |f_{\xi\xi}|^2\right)
    \cdot \xi^{\beta + 1}\cdot \,\mathrm{d}x \,\mathrm{d}\xi
\\
& {}
  + \int_{ B^{+}_R(x_0,0) }
    \left( |f_{x}|^2 + |f_{\xi}|^2\right)
    \cdot \xi^{\beta}\cdot \,\mathrm{d}x \,\mathrm{d}\xi
  + \int_{ B^{+}_R(x_0,0) } |f(x,\xi)|^2
    \cdot \xi^{\beta - 1}\cdot \,\mathrm{d}x \,\mathrm{d}\xi
\end{aligned}
\\
\nonumber
&
\begin{aligned}
  \equiv
    [f]_{ H^2(B^{+}_R(x_0,0); \mathfrak{w}) }^2
  + \| f\|_{ H^1(B^{+}_R(x_0,0); \mathfrak{w}) }^2
    < \infty \,,
\end{aligned}
\end{align}
where
\begin{math}
  [\,\cdot\,]_{ H^2(B^{+}_R(x_0,0); \mathfrak{w}) }
\end{math}
is a seminorm on
$H^2(B^{+}_R(x_0,0); \mathfrak{w})$ defined by
\begin{equation}
\label{loc:seminorm-H^2}
  [f]_{ H^2(B^{+}_R(x_0,0); \mathfrak{w}) }
  \eqdef
  \left( \int_{ B^{+}_R(x_0,0) }
    \left( |f_{xx}|^2 + |f_{x\xi}|^2 + |f_{\xi\xi}|^2\right)
    \cdot \xi^{\beta + 1}\cdot \,\mathrm{d}x \,\mathrm{d}\xi
  \right)^{1/2} \,.
\end{equation}
It is easy to see that
$\tilde{H}^2(B^{+}_R(x_0,0); \mathfrak{w})$
consists of all functions
$f\in W^{2,2}_{\mathrm{loc}}(B^{+}_R(x_0,0))$
that satisfy
\begin{math}
  \| f\|_{ H^2(B^{+}_R(x_0,0); \mathfrak{w}) }^{\flat} < \infty .
\end{math}
The continuous restriction mapping in \eqref{eq:restr-H^2}
is termed a {\em restricted imbedding\/}, by
{\sc R.~A.\ Adams} and {\sc J.~J.~F. Fournier}
\cite[Chapt.~6, {\S}6.1, p.~167]{AdamsFournier}.
In the course of the proof of this restriction imbedding,
we will obtain also certain boundary conditions
(i.e., trace results as $\xi\to 0+$)
on the boundary
\begin{math}
  \partial\HH\cap \overline{B}^{+}_{ R/\sqrt{2} }(x_0,0) .
\end{math}

Keeping in mind that some of the constants in our estimates below
may depend on the choice of $x_0\in \RR$,
we suppress the dependence on $x_0$ in the notation
for the half\--disc $B^{+}_R(x_0,0)\subset \HH$ and, thus, write only
$B^{+}_R\equiv B^{+}_R(x_0,0)$.
We further denote by
$Q^{+}_r\equiv Q^{+}_r(x_0,0)$ the open rectangle
\begin{equation*}
  Q^{+}_r(x_0,0)\eqdef (x_0 - r,\, x_0 + r)\times (0,r)\subset \HH
\end{equation*}
(a ``half\--square'') with side lengths $2r$ and $r\in (0,\infty)$.
Its closure in $\mathbb{R}^2$ is denoted by $\overline{Q}^{+}_r$.

Our first lemma is an essential estimate for obtaining
the boundary trace as $\xi\to 0+$.
Given a function
$u\in W^{1,1}_{\mathrm{loc}}(B^{+})$,
we abbreviate the gradient $\nabla u = (u_x,u_{\xi})$.

%%%%%%%%%%%%%%%%%%%%%%%%%%%%%%%%%%%%%%%%%%%%%%%%%%%%%%%%%%%%%%%%%%%%%%%
%%%%%    A pointwise (x) trace inequality (Lemma)    %%%%%%%%%%%%%%%%%%
%%%%%%%%%%%%%%%%%%%%%%%%%%%%%%%%%%%%%%%%%%%%%%%%%%%%%%%%%%%%%%%%%%%%%%%
\begin{lemma}\label{lem-Trace_ineq}
{\rm ($\xi$-derivative inequalities.)}$\;$
Let\/ $\beta > 0$ and\/ $R > 0$, and set\/ $r = R / \sqrt{2}$.
Assume that\/
$u\in \tilde{H}^2(B^{+}_R; \mathfrak{w})$.
Then $Q^{+}_r\subset B^{+}_R$
and the following inequalities hold at almost every point\/
$(x,\xi)\in Q^{+}_r$,
\begin{equation}
\label{ineq:ano_trace_pt:u}
\begin{aligned}
&   \frac{\beta}{2}\, \xi^{\beta - 1}\cdot |\nabla u(x,\xi)|^2
  - \frac{2}{\beta}\, \xi^{\beta + 1}\cdot
    \left| \partial_{\xi}\nabla u\right|^2
  \leq
  \frac{\partial}{\partial\xi}
  \left( \xi^{\beta}\cdot |\nabla u(x,\xi)|^2 \right)
\\
& {}
  \leq
    \frac{3\beta}{2}\, \xi^{\beta - 1}\cdot |\nabla u(x,\xi)|^2
  + \frac{2}{\beta}\, \xi^{\beta + 1}\cdot
    \left| \partial_{\xi}\nabla u\right|^2 \,.
\end{aligned}
\end{equation}
\end{lemma}
%%%%%%%%%%%%%%%%%%%%%%%%%%%%%%%%%%%%%%%%%%%%%%%%%%%%%%%%%%%%%%%%%%%%%%%
%\par\vskip 10pt

%%%%%%%%%%%%%%%%%
\proof
The following partial derivatives exist almost everywhere in
$B^{+}_R$; we first calculate
\begin{equation}
\label{e:trace:u}
  \frac{\partial}{\partial\xi}
  \left( \xi^{\beta}\cdot |\nabla u(x,\xi)|^2 \right)
  = \beta\, \xi^{\beta - 1}\cdot |\nabla u(x,\xi)|^2
  + 2\xi^{\beta}\cdot ( \nabla u\cdot \partial_{\xi}\nabla u ) \,,
\end{equation}
with the scalar product
\begin{math}
  \nabla u\cdot \partial_{\xi}\nabla u
  = u_x u_{x\xi} + u_{\xi} u_{\xi\xi}
\end{math}
in $\RR^2\subset \CC^2$, then estimate the scalar product
on the right\--hand side by Cauchy's inequality,
%
%\begin{multline}
\begin{equation}
\label{ineq_Cauchy:trace:u}
{}\pm 2( \nabla u\cdot \partial_{\xi}\nabla u )
  \leq 2\, |\nabla u|\cdot \left| \partial_{\xi}\nabla u \right|
  \leq
    \frac{\beta}{2}\, \xi^{-1}\cdot |\nabla u|^2
  + \frac{2}{\beta}\, \xi\cdot \left| \partial_{\xi}\nabla u\right|^2
\end{equation}
%\end{multline}
%
for a.e.\ $(x,\xi)\in B^{+}_R$.
We apply Ineq.~\eqref{ineq_Cauchy:trace:u} to estimate
the right\--hand side of Eq.~\eqref{e:trace:u}, thus arriving at
\eqref{ineq:ano_trace_pt:u} as desired.
%\null\hfill\qed
\qed
%%%%%%%%%%%%%%%%%
\par\vskip 10pt

%%%%%%%%%%%%%%%%%%%%%%%%%%%%%%%%%%%%%%%%%%%%%%%%%%%%%%%%%%%%%%%%%%%%%%%
%%%%%    A pointwise trace inequality (Lemma)    %%%%%%%%%%%%%%%%%%%%%%
%%%%%%%%%%%%%%%%%%%%%%%%%%%%%%%%%%%%%%%%%%%%%%%%%%%%%%%%%%%%%%%%%%%%%%%
\begin{lemma}\label{lem-Trace_pt}
{\rm (Pointwise trace inequalities.)}$\;$
Let\/ $\beta > 0$ and\/ $R > 0$, and set\/ $r = R / \sqrt{2}$.
Assume that\/
$u\in \tilde{H}^2(B^{+}_R; \mathfrak{w})$.
Then the following inequalities hold at almost every point\/
$x\in (-r,r)$, for every\/ $\xi\in (0,r)$:
\begin{equation}
\label{e:int_trace_pt:u}
\begin{aligned}
&   \frac{\beta}{2}\int_0^{\xi}
      |\nabla u(x,\xi')|^2\cdot (\xi')^{\beta - 1} \,\mathrm{d}\xi'
  - \frac{2}{\beta}\int_0^{\xi} 
    \left| \partial_{\xi}\nabla u(x,\xi')\right|^2\cdot
      (\xi')^{\beta + 1} \,\mathrm{d}\xi'
\\
& {}
  \leq
    \xi^{\beta}\cdot |\nabla u(x,\xi)|^2
  - \lim_{\xi'\to 0+}
    \left[ (\xi')^{\beta}\cdot |\nabla u(x,\xi')|^2 \right]
\\
& {}
  \leq
    \frac{3\beta}{2}\int_0^{\xi}
      |\nabla u(x,\xi')|^2\cdot (\xi')^{\beta - 1} \,\mathrm{d}\xi'
  + \frac{2}{\beta}\int_0^{\xi} 
    \left| \partial_{\xi}\nabla u(x,\xi')\right|^2\cdot
      (\xi')^{\beta + 1} \,\mathrm{d}\xi' \,.
\end{aligned}
\end{equation}

At almost every point\/ $x\in (-r,r)$,
all (Lebesgue) integrals above are finite and the limit\/
(viewed as {\rm a boundary condition})
\begin{equation}
\label{lim:trace_pt:v=0}
  L_0(x)\eqdef \lim_{\xi'\to 0+}
    \left[ (\xi')^{\beta}\cdot |\nabla u(x,\xi')|^2 \right]
  = 0 \quad\mbox{ exists. }
\end{equation}
\end{lemma}
%%%%%%%%%%%%%%%%%%%%%%%%%%%%%%%%%%%%%%%%%%%%%%%%%%%%%%%%%%%%%%%%%%%%%%%
%\par\vskip 10pt

%%%%%%%%%%%%%%%%%
\proof
Let us set
\begin{equation*}
%\label{lim_inf:trace_pt:v=0}
  \ell_0(x)\eqdef \liminf_{\xi'\to 0+}
    \left[ (\xi')^{\beta}\cdot |\nabla u(x,\xi')|^2 \right]
    \quad\mbox{ for every }\, x\in (-r,r) \,;
\end{equation*}
hence, $0\leq \ell_0(x)\leq \infty$.
Clearly, the function
\begin{math}
  \ell_0\colon x\mapsto \ell_0(x)\colon (-r,r)\to [0, +\infty]
\end{math}
is Lebesgue\--measurable.
Fatou's lemma yields
\begin{equation*}
%\label{lim_inf:Fatou:v=0}
  \int_{-r}^{r} \ell_0(x) \,\mathrm{d}x
  \leq \hat{\ell}_0\eqdef \liminf_{\xi'\to 0+}
    \left[ (\xi')^{\beta} \int_{-r}^{r}
           |\nabla u(x,\xi')|^2 \,\mathrm{d}x
    \right] \leq \infty \,.
\end{equation*}

Furthermore, from the hypothesis
$u\in \tilde{H}^2(B^{+}_R; \mathfrak{w})$ combined with
$Q^{+}_r\subset B^{+}_R$, we deduce that
\begin{equation}
\label{int:Fatou:v=0}
    \int_0^r
    \left| \partial_{\xi}\nabla u(x,\xi')\right|^2\cdot
      (\xi')^{\beta + 1} \,\mathrm{d}\xi'
  + \int_0^r
    |\nabla u(x,\xi')|^2\cdot (\xi')^{\beta} \,\mathrm{d}\xi'
  < \infty
\end{equation}
holds for every $x\in (-r,r)\setminus M_0$,
where $M_0\subset (-r,r)$ is a set of Lebesgue measure zero.
As an easy consequence, we observe that, due to the change of weight
$(\xi')^{\beta} \,\leftrightarrow\, (\xi')^{\beta - 1}$,
for every $x\in (-r,r)\setminus M_0$ we have
\begin{equation}
\label{infty:Fatou:v=0}
\begin{aligned}
&   \int_0^r |\nabla u(x,\xi')|^2\cdot
      (\xi')^{\beta - 1} \,\mathrm{d}\xi' = \infty
  \quad\mbox{ if and only if }\quad
\\
&   \int_0^{\xi} |\nabla u(x,\xi')|^2\cdot
      (\xi')^{\beta - 1} \,\mathrm{d}\xi' = \infty
  \quad\mbox{ holds for all $\xi\in (0,r]$. }
\end{aligned}
\end{equation}

Integrating the first inequality in \eqref{ineq:ano_trace_pt:u}
(in Lemma~\ref{lem-Trace_ineq} above),
for every $x\in (-r,r)\setminus M_0$ and every $\xi\in (0,r)$ we obtain
\begin{equation}
\label{e_left:int_trace_pt:u}
\begin{aligned}
&   \frac{\beta}{2}\int_0^{\xi}
      |\nabla u(x,\xi')|^2\cdot (\xi')^{\beta - 1} \,\mathrm{d}\xi'
  + \ell_0(x)
\\
& {}
  \leq
    \xi^{\beta}\cdot |\nabla u(x,\xi)|^2
  + \frac{2}{\beta}\int_0^{\xi} 
    \left| \partial_{\xi}\nabla u(x,\xi')\right|^2\cdot
      (\xi')^{\beta + 1} \,\mathrm{d}\xi'
\end{aligned}
\end{equation}
with the limit $0\leq \ell_0(x)\leq \infty$.
Thus, if the integral over $(0,r)$ in \eqref{infty:Fatou:v=0}
were infinite, so would be the integral over $(0,\xi)$
for every $0 < \xi\leq r$.
As the same integral over $(0,\xi)$ appears in
Ineq.~\eqref{e_left:int_trace_pt:u} as well, thanks to
\eqref{int:Fatou:v=0} this would force
$\xi^{\beta}\cdot |\nabla u(x,\xi)|^2 = \infty$
for every $0 < \xi\leq r$, thus contradicting
\eqref{int:Fatou:v=0}.
We conclude that, for every $x\in (-r,r)\setminus M_0$,
all integrals in \eqref{infty:Fatou:v=0} must be finite, whenever
$0 < \xi\leq r$.
Moreover, also $\ell_0(x) < \infty$ must hold.
However, if $\ell_0(x) > 0$ then all integrals in
\eqref{infty:Fatou:v=0} would have to be infinite,
another contradiction.
It follows that $\ell_0(x) = 0$.

Similarly, integrating both inequalities in
\eqref{ineq:ano_trace_pt:u}, combined with $\ell_0(x) = 0$,
for every $x\in (-r,r)\setminus M_0$ and every pair
$\xi_1, \xi_2\in \RR$, $0 < \xi_1 < \xi_2\leq r$, we get
\begin{equation}
\label{e_right:int_trace_pt:u}
\begin{aligned}
& \left|
    \xi_2^{\beta}\cdot |\nabla u(x,\xi_2)|^2
  - \xi_1^{\beta}\cdot |\nabla u(x,\xi_1)|^2
  \right|
\\
& {}
  \leq
    \frac{3\beta}{2}\int_0^{\xi}
      |\nabla u(x,\xi')|^2\cdot (\xi')^{\beta - 1} \,\mathrm{d}\xi'
  + \frac{2}{\beta}\int_0^{\xi}
    \left| \partial_{\xi}\nabla u(x,\xi')\right|^2\cdot
      (\xi')^{\beta + 1} \,\mathrm{d}\xi' \,.
\end{aligned}
\end{equation}
Consequently, for every $x\in (-r,r)\setminus M_0$, the function
\begin{equation*}
%\label{e_right:int_trace_pt:u}
  \xi \,\longmapsto\, \xi^{\beta}\cdot |\nabla u(x,\xi)|^2
  \colon (0,r]\to \RR_+
\end{equation*}
is absolutely continuous with the vanishing limit
$L_0(x) = \ell_0(x) = 0$ in Eq.~\eqref{lim:trace_pt:v=0}
(as $\xi\to 0+$).
In particular, the inequalities in \eqref{e:int_trace_pt:u}
are valid for every $x\in (-r,r)\setminus M_0$
and almost every $\xi\in (0,r)$,
with the function
\begin{math}
%\label{e_right:int_trace_pt:u}
  \xi \,\longmapsto\, \xi^{\beta}\cdot |\nabla u(x,\xi)|^2
  \colon (0,r]\to \RR_+
\end{math}
being absolutely continuous on $[0,r]$ with the limit $L_0(x) = 0$.
%\null\hfill\qed
\qed
%%%%%%%%%%%%%%%%%
\par\vskip 10pt

Finally, we integrate all equations and inequalities
\eqref{int:Fatou:v=0} -- \eqref{e_right:int_trace_pt:u}
with respect to $x\in (-r,r)$ to derive the following corollary of
Lemma~\ref{lem-Trace_pt}.

%%%%%%%%%%%%%%%%%%%%%%%%%%%%%%%%%%%%%%%%%%%%%%%%%%%%%%%%%%%%%%%%%%%%%%%
%%%%%    A pointwise trace inequality (Corollary)    %%%%%%%%%%%%%%%%%%
%%%%%%%%%%%%%%%%%%%%%%%%%%%%%%%%%%%%%%%%%%%%%%%%%%%%%%%%%%%%%%%%%%%%%%%
\begin{corollary}\label{cor-Trace_pt}
{\rm (Global trace inequalities.)}$\;$
Let\/ $\beta > 0$ and\/ $R > 0$, and set\/ $r = R / \sqrt{2}$.
Assume that\/
$u\in \tilde{H}^2(B^{+}_R; \mathfrak{w})$.
Then the following inequalities hold for every\/ $\xi\in (0,r)$:
\begin{equation}
\label{int:int_trace_pt:u}
\begin{aligned}
&   \frac{\beta}{2}\int_0^{\xi} (\xi')^{\beta - 1}
    \int_{-r}^{r}
      |\nabla u(x,\xi')|^2\cdot \,\mathrm{d}x \,\mathrm{d}\xi'
\\
& {}
  - \frac{2}{\beta}\int_0^{\xi} (\xi')^{\beta + 1}
    \int_{-r}^{r}
    \left| \partial_{\xi}\nabla u(x,\xi')\right|^2\cdot
      \,\mathrm{d}x \,\mathrm{d}\xi'
\\
& {}
  \leq
    \xi^{\beta} \int_{-r}^{r} |\nabla u(x,\xi)|^2 \,\mathrm{d}x
  - \lim_{\xi'\to 0+}
    \left[ (\xi')^{\beta}
      \int_{-r}^{r} |\nabla u(x,\xi')|^2 \,\mathrm{d}x
    \right]
\\
& {}
  \leq
    \frac{3\beta}{2}\int_0^{\xi} (\xi')^{\beta - 1}
    \int_{-r}^{r}
      |\nabla u(x,\xi')|^2\cdot \,\mathrm{d}x \,\mathrm{d}\xi'
\\
& {}
  + \frac{2}{\beta}\int_0^{\xi} (\xi')^{\beta + 1}
    \int_{-r}^{r}
    \left| \partial_{\xi}\nabla u(x,\xi')\right|^2\cdot
      \,\mathrm{d}x \,\mathrm{d}\xi' \,.
\end{aligned}
\end{equation}

All (Lebesgue) integrals above are finite and the limit\/
(viewed as {\rm a boundary condition})
\begin{equation}
\label{lim:trace:v=0}
  \hat{L}_0\eqdef \lim_{\xi'\to 0+}
    \left[ (\xi')^{\beta}
      \int_{-r}^{r} |\nabla u(x,\xi')|^2 \,\mathrm{d}x
    \right]
  = 0 \quad\mbox{ exists. }
\end{equation}
In addition, the restriction mapping \eqref{eq:restr-H^2} from
$\tilde{H}^2(B^{+}_R; \mathfrak{w})$ to
$H^2( B^{+}_{ R/\sqrt{2} }; \mathfrak{w} )$
is continuous.
\end{corollary}
%%%%%%%%%%%%%%%%%%%%%%%%%%%%%%%%%%%%%%%%%%%%%%%%%%%%%%%%%%%%%%%%%%%%%%%
%\par\vskip 10pt

%%%%%%%%%%%%%%%%%
\proof
The inequalities in \eqref{int:int_trace_pt:u}
follow directly from those in \eqref{e:int_trace_pt:u}.
Of course, Ineq.~\eqref{int:Fatou:v=0}
has to be replaced by
\begin{equation*}
%\label{int:Fatou:v=0}
    \iint_{Q^{+}_r}
    \left| \partial_{\xi}\nabla u(x,\xi')\right|^2\cdot
      (\xi')^{\beta + 1} \,\mathrm{d}x \,\mathrm{d}\xi'
  + \iint_{Q^{+}_r}
    |\nabla u(x,\xi')|^2\cdot (\xi')^{\beta}
    \,\mathrm{d}x \,\mathrm{d}\xi'
  < \infty \,.
\end{equation*}
The vanishing limit $\hat{L}_0 = 0$ in Eq.~\eqref{lim:trace:v=0}
(as $\xi\to 0+$)
is derived from \eqref{int:int_trace_pt:u}
in an analogous way as is $L_0(x)= 0$ from Eq.~\eqref{lim:trace_pt:v=0}
in our proof of Lemma~\ref{lem-Trace_pt} above.
Finally, we employ
the inequalities in \eqref{int:int_trace_pt:u}
to compare the norms on
$H^2( B^{+}_{ R/\sqrt{2} }; \mathfrak{w} )$ and
$\tilde{H}^2(B^{+}_R; \mathfrak{w})$,
defined by Eqs.\ \eqref{e:norm-H^2} and \eqref{eq:norm-H^2},
respectively.
Recall that
\begin{math}
  B^{+}_{ R/\sqrt{2} }\subset Q^{+}_{ R/\sqrt{2} }\subset B^{+}_R .
\end{math}
The continuity of the restriction mapping \eqref{eq:restr-H^2}
follows.
%\null\hfill\qed
\qed
%%%%%%%%%%%%%%%%%
\par\vskip 10pt

Our results in Corollary~\ref{cor-Trace_pt}
above will lead us to a restricted imbedding lemma,
Lemma~\ref{lem-Sobol-Morrey}, needed in
Section~\ref{s:smooth_Hoelder}, {\S}\ref{ss:smooth_Hoelder-2}.
This lemma will be derived from the following
Hardy\--Sobolev\--type inequality proved in
{\sc H.\ Castro} \cite[Theorem~4, p.~594]{Castro-2017}.

%%%%%%%%%%%%%%%%%%%%%%%%%%%%%%%%%%%%%%%%%%%%%%%%%%%%%%%%%%%%%%%%%%%%%%%
%%%%%    Castro: Hardy-Sobolev Imbedding (Lemma)    %%%%%%%%%%%%%%%%%%%
%%%%%%%%%%%%%%%%%%%%%%%%%%%%%%%%%%%%%%%%%%%%%%%%%%%%%%%%%%%%%%%%%%%%%%%
\begin{lemma}\label{lem-Hardy-Sobol}
{\rm (A Hardy\--Sobolev\--type inequality.)}$\;$
Let\/ $2\leq p < \infty$, $R > 0$, and set\/
$r = R / \sqrt{2}$.
Assume that\/
$a, a', b\in \RR$ satisfy the following inequalities,
\begin{itemize}
\item[{\rm (i)}]
$a > 0$, $0\leq a-b < 1$, $0 < a'< 1$, and\/
\item[{\rm (ii)}]
$1 - \frac{2}{p} < (a+a') - b\leq 1$.
\end{itemize}
Define $p^{\ast}\in (0,\infty)$ by\/
\begin{align}
\label{def:p^*}
  \frac{1}{p^{\ast}} + \frac{b+1}{2} =
  \frac{1}{p} + \frac{a+a'}{2} \,,
  \;\mbox{ whence }\;
  \frac{1}{p^{\ast}} =
  \frac{1}{p} - \frac{ 1 + b - (a + a')}{2} \leq \frac{1}{p} \,.
\end{align}
Then there exists a constant\/
$C\equiv C(R;p)\in (0,\infty)$ such that\/
\begin{equation}
\label{ineq:Castro_p}
\begin{aligned}
  \| u\|_{ L^{p^{*}}( B^{+}_r; \xi^{b p^{*}} ) }
  \eqdef
  \left( \int_{B^{+}_r} |u(x,\xi)|^{p^{*}}
         \cdot \xi^{b p^{*}}\cdot \,\mathrm{d}x \,\mathrm{d}\xi
  \right)^{1/p^{*}}
  \leq C\cdot
  \| u\|_{ W^{1,p}( B^{+}_R; \xi^{ap} ) }
\end{aligned}
\end{equation}
holds for all\/
$u\in W^{1,p}( B^{+}_R; \xi^{ap} )$, i.e., for all\/
$u\in W^{1,p}_{\mathrm{loc}}(B^{+}_R)$ with the norm
\begin{equation*}
%\label{norm:W^1,p}
  \| u\|_{ W^{1,p}( B^{+}_R; \xi^{ap} ) } \eqdef
  \left(
    \int_{B^{+}_R}
    \left( |\nabla u(x,\xi)|^p + |u(x,\xi)|^p\right)
    \cdot \xi^{ap}\cdot \,\mathrm{d}x \,\mathrm{d}\xi
  \right)^{1/p} < \infty \,.
\end{equation*}

In particular, for\/ $p = 2$ and\/
$p^{*} = 2^{*} = \frac{2}{(a + a') - b}$ $({}\geq 2)$
the following analogue of the restricted imbedding \eqref{eq:restr-H^2}
is continuous, this time considered as a linear mapping
\begin{equation*}
%\label{eq:restr-W^1,p->L^p}
  u \,\longmapsto\, u\vert_{B^{+}_r}
  \colon W^{1,2}( B^{+}_R; \xi^{2a} ) \,\longrightarrow\,
         L^{2^{*}}( B^{+}_r; \xi^{2^{*} b} ) \,.
\end{equation*}
\end{lemma}
%%%%%%%%%%%%%%%%%%%%%%%%%%%%%%%%%%%%%%%%%%%%%%%%%%%%%%%%%%%%%%%%%%%%%%%
\par\vskip 10pt

%\par\vskip 10pt
%%%%%%%%%%%%%%%%%
\proof
This lemma follows directly from
Theorem~4 in {\sc H.\ Castro} \cite[p.~594]{Castro-2017}.
Since Theorem~4 in \cite{Castro-2017}
is formulated for a $C^1$ function $u\colon \RR^2\to \RR$
with compact support, we have to apply it to the function
$\phi u$, where $\phi\colon \RR^2\to \RR$ is a $C^{\infty}$ function
with the following properties:
$\phi(x,\xi) = 1$ in $B_r(x_0,0)$,
$\phi(x,\xi) = 0$ in $\RR^2\setminus B_R(x_0,0)$, and
$0\leq \phi(x,\xi)\leq 1$ in $B_R(x_0,0)\setminus B_r(x_0,0)$.
Recall that $B_R\equiv B_R(x_0,\xi_0)$ denotes the open disc in $\RR^2$
with radius $R > 0$ centered at the point $(x_0,\xi_0)\in \RR^2$.
We have abbreviated the upper half\--disc by
$B^{+}_R\equiv B^{+}_R(x_0,\xi_0)$.
Applying \cite[Theorem~4, p.~594]{Castro-2017}
with the compactly support product function
$\phi u\in W^{1,p}( B^{+}_R; |x|^{a'}\xi^{ap} )$,
we obtain the following Hardy\--Sobolev\--type inequality,
\begin{equation}
\label{grad:Castro_p}
\begin{aligned}
& \| (\phi u)\|_{ L^{p^{*}}( B^{+}_r; \xi^{b p^{*}} ) } =
  \left( \int_{B^{+}_r} |\phi u|^{p^{*}}
         \cdot \xi^{b p^{*}}\cdot \,\mathrm{d}x \,\mathrm{d}\xi
  \right)^{1/p^{*}}
\\
& {}
  \leq C'\cdot
  \| (\phi u)\|_{ W^{1,p}( B^{+}_R; |x|^{a'}\xi^{ap} ) }
  = C'\cdot
  \left(
    \int_{B^{+}_R} |\nabla (\phi u)|^p
    \cdot |x|^{a'}\xi^{ap}\cdot \,\mathrm{d}x \,\mathrm{d}\xi
  \right)^{1/p}
\end{aligned}
\end{equation}
for all\/
$u\in W^{1,p}_{\mathrm{loc}}(B^{+}_R)$ with\/
\begin{math}
  \| (\phi u)\|_{ W^{1,p}( B^{+}_R; |x|^{a'}\xi^{ap} ) } < \infty .
\end{math}
Here, $C'\equiv C'(R;p)\in (0,\infty)$
is a constant independent from the product function $\phi u$.
Thanks to
$\nabla (\phi u) = \phi\, \nabla u + u\, (\nabla\phi)$ with both
$\phi ,\, |\nabla\phi| \in L^{\infty}(\RR^2)$
and $|x|^{a'}\leq R^{a'} < \infty$ for $x\in (-R,R)$, as well,
we can apply the triangle inequality in $L^p(B^{+}_R)$
to the right\--hand side of Ineq.~\eqref{grad:Castro_p}
in order to derive the desired inequality \eqref{ineq:Castro_p}.
%\null\hfill\qed
\qed
%%%%%%%%%%%%%%%%%
\par\vskip 10pt

Unfortunately, earlier results of this kind
({\sc P.~M.~N.\ Feehan} and {\sc C.~A.\ Pop}
 \cite{Feehan-Pop-17}, Lemma 2.2, Eq.\ (2.2), on p.~1091, and
 {\sc H.\ Koch} \cite[Lemma 4.2.4, p.~62]{Koch-1999})
seem to be useless in our case due to the hypothesis
$u\in H^1(B^{+}_R; \mathfrak{w})$ that is weaker than
$u\in W^{1,2}( B^{+}_R; \xi^{\beta - 1} )$
owing to the seminorm
\begin{equation*}
%\label{loc:norm-H^1}
\textstyle
    \|\nabla u\|_{ L^2( B^{+}_R; \xi^{\beta} ) } \eqdef
  \left( \int_{B^{+}_R}
    \left( |u_{x}|^2 + |u_{\xi}|^2\right)
    \cdot \xi^{\beta}\cdot \,\mathrm{d}x \,\mathrm{d}\xi
  \right)^{1/2}
\end{equation*}
in Eqs.\ \eqref{loc:norm-H^1} and \eqref{eq:norm-H^2}
being weaker than
\begin{equation*}
%\label{loc:norm-H^1}
\textstyle
    \|\nabla u\|_{ L^2( B^{+}_R; \xi^{\beta - 1} ) } \eqdef
  \left( \int_{B^{+}_R}
    \left( |u_{x}|^2 + |u_{\xi}|^2\right)
    \cdot \xi^{\beta - 1}\cdot \,\mathrm{d}x \,\mathrm{d}\xi
  \right)^{1/2}
\end{equation*}
which appears in Eq.~\eqref{e:norm-H^2}, thanks to
$\xi^{\beta} / \xi^{\beta - 1} = \xi$.

In contrast to these results, the next lemma enables us
to establish the restricted Sobolev imbedding \eqref{e:H^2->L^p}
(see {\S}\ref{ss:smooth_Hoelder-2})
by replacing the pair $(p,p^{*})$ by $(2,p) = (2,2^{*})$.
In this pair we allow for $p = 2^{*}\in (2,\infty)$ arbitrary
to which we associate suitable constants
$a, a', b\in (0,\infty)$ that verify conditions
{\rm (i)}, {\rm (ii)}, and Eq.~\eqref{def:p^*}
of Lemma~\ref{lem-Hardy-Sobol}.

%%%%%%%%%%%%%%%%%%%%%%%%%%%%%%%%%%%%%%%%%%%%%%%%%%%%%%%%%%%%%%%%%%%%%%%
%%%%%    Sobolev-Morrey Imbedding (Lemma)    %%%%%%%%%%%%%%%%%%%%%%%%%%
%%%%%%%%%%%%%%%%%%%%%%%%%%%%%%%%%%%%%%%%%%%%%%%%%%%%%%%%%%%%%%%%%%%%%%%
\begin{lemma}\label{lem-Sobol-Morrey}
{\rm (Two Sobolev\--type imbeddings.)}$\;$
Let\/ $2 < p < \infty$ and\/ $R > 0$ be arbitrary, and set\/
$r = R / \sqrt{2}$.
Let\/ $\beta\in \RR$ satisfy\/
\begin{equation}
\label{cond:beta}
  0 < \beta - 1 < \frac{4}{p-2} \,.
\end{equation}
Then there exists a constant\/
$C_{\beta}\equiv C_{\beta}(R;p)\in (0,\infty)$ such that\/
\begin{equation}
\label{ineq:W^1,2->Castro_p}
\begin{aligned}
& \| u\|_{ L^p( B^{+}_r; \xi^{\beta - 1} ) }
  =
  \left( \int_{B^{+}_r} |u(x,\xi)|^p
         \cdot \xi^{\beta - 1}\cdot \,\mathrm{d}x \,\mathrm{d}\xi
  \right)^{1/p}
\\
& {}
  \leq C_{\beta}\cdot
  \| u\|_{ W^{1,2}( B^{+}_R; \xi^{\beta - 1} ) }
  = C_{\beta}\cdot
  \left(
    \int_{B^{+}_R} \left( |\nabla u|^2 + |u|^2\right)
    \cdot \xi^{\beta - 1}\cdot \,\mathrm{d}x \,\mathrm{d}\xi
  \right)^{1/2}
\end{aligned}
\end{equation}
holds for all\/
$u\in W^{1,2}( B^{+}_R; \xi^{\beta - 1} )$.

Furthermore, there exists another constant\/
$C_{\beta}^{\prime}\equiv C_{\beta}^{\prime}(R;p)\in (0,\infty)$
such that\/
\begin{equation}
\label{ineq:H^2->Castro_p}
  \| u\|_{ L^p( B^{+}_{R/2}; \xi^{\beta - 1} ) }
  \leq C_{\beta}^{\prime}\cdot
  \| u\|_{ H^2( B^{+}_R; \mathfrak{w} ) }
    \quad\mbox{ for all }\,
    u\in H^2( B^{+}_R; \mathfrak{w} ) \,.
\end{equation}
In particular, the {\em restricted Sobolev imbedding\/}
(cf.\ {\rm Eq.}~\eqref{e:H^2->L^p})
\begin{equation}
\label{eq:H^2->L^p}
  u\vert_{ B_R^{+} }\longmapsto
  u\vert_{ B_{R/2}^{+} }\colon
  H^2\left( B_R^{+}; \mathfrak{w}\right) \hookrightarrow
  L^p\left( B_{R/2}^{+}; \mathfrak{w}\right)
\end{equation}
is continuous.
\end{lemma}
%%%%%%%%%%%%%%%%%%%%%%%%%%%%%%%%%%%%%%%%%%%%%%%%%%%%%%%%%%%%%%%%%%%%%%%
%\par\vskip 10pt

%%%%%%%%%%%%%%%%%
\proof
We wish to apply Lemma~\ref{lem-Hardy-Sobol} stated above
with the weight $\xi^{\beta - 1}$ as indicated in \eqref{e:H^2->L^p}.
We replace the pair $(p,p^{*})$ by $(2,2^{*})$ and forget
the former one entirely; thus, from now on, we may write
$p = 2^{*}$ with $2 < p < \infty$.
We need to fix the constant $\beta\in (1,\infty)$ in such a way that
Lemma~\ref{lem-Hardy-Sobol} is applicable with suitable constants
$a, a', b\in \RR$.
Consequently, we choose the constants
\begin{equation*}
  a = \frac{\beta - 1}{2} \quad\mbox{ and }\quad
  b = \frac{\beta - 1}{2^{*}} = \frac{\beta - 1}{p} \,.
\end{equation*}
Clearly, we have $a > 0$, $b > 0$, and $a-b > 0$, by $p > 2$.
In order to fulfill also the condition $a-b < 1$,
we have to choose $\beta\in (1,\infty)$ such that
\begin{math}
  a-b = (\beta - 1)\left( \frac{1}{2} - \frac{1}{p}\right) < 1
\end{math}
or, equivalently,
\begin{math}
  1 < \beta < 1 + \frac{2p}{p-2} \,.
\end{math}
These inequalities follow from our choice of $p > 2$ and $\beta$
obeying the conditions in \eqref{cond:beta}.
We have no other restriction on $p\in (2,\infty)$.
The remaining constant, $a'\in (0,1)$, must be chosen in such a way that
Eq.~\eqref{def:p^*} holds with the pair
$(2,2^{*})$ in place of $(p,p^{*})$, i.e.,
\begin{math}
  a'= \frac{2}{2^{*}} - (a-b) ,
\end{math}
together with the inequalities $a'> 0$ and
$(a+a') - b\leq 1$.
Since $2^{*} = p > 2$, we get
\begin{math}
  a'= \frac{2}{2^{*}} - (a-b) < 1 - (a-b) ,
\end{math}
i.e., $a+a'< b+1$.
It remains to verify $a'> 0$ which is equivalent with
(from now on we write $2^{*} = p > 2$)
\begin{equation*}
  \frac{2}{p} > (\beta - 1)\left( \frac{1}{2} - \frac{1}{p}\right) \,.
\end{equation*}
This inequality is equivalent with the condition in \eqref{cond:beta}.

The desired inequality in \eqref{ineq:W^1,2->Castro_p}
now follows directly from Lemma~\ref{lem-Hardy-Sobol}.
Finally, we apply
Corollary~\ref{cor-Trace_pt}, Eq.~\eqref{int:int_trace_pt:u}
with $R$ and $r = R / \sqrt{2}$,
to the right\--hand side of Ineq.~\eqref{ineq:W^1,2->Castro_p}
with $r$ and $R/2 = r / \sqrt{2}$
(in place of $R$ and $r = R / \sqrt{2}$, respectively)
to derive \eqref{ineq:H^2->Castro_p}.
%
%\null\hfill\qed
\qed
%%%%%%%%%%%%%%%%%
\par\vskip 10pt

%%%%%%%%%%%%%%%%%%%%%%%%%%%%%%%%%%%%%%%%%%%%%%%%%%%%%%%%%%%%%%%%%%%%%%%
%%%%%    Some known elliptic regularity results    %%%%%%%%%%%%%%%%%%%%
%%%%%%%%%%%%%%%%%%%%%%%%%%%%%%%%%%%%%%%%%%%%%%%%%%%%%%%%%%%%%%%%%%%%%%%

\section{Appendix: Some known elliptic regularity results}
\label{s:ellipt_regul}

In this appendix we collect a few known results on
the local regularity of a weak solution $u\in V$
to the degenerate elliptic problem
$(\lambda I + \mathcal{A}) u = f\in H$, i.e., for
$u = (\lambda I + \mathcal{A})^{-1} f\in V$.
Recall that $\lambda > \lambda_0$ with the constant $\lambda_0 > 0$
determined by Ineq.~\eqref{norm:e^(-tA)}.

The first regularity result is due to
{\sc P.~M.~N.\ Feehan} and {\sc C.~A.\ Pop}
\cite{Feehan-Pop-15}, Theorem 3.16, Eq.\ (3.12), on p.~385.

%%%%%%%%%%%%%%%%%%%%%%%%%%%%%%%%%%%%%%%%%%%%%%%%%%%%%%%%%%%%%%%%%%%%%%%
%%%%%    Smoothing properties of Heston's semigroup (Lemma)    %%%%%%%%
%%%%%%%%%%%%%%%%%%%%%%%%%%%%%%%%%%%%%%%%%%%%%%%%%%%%%%%%%%%%%%%%%%%%%%%
\begin{lemma}\label{lem-smooth_Hoelder-1}
{\rm ($H^2$-smoothing property.)}$\;$
Let\/
$\rho$, $\sigma$, $\theta$, $q_r$, and\/ $\gamma$
be given constants in $\RR$,
$\rho\in (-1,1)$, $\sigma > 0$, $\theta > 0$, and\/ $\gamma > 0$.
Assume that\/ $\beta$, $\gamma$, $\kappa$, and\/ $\mu$
are chosen as specified in {\rm Proposition~\ref{prop-Lions}}
and\/ $\lambda > \lambda_0$.
Then, given any\/ $x_0\in \RR$ and\/
$R_0, R_1\in \RR$ with $0 < R_1 < R_0$, and any function $f\in H$,
the restriction
$u\vert_{ B^{+}_{R_1}(x_0,0) }$ of the function
$u = (\lambda I + \mathcal{A})^{-1} f\in V$
to the open half\--disc $B^{+}_{R_1}(x_0,0)$ satisfies\/
\begin{math}
  u\vert_{ B^{+}_{R_1}(x_0,0) }\in 
  H^2\left( B_{R_1}^{+}(x_0,0); \mathfrak{w}\right) .
\end{math}
Furthermore, there is a constant\/ $C_1\in (0,\infty)$
independent from $f$ and\/ $u$, such that
\begin{equation}
\label{est:u_1,k(t)}
\begin{aligned}
  \| u\|_{ H^2\left( B_{R_1}^{+}(x_0,0); \mathfrak{w}\right) }
  \leq C_1
  \left(
    \| u\|_{ L^2\left( B_{R_0}^{+}(x_0,0); \mathfrak{w}\right) }
  + \| f\|_{ L^2\left( B_{R_0}^{+}(x_0,0); \mathfrak{w}\right) }
  \right) \,.
\end{aligned}
\end{equation}
(The weighted Sobolev norm on the left\--hand side has been introduced
 in {\rm Eq.}~\eqref{e:norm-H^2}.)
\end{lemma}
%%%%%%%%%%%%%%%%%%%%%%%%%%%%%%%%%%%%%%%%%%%%%%%%%%%%%%%%%%%%%%%%%%%%%%%
\par\vskip 10pt

This lemma gets us from $H = L^2(\HH;\mathfrak{w})$
into the (local) interior regularity of
weighted $H^2$-type over an open half\--disc $B_{R_1}^{+}(x_0,0)$.
Our restricted Hardy\--Sobolev\--type imbedding
(Lemma \ref{lem-Sobol-Morrey})
brings an $H^2$-type function into a weighted $L^p$-space over
a smaller open half\--disc $B_{R_1'}^{+}(x_0,0)$ with the radius
$R_1'= R_1/2$ ($0 < R_1'< R_1 < R_0)$.
(This step will require an additional upper bound on $\beta > 1$,
 in addition to Ineq.~\eqref{ineq:beta-1<mu},
 in order to allow for $p > 2$ large enough in
 Lemma~\ref{lem-smooth_Hoelder-2} below
 and still fulfill Ineq.~\eqref{cond:beta} in
 Lemma~\ref{lem-Sobol-Morrey}.)

Now we continue with another local regularity result for
a weak solution
$u = (\lambda I + \mathcal{A})^{-1} f$ $\in V$,
this time for $u\in H$ with $f\in V$ satisfying also
\begin{math}
  f\vert_{ B^{+}_{R_1}(x_0,0) } \in
  L^p\left( B_{R_1}^{+}(x_0,0); \mathfrak{w}\right) .
\end{math}
(This weighted $L^p$-space has been introduced
 in Eq.~\eqref{e:norm_L^p}.)

%%%%%%%%%%%%%%%%%%%%%%%%%%%%%%%%%%%%%%%%%%%%%%%%%%%%%%%%%%%%%%%%%%%%%%%
%%%%%    Smoothing properties of Heston's semigroup (Lemma)    %%%%%%%%
%%%%%%%%%%%%%%%%%%%%%%%%%%%%%%%%%%%%%%%%%%%%%%%%%%%%%%%%%%%%%%%%%%%%%%%
\begin{lemma}\label{lem-smooth_Hoelder-2}
{\rm ($C_s^{\alpha}$-smoothing property.)}$\;$
Let\/
$\rho$, $\sigma$, $\theta$, $q_r$, and\/ $\gamma$
be given constants in $\RR$,
$\rho\in (-1,1)$, $\sigma > 0$, $\theta > 0$, and\/ $\gamma > 0$.
Assume that\/ $\beta$, $\gamma$, $\kappa$, and\/ $\mu$
are chosen as specified in {\rm Proposition~\ref{prop-Lions}}
and\/ $\lambda > \lambda_0$.
Finally, let\/ $p$ satisfy\/
$\max\{ 4,\, 2+\beta\} < p < \infty$.
Then, given any\/ $x_0\in \RR$ and\/ $R_1'\in (0,\infty)$, 
there are constants\/
$R_2\equiv R_2(R_1')$, which depends on $R_1'$, $\alpha\in (0,1)$, and\/
$C_2\in (0,\infty)$ with the following properties:

{\rm (a)}$\;$
$0 < R_2 < R_1'$,

{\rm (b)}$\;$
given any function $f\in H$ with the restriction
\begin{math}
  f\vert_{ B^{+}_{R_1'}(x_0,0) } \in
  L^p\left( B_{R_1'}^{+}(x_0,0); \mathfrak{w}\right) ,
\end{math}
the restriction
$u\vert_{ \overline{B}^{+}_{R_2}(x_0,0) }$ of the function
$u = (\lambda I + \mathcal{A})^{-1} f\in V$
to the closed half\--disc
$\overline{B}^{+}_{R_2}(x_0,0)$ satisfies\/
\begin{math}
  u\vert_{ \overline{B}^{+}_{R_2}(x_0,0) }\in 
  C_s^{\alpha}
  \left( \overline{B}_{R_2}^{+}(x_0,0); \mathfrak{w}\right) ,
\end{math}
and\/

{\rm (c)}$\;$
for all pairs $f$ and\/ $u$ from {\rm Part (b)},
the following inequality holds,
\begin{equation}
\label{est:u_2,k(t)}
\begin{aligned}
  \| u\|_{ C_s^{\alpha}\left( \overline{B}_{R_2}^{+}(x_0,0)\right) }
  \leq C_2
  \left(
    \| u\|_{ L^2\left( B_{R_1'}^{+}(x_0,0); \mathfrak{w}\right) }
  + \| f\|_{ L^p\left( B_{R_1'}^{+}(x_0,0); \mathfrak{w}\right) }
  \right) \,.
\end{aligned}
\end{equation}
(The weighted H\"older norm on the left\--hand side has been introduced
 in {\rm Eq.}~\eqref{norm:Hoelder}.)
\end{lemma}
%%%%%%%%%%%%%%%%%%%%%%%%%%%%%%%%%%%%%%%%%%%%%%%%%%%%%%%%%%%%%%%%%%%%%%%
\par\vskip 10pt

This lemma improves the (local) interior regularity of $u\in H$ from
\begin{math}
  L^p\left( B_{R_1'}^{+}(x_0,0); \mathfrak{w}\right)
\end{math}
to the weighted H\"older space
\begin{math}
  C_s^{\alpha}\left( \overline{B}_{R_2}^{+}(x_0,0)\right) ,
\end{math}
$0 < R_2 < R_1'$ $({} < R_1 < R_0)$.
The proof of this lemma is given in
{\sc P.~M.~N.\ Feehan} and {\sc C.~A.\ Pop}
\cite{Feehan-Pop-17}, Theorem 1.11, Eq.\ (1.31), on p.~1083;
see also
\cite{Feehan-Pop-15}, Theorem 2.5, Eq.\ (2.12), pp.\ 375--376.
We stress that the constant $R_2\equiv R_2(R_1')$ depends on $R_1'$,
while $R_1'\in (0,\infty)$ is arbitrary.

The last (local) interior regularity results for $u\in H$
brings $u$ from
\begin{math}
  C_s^{\alpha}\left( \overline{B}_{R_2}^{+}(x_0,0)\right) ,
\end{math}
to another weighted H\"older space
\begin{math}
  C_s^{2+\alpha}\left( \overline{B}_{R_2'}^{+}(x_0,0)\right) ,
\end{math}
$0 < R_2'< R_2$ $({} < R_1'< R_1 < R_0)$.
(The weighted H\"older space above has been introduced
 in Eq.~\eqref{norm:2+alpha}.)
Here, the constants $R_2$ and $R_2'$ are arbitrary with
$0 < R_2'< R_2 < \infty$.

%%%%%%%%%%%%%%%%%%%%%%%%%%%%%%%%%%%%%%%%%%%%%%%%%%%%%%%%%%%%%%%%%%%%%%%
%%%%%    Smoothing properties of Heston's semigroup (Lemma)    %%%%%%%%
%%%%%%%%%%%%%%%%%%%%%%%%%%%%%%%%%%%%%%%%%%%%%%%%%%%%%%%%%%%%%%%%%%%%%%%
\begin{lemma}\label{lem-smooth_Hoelder-3}
{\rm ($C_s^{2+\alpha}$-smoothing property.)}$\;$
Let\/
$\rho$, $\sigma$, $\theta$, $q_r$, and\/ $\gamma$
be given constants in $\RR$,
$\rho\in (-1,1)$, $\sigma > 0$, $\theta > 0$, and\/ $\gamma > 0$.
Assume that\/ $\beta$, $\gamma$, $\kappa$, and\/ $\mu$
are chosen as specified in {\rm Proposition~\ref{prop-Lions}}
and\/ $\lambda > \lambda_0$.
Finally, let\/ $\alpha\in (0,1)$ be arbitrary.
Then, given any\/ $x_0\in \RR$ and\/
$R_2, R_2'\in \RR$ with $0 < R_2'< R_2$, and any function $f\in H$
with the restriction
\begin{math}
  f\vert_{ \overline{B}^{+}_{R_2}(x_0,0) } \in
  C_s^{\alpha}\left( \overline{B}_{R_2}^{+}(x_0,0)\right) ,
\end{math}
the restriction
$u\vert_{ \overline{B}^{+}_{R_2'}(x_0,0) }$ of the function
$u = (\lambda I + \mathcal{A})^{-1} f\in V$
to the closed half\--disc
$\overline{B}^{+}_{R_2'}(x_0,0)$ satisfies\/
\begin{math}
  u\vert_{ \overline{B}^{+}_{R_2'}(x_0,0) }\in 
  C_s^{2+\alpha}\left( \overline{B}_{R_2'}^{+}(x_0,0)\right) .
\end{math}
Furthermore, there is a constant\/ $C_3\in (0,\infty)$
independent from $f$ and\/ $u$, such that
\begin{equation}
\label{est:u_3,k(t)}
\begin{aligned}
&   \| u\|_{ C_s^{2+\alpha}\left( \overline{B}_{R_2'}^{+}(x_0,0)\right) }
  \leq C_3
  \left(
    \| f\|_{ C_s^{2+\alpha}\left( \overline{B}_{R_2}^{+}(x_0,0)\right) }
  + \| u\|_{ C\left( \overline{B}_{R_2}^{+}(x_0,0)\right) }
  \right) \,.
\end{aligned}
\end{equation}
\end{lemma}
%%%%%%%%%%%%%%%%%%%%%%%%%%%%%%%%%%%%%%%%%%%%%%%%%%%%%%%%%%%%%%%%%%%%%%%
\par\vskip 10pt

This lemma is proved in
{\sc P.~M.~N.\ Feehan} and {\sc C.~A.\ Pop}
\cite{Feehan-Pop-14}, Theorem 8.1, Eq.\ (8.4), pp.\ 937--938
(see also
 {\sc P.~M.~N.\ Feehan}
 \cite{Feehan-17}, Theorem 1.1, Part~2, on pp.\ 2487--2488).

Lemma~\ref{lem-smooth_Hoelder-3}
has the following important consequence for the boundary limits
(as $\xi\to 0+$) of the functions
$\xi\cdot f_{xx}$, $\xi\cdot f_{x\xi}$, and $\xi\cdot f_{\xi\xi}$,
provided
$f\in C_s^{2+\alpha}\left( \overline{B}_{R_2'}^{+}(x_0,0)\right)$.

%%%%%%%%%%%%%%%%%%%%%%%%%%%%%%%%%%%%%%%%%%%%%%%%%%%%%%%%%%%%%%%%%%%%%%%
%%%%%    Boundary limits of $f\in C_s^{2+\alpha} (Corollary)    %%%%%%%
%%%%%%%%%%%%%%%%%%%%%%%%%%%%%%%%%%%%%%%%%%%%%%%%%%%%%%%%%%%%%%%%%%%%%%%
\begin{corollary}\label{cor-boundary_lim}
{\rm (Boundary limits in $C_s^{2+\alpha}$.)}$\;$
Let\/ $\alpha\in (0,1)$ be arbitrary,
$x_0\in \RR$ and\/ $R\in (0,\infty)$.
Then every function
$f\in C_s^{2+\alpha}\left( \overline{B}_R^{+}(x_0,0)\right)$
has the following behavior near the boundary\/
$\partial\HH = \RR\times \{ 0\}$ of the half\--plane
$\HH = \RR\times (0,\infty)\subset \RR^2$:
\begin{equation}
\label{lim_xi:2+alpha}
  \lim_{\xi\to 0+}
  \left[ \xi\cdot
  \left( |f_{xx}(x,\xi)| + |f_{x\xi}(x,\xi)| + |f_{\xi\xi}(x,\xi)|
  \right)
  \right] = 0
\end{equation}
for every $x\in (x_0 - R,\, x_0 + R)$.
In addition, there exists a constant\/
$c_{\alpha}\in (0,\infty)$ such that\/
\begin{equation}
\label{e:lim_xi:2+alpha}
  |f_{xx}(x,\xi)| + |f_{x\xi}(x,\xi)| + |f_{\xi\xi}(x,\xi)|
  \leq c_{\alpha}\,
    \| f\|_{ C_s^{2+\alpha}\left( \overline{B}_R^{+}(x_0,0)\right) }
    \cdot \xi^{ - [1 - (\alpha / 2)] }
\end{equation}
for all\/ $(x,\xi)\in B_R^{+}(x_0,0)$.
\end{corollary}
%%%%%%%%%%%%%%%%%%%%%%%%%%%%%%%%%%%%%%%%%%%%%%%%%%%%%%%%%%%%%%%%%%%%%%%
\par\vskip 10pt

%\par\vskip 10pt
%%%%%%%%%%%%%%%%%
\proof
In a somewhat stronger version stated in Eq.~\eqref{lim:2+alpha},
the limit \eqref{lim_xi:2+alpha} in this corollary is proved in
{\sc P.~M.~N.\ Feehan} and {\sc C.~A.\ Pop}
\cite{Feehan-Pop-13}, Lemma 3.1, Eq.\ (3.1), on p.~4409
(see also
 {\sc P.\ Daskalopoulos} and {\sc R.\ Hamilton}
 \cite{Daska-Hamil-98}, Prop.\ I.12.1 on p.~940).
The estimate \eqref{e:lim_xi:2+alpha}
is derived from this limit combined with
$f\in C_s^{2+\alpha}\left( \overline{B}_R^{+}(x_0,0)\right)$; cf.\
Eq.~\eqref{norm:2+alpha} with all
\begin{math}
%\label{norm:2+alpha}
  \xi\cdot f_{xx} ,\, \xi\cdot f_{x\xi} ,\, \xi\cdot f_{\xi\xi} \in
  C_s^{\alpha}\left( \overline{B}_R^{+}(x_0,0)\right) .
\end{math}

As far as the H\"older norm in
$C_s^{\alpha}\left( \overline{B}_R^{+}(x_0,0)\right)$,
given by Eq.~\eqref{norm:Hoelder}, is concerned, notice that for
every pair of points
$P^{\ast} = (x,0)\in \partial\HH$ and $P = (x,\xi)\in \HH$, with
$x\in (x_0 - R, x_0 + R)$ and
$(x,\xi)\in B_R^{+}(x_0,0)$, we have the $s$-distance
\begin{math}
%\label{def:cyclo}
  s(P,P^{\ast}) = \sqrt{\xi / 2}
\end{math}
(see Eq.~\eqref{def:cyclo}).
Hence, the inequality in \eqref{e:lim_xi:2+alpha}
follows.
%\null\hfill\qed
\qed
%%%%%%%%%%%%%%%%%
\par\vskip 10pt

%%%%%%%%%%%%%%%%%%%%%%%%%%%%%%%%%%%%%%%%%%%%%%%%%%%%%%%%%%%%%%%%%%%%%%%

\section*{Acknowledgment.}
A part of this research was performed while
the second author (P.T.) was a visiting professor at
Toulouse School of Economics, I.M.T.,
Universit\'e de Toulouse -- Capitole, Toulouse, France.

%The authors express their thanks to
%an anonymous referee for
%very careful reading of the manuscript and
%a number of helpful suggestions.

%\newpage
%%%%%%%%%%%%%%%%%%%%%%%%%%%%%%%%%%%%%%%%%%%%%%%%%%%%%%%%%%%%%%%%%%%%%%%
%%%%%    BIBLIOGRAPHY    %%%%%%%%%%%%%%%%%%%%%%%%%%%%%%%%%%%%%%%%%%%%%%
%%%%%%%%%%%%%%%%%%%%%%%%%%%%%%%%%%%%%%%%%%%%%%%%%%%%%%%%%%%%%%%%%%%%%%%

%%%%%%%%%%%%%%%%%%%%%%%%%%%%%%%%%%%%%%%%%%%%%%%%%%%%%%%%%%%%%%%%%%%%%%%%%%%%
%               AMS-LaTeX Version 1.1 file for electronic submission       %
% PAPER.BBL                  			          October 1999     %
%                                                                          %
%       Ground\--State Positivity, Negativity, and Compactness		   %
%       for a Schr\"odinger Operator in $\mathbb{R}^N$			   %
%                                                                          %
%%%%%%%%%%%%%%%%%%%%%%%%%%%%%%%%%%%%%%%%%%%%%%%%%%%%%%%%%%%%%%%%%%%%%%%%%%%%
%
%\bibliographystyle{amsplain}        %% if within the text
%
\makeatletter \renewcommand{\@biblabel}[1]{\hfill#1.} \makeatother
\ifx\undefined\bysame
\newcommand{\bysame}{\leavevmode\hbox to3em{\hrulefill}\,}
\fi
%

%
%%%%%%%%%%%%%%%%%%%%%%%%%%%%%%%%%%%%%%%%%%%%%%%%%%%%%%%%%%%%%%%%%%%%%%%
%%%%%%%%%%%%%%%%%%%%%%%%%%%%%%%%%%%%%%%%%%%%%%%%%%%%%%%%%%%%%%%%%%%%%%%
%%%%%%%%%%%%%%%%%%%%%%%%%%%%%%%%%%%%%%%%%%%%%%%%%%%%%%%%%%%%%%%%%%%%%%%
%
%%%%%%%%%%%%%%%%%%%%%%%%%%%
%
\end{document}